# Groupes $p$-divisibles, groupes finis et modules filtrés

Par Christophe Breuil

*A Luce et Pascale Breuil*

Sommaire

# 1. Introduction

Fixons $p$ un nombre premier, $k$ un corps parfait de caractéristique $p > 0$, $W = W(k)$ l'anneau des vecteurs de Witt à coefficients dans $k$, $K_0 = \mathrm{Frac}(W)$, $K$ une extension finie totalement ramifiée de $K_0$ de degré $e \geq 1$, $\mathcal{O}_K$ l'anneau des entiers de $K$, $\overline{K}$ une clôture algébrique de $K$, $\mathcal{O}_{\overline{K}}$ l'anneau des entiers de $\overline{K}$ et $\pi$ une uniformisante de $\mathcal{O}_K$. Tous les schémas en groupes considérés sont commutatifs. On appelle parfois $p$-groupe un schéma en groupes commutatifs annulés par une puissance de $p$.

Le but de cet article est de donner une classification, pour $p \neq 2$ mais sans restriction sur la ramification $e$, des groupes $p$-divisibles sur $\mathcal{O}_K$ et des $p$-groupes finis et plats sur $\mathcal{O}_K$, puis d'appliquer cette classification pour démontrer dans le cas où $k \subseteq \overline{\mathbf{F}}_p$ (et $p \neq 2$) la conjecture de Fontaine qui dit qu'une représentation $p$-adique de $\mathrm{Gal}(\overline{K}/K)$ est cristalline à poids de Hodge-Tate dans $\{0, 1\}$ si et seulement si elle provient d'un groupe $p$-divisible sur $\mathcal{O}_K$ (cf. th. 1.4 ci-dessous).

Il y a un peu plus de vingt ans, J.-M. Fontaine a classifié tous les groupes $p$-divisibles sur $\mathcal{O}_K$ pour $e \leq p - 2$ et tous les groupes $p$-divisibles connexes sur $\mathcal{O}_K$ pour $e = p - 1$ ([Fo1], [Fo3]) ainsi que tous les $p$-groupes finis et plats sur $W$ annulés par une puissance de $p$ ([Fo2]) au moyen de ce qu'il a appelé des "systèmes de Honda". Pour $e = 1$, cette classification a été traduite quelques années après en termes de modules filtrés dans [FL, 9]. Récemment, Conrad ([Co]) a étendu la classification par les systèmes de Honda de [Fo2] aux $p$-groupes finis et plats sur $\mathcal{O}_K$ pour $e \leq p - 2$ (et pour $e \leq p - 1$ si l'on se restreint à des $p$-groupes connexes). Dans ce travail, nous généralisons la classification en termes de modules filtrés de [FL] au cas $p \neq 2$ et $e$ *quelconque* en donnant une construction directe à partir des schémas en groupes de modules filtrés "généralisés". Avant de donner plus de précisions, signalons que les difficultés dans le cas $e \geq p - 1$ sont de deux types par rapport au cas $e \leq p - 2$: d'une part la catégorie des $p$-groupes finis et plats sur $\mathcal{O}_K$ n'est plus abélienne, d'autre part le foncteur "fibre générique" de la catégorie des $p$-groupes finis et plats sur $\mathcal{O}_K$ dans celle des $p$-groupes sur $K$ n'est plus pleinement fidèle (par exemple, la représentation galoisienne associée à un $p$-groupe sur $\mathcal{O}_K$ ne détecte plus toutes les torsions "à la Tate"). Il faut donc s'abstraire totalement des représentations galoisiennes dans la construction de la classification.

Soient $u$ une indéterminée, $S$ le complété $p$-adique de l'enveloppe aux puissances divisées de $W[u]$ par rapport à l'idéal engendré par le polynôme minimal de $\pi$ sur $K_0$ et $\mathrm{Fil}^1 S$ le complété $p$-adique de cet idéal à puissances divisées. Les modules filtrés "généralisés" considérés ici sont des $S$-modules de type fini $\mathcal{M}$ munis de structures additionnelles: un sous-$S$-module $\mathrm{Fil}^1 \mathcal{M}$ contenant $\mathrm{Fil}^1 S \cdot \mathcal{M}$ et une application "semi-linéaire" $\phi_1 : \mathrm{Fil}^1 \mathcal{M} \to \mathcal{M}$ dont l'image engendre $\mathcal{M}$ sur $S$ (voir (2.1.1) pour plus de détails). Disons qu'un tel module $\mathcal{M}$ est fortement divisible s'il est libre sur $S$ et si $\mathcal{M}/\text{Fil}^1\mathcal{M}$ est sans $p$-torsion (cf. 2.1.1.8). Un des résultats principaux de cet article est (cf. théorème 4.2.2.9):

THÉORÈME 1.1. *Supposons $p \neq 2$, il y a une anti-équivalence de catégories entre la catégorie des groupes $p$-divisibles sur $\mathcal{O}_K$ et la catégorie des modules fortement divisibles.*

Comme ces modules fortement divisibles sont automatiquement "filtered free" au sens de [Fa2, 3] (2.1.1.9), l'existence d'un foncteur associant à un groupe $p$-divisible sur $\mathcal{O}_K$ un tel module était déjà connue de Faltings [Fa2, 6]. Le théorème (1.1) se déduit, par passage à la limite, de théorèmes de classification sur les $p$-groupes finis et plats sur $\mathcal{O}_K$. Le premier de ces théorèmes est le suivant (cf. théorème 3.3.7):

THÉORÈME 1.2. *Supposons $p \neq 2$, il y a une anti-équivalence de catégories entre la catégorie des schémas en groupes finis et plats sur $\mathcal{O}_K$ annulés par $p$ et la catégorie des objets $(\mathcal{M}, \text{Fil}^1\mathcal{M}, \phi_1)$ où $\mathcal{M}$ est un $S/pS$-module libre de type fini, $\text{Fil}^1\mathcal{M}$ un sous $S/pS$-module contenant $\text{Fil}^1 S \cdot \mathcal{M}$ et $\phi_1$ une application "semi-linéaire" de $\text{Fil}^1\mathcal{M}$ dans $\mathcal{M}$ dont l'image engendre $\mathcal{M}$ sur $S/pS$.*

Pour prouver (1.2), on utilise la topologie syntomique (sur $\text{Spec}(\mathcal{O}_K)$), introduite par B. Mazur dans [Ma] et utilisée pour la première fois par Fontaine et Messing dans [FM]. Une des raisons en est que le topos syntomique est très pratique pour définir et manier l'opérateur $\phi_1$. La première partie de la preuve consiste à associer à un objet $(\mathcal{M}, \text{Fil}^1\mathcal{M}, \phi_1)$ comme en (1.2) un faisceau en groupes canonique sur le site syntomique de $\text{Spec}(\mathcal{O}_K)$ et à montrer que ce faisceau est représentable (3.1). La deuxième partie consiste à associer un objet du type $(\mathcal{M}, \text{Fil}^1\mathcal{M}, \phi_1)$ à un schéma en groupes tués par $p$ et à montrer que $\mathcal{M}$ est bien libre sur $S/pS$ et que l'image de $\phi_1$ engendre bien $\mathcal{M}$ (3.2). Elle utilise de façon cruciale plusieurs résultats de Berthelot, Breen et Messing ([BBM]), valables souvent dans des situations plus générales. Par exemple, le module $\mathcal{M}$ s'interprète comme l'évaluation du cristal de Berthelot-Breen-Messing associé au schéma en groupes sur l'épaississement à puissance divisées $\text{Spec}(\mathcal{O}_K/p\mathcal{O}_K) \hookrightarrow \text{Spec}(S/pS)$ associé au choix de $\pi$, il ne dépend en particulier que de la réduction modulo $p$ du schéma en groupes, de même que $\text{Fil}^1\mathcal{M}$. Par contre, l'opérateur $\phi_1$, qui doit être pensé comme "$\frac{\phi}{p}$" ou comme un "inverse" du Verschiebung, dépend de la réduction modulo $p^2$ du schéma en groupes (alors que le Verschiebung ne dépend que de la réduction modulo $p$). Ces techniques syntomiques appliquées à l'étude des schémas en groupes étaient connues de Fontaine et Messing dans le cas $e = 1$ ([Fo4]).

Le théorème 1.2 se généralise par dévissage au cas des $p$-groupes finis et plats sur $\mathcal{O}_K$ en associant à de tels schémas en groupes certains objets $(\mathcal{M}, \text{Fil}^1\mathcal{M}, \phi_1)$ définis comme "extensions successives" d'objets annulés par $p$ du type (1.2), voir le théorème 4.2.1.6 pour un énoncé précis. Nous n'avons pas cherché à rendre ces $S$-modules plus explicites dans le cas général (ce qu'il sera peut-être intéressant de faire un jour), mais seulement dans le cas particulier important suivant (cf. théorème 4.2.2.5):

THÉORÈME 1.3. *Supposons $p \neq 2$, il y a une anti-équivalence de catégories entre la catégorie des $p$-groupes finis et plats sur $\mathcal{O}_K$ dont le noyau de la multiplication par $p^n$ est encore plat pour tout $n$ et la catégorie des objets $(\mathcal{M}, \text{Fil}^1\mathcal{M}, \phi_1)$ où $\mathcal{M}$ est un $S$-module de la forme $\oplus_{i \in I} S/p^{n_i}S$ pour $I$ fini et $n_i \in \mathbf{N}^*$, $\text{Fil}^1\mathcal{M}$ un sous $S$-module contenant $\text{Fil}^1 S \cdot \mathcal{M}$ et $\phi_1$ une application "semi-linéaire" de $\text{Fil}^1\mathcal{M}$ dans $\mathcal{M}$ dont l'image engendre $\mathcal{M}$ sur $S$.*

Par des résultats de Raynaud, la condition de platitude ci-dessus est automatique si $e < p - 1$, de sorte que dans ce cas, (1.3) fournit une nouvelle classification explicite de tous les $p$-groupes finis et plats. Cette condition est aussi satisfaite par tous les $p$-groupes provenant du noyau de la multiplication par une puissance de $p$ sur les schémas abéliens sur $\mathcal{O}_K$.

Les classifications précédentes présentent l'avantage de s'insérer naturellement dans une théorie plus vaste (voir [FL], [Br3], [Br5] ou la dernière partie) dont le but est la construction et l'étude de certaines représentations galoisiennes $p$-adiques ou de $p$-torsion, théorie qui, dans le cas présent, redonne la représentation associée aux points à valeurs dans $\overline{K}$ du schéma en groupes. Signalons que plusieurs questions restent à traiter pour avoir une vision complète de la situation: étendre ces classifications au cas $p = 2$ quitte à se restreindre à des groupes connexes, expliciter côté $S$-modules les manipulations classiques des schémas en groupes finis et plats (changement de base, dualité de Cartier), faire le lien direct avec les systèmes de Honda pour $e \leq p-2$ (classification de Fontaine et Conrad). En ce qui concerne les groupes $p$-divisibles, Zink a récemment fait le lien entre la classification (1.1) et la sienne ([Zi1]) tout en généralisant (1.1) et en incluant pour $p = 2$ le cas des groupes 2-divisibles connexes (voir [Zi2]). D'autres directions de généralisation sont encore possibles (voir la fin de (4.2.2)).

Dans la dernière partie, on combine le théorème (1.1) précédent avec les résultats de [Br5] pour démontrer (cf. théorème 5.3.2):

THÉORÈME 1.4. *Supposons $p \neq 2$, $k \subseteq \overline{\mathbf{F}}_p$ et soit $V$ une représentation $p$-adique cristalline de $\text{Gal}(\overline{K}/K)$ à poids de Hodge-Tate entre $0$ et $1$. Alors il existe un groupe $p$-divisible $H$ sur $\mathcal{O}_K$ tel que $V \simeq T_pH \otimes_{\mathbf{Z}_p} \mathbf{Q}_p$ où $T_pH$ est le module de Tate de $H$.*

Cela résout, pour $p \neq 2$ et $k \subseteq \overline{\mathbf{F}}_p$, une conjecture de Fontaine [Fo3]. Ce théorème n'était connu jusqu'à présent que pour $e < p$ (mais sans restriction sur le corps résiduel, cf. [Laf1]). L'auteur s'excuse de n'être pas arrivé à remplacer "$k \subseteq \overline{\mathbf{F}}_p$" par "$k$ parfait". L'hypothèse "$k \subseteq \overline{\mathbf{F}}_p$" est néanmoins satisfaite dans la plupart des applications. Quant à l'hypothèse $p \neq 2$, on devrait pouvoir s'en débarrasser en utilisant l'extension de (1.1) au cas des groupes 2-divisibles connexes due à Zink (cf. précédemment), le corollaire (1.12) de [Laf2] et une généralisation au cas "Fil$^{p-1}$" des résultats de [Br5] et de la dernière section (et un peu de courage). Signalons que la même méthode de preuve fournit également (cf. théorème 5.2.4):

THÉORÈME 1.5. *Supposons $k \subseteq \overline{\mathbf{F}}_p$ et soit $D$ un $(\phi, N)$-module filtré faiblement admissible de Fontaine tel que* $\text{Fil}^0(D \otimes_{K_0} K) = D \otimes_{K_0} K$ *et* $\text{Fil}^{p-1}(D \otimes_{K_0} K) = 0$, *alors $D$ est admissible.*

Ce dernier résultat vient d'être obtenu sans la restriction sur la filtration et pour $k$ parfait par Colmez et Fontaine ([CF]) grâce à des techniques $p$-adiques différentes introduites par J.-M. Fontaine. Enfin, un corollaire intéressant de (1.4) est (cf. corollaire 5.3.4):

COROLLAIRE 1.6. *Supposons $p \neq 2$ et $k \subseteq \overline{\mathbf{F}}_p$. Soient $A$ une variété abélienne définie sur $K$ et $V_p(A)$ son module de Tate tensorisé par $\mathbf{Q}_p$. Alors $A$ a réduction semi-stable (resp. bonne réduction) sur $\mathcal{O}_K$ si et seulement si $V_p(A)$ est une représentation semi-stable (resp. cristalline) de $\text{Gal}(\overline{K}/K)$.*

Plusieurs cas de ce corollaire étaient déjà connus, voir la dernière section.

Cet article est une compilation remaniée des deux prépublications [Br1] et [Br6]. Un court résumé de [Br1], i.e. des parties 2 à 4, a été publié dans [Br2]. L'auteur tient à remercier J.-M. Fontaine pour toutes ses remarques, W. Messing pour ses explications sur la théorie de Dieudonné cristalline, M. Raynaud pour ses réponses à plusieurs questions et A. Tamagawa et T. Tsuji pour lui avoir signalé le corollaire ci-dessus.

Signalons enfin qu'une variante de (1.2) a trouvé une jolie application dans [BCDT].

## 2. Préliminaires

On introduit les objets étudiés: catégories de modules, schémas en groupes finis et plats, faisceaux $\mathcal{O}^{\text{cris}}_{\infty,\pi}$ et $\mathcal{J}^{\text{cris}}_{\infty,\pi}$ et on donne quelques premières propriétés.

2.1. *Les catégories de modules.* On suppose $p \neq 2$. On introduit plusieurs catégories de modules dont on étudie en détail les objets annulés par $p$.

2.1.1. *Définitions des catégories*. Soit $W[u]$ l'anneau des polynômes en l'indéterminée $u$ et $s : W[u] \to \mathcal{O}_K$ l'unique surjection de $W$-algèbres telle que $s(u) = \pi$. On a $\mathrm{Ker}(s) = (E(u))$ où $E(u)$ est le polynôme minimal sur $K_0$ de l'uniformisante $\pi$. Notons $S$ le complété $p$-adique de l'enveloppe aux puissances divisées de $W[u]$ par rapport à $\mathrm{Ker}(s)$. Si $i \in \mathbf{N}$, notons $q(i)$ le quotient dans la division euclidienne de $i$ par $e$, on voit que $S$ s'identifie à la sous-$W$-algèbre de $K_0[[u]]$ formée des $\{\sum_{i=0}^{\infty} w_i \frac{u^i}{q(i)!}, \; w_i \in W, \lim_{i \to \infty} w_i = 0\}$. Soit $\mathrm{Fil}^1 S$ la complétion $p$-adique de l'idéal engendré par les $\gamma_i(E(u)) = E(u)^i/i!$ pour $i \geq 1$, on a $S/\mathrm{Fil}^1 S \overset{\sim}{\to} \mathcal{O}_K$, $u \mapsto \pi$. On munit $S$ de l'unique opérateur $\phi$ semi-linéaire par rapport au Frobenius sur $W$ et continu pour la topologie $p$-adique tel que $\phi(u^i/q(i)!) = u^{pi}/q(i)!$. On vérifie que $\phi(\mathrm{Fil}^1 S) \subset pS$ et on pose $\phi_1 = \frac{\phi}{p}|_{\mathrm{Fil}^1 S}$. Pour alléger les notations, on pose $S_n = S/p^n$, $v = E(u)$ et $c = \phi_1(v) \in S^*$. On note enfin $\mathrm{Fil}^p S_n$ l'idéal de $S_n$ engendré par l'image des $\gamma_i(E(u))$ pour $i \geq p$.

Soit $'(\mathrm{Mod}/S)$ la catégorie suivante: les objets sont la donnée:

- d'un $S$-module $\mathcal{M}$,

- d'un sous-$S$-module $\mathrm{Fil}^1 \mathcal{M}$ de $\mathcal{M}$ contenant $\mathrm{Fil}^1 S \cdot \mathcal{M}$,

- d'une flèche $\phi$-semi-linéaire $\phi_1 : \mathrm{Fil}^1 \mathcal{M} \to \mathcal{M}$ telle que pour tout $s \in \mathrm{Fil}^1 S$, et $x \in \mathcal{M}$, $\phi_1(sx) = \phi_1(s)\phi(x)$ où $\phi(x) = \frac{1}{c}\phi_1(vx)$.

Les flèches sont les morphismes $S$-linéaires qui préservent $\mathrm{Fil}^1 \mathcal{M}$ et commutent à $\phi_1$. On note $'(\mathrm{Mod}/S_1)$ la sous-catégorie pleine des objets tués par $p$. La catégorie $'(\mathrm{Mod}/S)$ est munie d'une notion de suite exacte courte: $0 \to \mathcal{M}' \to \mathcal{M} \to \mathcal{M}'' \to 0$ est une suite exacte dans $'(\mathrm{Mod}/S)$ si les deux suites de $S$-modules $0 \to \mathcal{M}' \to \mathcal{M} \to \mathcal{M}'' \to 0$ et $0 \to \mathrm{Fil}^1 \mathcal{M}' \to \mathrm{Fil}^1 \mathcal{M} \to \mathrm{Fil}^1 \mathcal{M}'' \to 0$ sont exactes.

On note $(\mathrm{Mod}/S_1)$ la sous-catégorie pleine de $'(\mathrm{Mod}/S_1)$ formée des objets $\mathcal{M}$ qui vérifient les deux conditions supplémentaires:

- $\mathcal{M}$ est un $S_1$-module libre de type fini,

- $\phi_1(\mathrm{Fil}^1 \mathcal{M})$ engendre $\mathcal{M}$ sur $S_1$.

On note $(\mathrm{Mod}/S)$ la sous-catégorie pleine de $'(\mathrm{Mod}/S)$ stable par extension engendrée par $(\mathrm{Mod}/S_1)$. Par (2.1.1.2), les objets tués par $p$ de $(\mathrm{Mod}/S)$ sont encore les objets de $(\mathrm{Mod}/S_1)$.

LEMME 2.1.1.1. *Tout objet $\mathcal{M}$ de $(\mathrm{Mod}/S)$ est un $S$-module de type fini annulé par une puissance de $p$ tel que $\phi_1(\mathrm{Fil}^1 \mathcal{M})$ engendre $\mathcal{M}$ sur $S$.*

*Preuve.* Si $0 \to \mathcal{M}' \to \mathcal{M} \to \mathcal{M}'' \to 0$ est une extension dans $'(\text{Mod}/S)$ et si $\phi_1(\text{Fil}^1\mathcal{M}')$ engendre $\mathcal{M}'$ (resp. avec $\mathcal{M}''$), il en est encore de même de $\mathcal{M}$, d'où le résultat. $\square$

LEMME 2.1.1.2. *La catégorie* $(\text{Mod}/S_1)$ *est stable par extension dans* $'(\text{Mod}/S_1)$.

*Preuve.* Si $0 \to \mathcal{M}' \to \mathcal{M} \to \mathcal{M}'' \to 0$ est une extension dans $'(\text{Mod}/S_1)$ et si $\mathcal{M}'$ et $\mathcal{M}''$ sont libres sur $S_1$, il en est de même de $\mathcal{M}$, d'où le résultat par (2.1.1.1). $\square$

Soit $(\text{Mod FI}/S)$ (Modules à "Facteurs Invariants") la sous-catégorie pleine de $'(\text{Mod}/S)$ formée des objets $\mathcal{M}$ qui vérifient les deux conditions:

- le $S$-module $\mathcal{M}$ est de la forme $\mathcal{M} \simeq \oplus_{i \in I} S_{n_i}$ pour $I$ fini et $n_i \in \mathbf{N}^*$,
- $\phi_1(\text{Fil}^1\mathcal{M})$ engendre $\mathcal{M}$ sur $S$.

Les objets tués par $p$ de $(\text{Mod FI}/S)$ sont encore ceux de $(\text{Mod}/S_1)$.

LEMME 2.1.1.3. *Soit $\mathcal{M}$ un objet de* $(\text{Mod FI}/S)$, *alors pour tout* $r \in \mathbf{N}$, $p^r \text{Fil}^1\mathcal{M} = \text{Fil}^1\mathcal{M} \cap p^r\mathcal{M}$.

*Preuve.* Elle est repoussée à la fin de la section qui suit car elle nécessite l'étude des objets annulés par $p$. $\square$

COROLLAIRE 2.1.1.4. *La catégorie* $(\text{Mod FI}/S)$ *est une sous-catégorie pleine de* $(\text{Mod}/S)$ (*en général différente si* $e \geq p-1$).

*Preuve.* Par (2.1.1.3), on a $\text{Fil}^1\mathcal{M}/p^r\text{Fil}^1\mathcal{M} \hookrightarrow \mathcal{M}/p^r\mathcal{M}$ si $\mathcal{M}$ est dans $(\text{Mod FI}/S)$, donc $(\mathcal{M}/p^r\mathcal{M}, \text{Fil}^1\mathcal{M}/p^r\text{Fil}^1\mathcal{M}, \phi_1)$ est un objet de $(\text{Mod FI}/S)$ pour tout $r \in \mathbf{N}^*$. Comme $(p^r\mathcal{M}, p^r\text{Fil}^1\mathcal{M}, \phi_1)$ est aussi un objet de $(\text{Mod FI}/S)$, on déduit que $(\text{Mod FI}/S)$ est une sous-catégorie pleine de $(\text{Mod}/S)$. Dans le cas $e = p-1$, construisons un objet $\mathcal{M}$ de $(\text{Mod}/S)$ qui n'est pas dans $(\text{Mod FI}/S)$. Pour $e = p-1$, soient $S_1(1)$ le module $S_1 e_1$ muni de $\text{Fil}^1 S_1(1) = S_1 e_1$ et $\phi_1(e_1) = ce_1$, $S_1(0)$ le module trivial $S_1 e_2$ muni de $\text{Fil}^1 S_1(0) = \text{Fil}^1 S_1 e_2$ et $\phi_1(se_2) = \phi_1(s)e_2$ si $s \in \text{Fil}^1 S_1$, on vérifie que la flèche $S_1$-linéaire $S_1(1) \to S_1(0)$ induite par $e_1 \mapsto u^p e_2$ est un morphisme dans $'(\text{Mod}/S_1)$. Soit $S_2(1) = S_2 e_1$ avec $\text{Fil}^1 S_2(1) = S_2 e_1$ et $\phi_1(e_1) = ce_1$, on définit $\mathcal{M}$ comme le conoyau de $0 \to S_1(1) \to S_2(1) \oplus S_1(0), e_1 \mapsto pe_1 \oplus u^p e_2$ muni du $\text{Fil}^1$ conoyau des $\text{Fil}^1$ et du $\phi_1$ induit ($\mathcal{M}$ est bien défini). C'est un objet de $(\text{Mod}/S)$, car on a une extension dans $'(\text{Mod}/S)$:

$$\begin{array}{ccccccccc}
0 & \to & S_1(0) & \to & \mathcal{M} & \to & S_1(1) & \to & 0 \\
 & & e_2 & \mapsto & \overline{0 \oplus e_2} & \mapsto & 0 & & \\
 & & & & \overline{e_1 \oplus 0} & \mapsto & e_1 & &
\end{array}$$

et il n'est clairement pas dans $(\text{Mod FI}/S)$. $\square$

*Remarque* 2.1.1.5. Lorsque $c \equiv 1 \bmod p$, on verra que le morphisme $S_1(1) \to S_1(0)$ correspond à un morphisme non trivial de schémas en groupes $\mathbf{Z}/p\mathbf{Z} \to \mu_p$ (3.1).

*Remarque* 2.1.1.6. On peut montrer, par des manipulations d'algèbre linéaire pas toujours faciles, que $(\text{Mod FI}/S) = (\text{Mod}/S)$ lorsque $e \leq p - 2$ et que la catégorie $(\text{Mod FI}/S)$ est alors abélienne (pour ce dernier fait, la preuve est essentiellement la même qu'en [Br3, 2], voir aussi (4.2.2.6)). J'ignore s'il existe une description totalement explicite de tous les objets de $(\text{Mod}/S)$ lorsque $e \geq p - 1$ (mis à part les objets tués par $p$).

*Remarque* 2.1.1.7. Avec les notations de [Br3], on avait $'\underline{\mathcal{M}}_0^1$ (resp. $\underline{\mathcal{M}}_0^1$, $'\underline{\mathcal{M}}_{0,k}^1$, $\underline{\mathcal{M}}_{0,k}^1$) au lieu de $'(\text{Mod}/S)$ (resp. $(\text{Mod FI}/S)$, $'(\text{Mod}/S_1)$, $(\text{Mod}/S_1)$).

On introduit enfin une dernière définition:

*Définition* 2.1.1.8. On appelle module fortement divisible tout objet $\mathcal{M}$ de $'(\text{Mod}/S)$ vérifiant les trois conditions supplémentaires:

- le $S$-module $\mathcal{M}$ est libre de rang fini,
- le $S$-module $\mathcal{M}/\text{Fil}^1\mathcal{M}$ est sans $p$-torsion,
- le $S$-module $\mathcal{M}$ est engendré par $\phi_1(\text{Fil}^1\mathcal{M})$.

Ces modules sont tous "filtered free" au sens de [Fa2, 3]:

LEMME 2.1.1.9. *Soit $\mathcal{M}$ un module fortement divisible, il existe une base $(e_1, \ldots, e_d)$ de $\mathcal{M}$ sur $S$ et un entier $d_1 \in \{1, \ldots, d\}$ tels que:*

$$\text{Fil}^1\mathcal{M} = (\oplus_{i=1}^{d_1}\text{Fil}^1Se_i) \oplus (\oplus_{i=d_1+1}^{d}Se_i).$$

*Preuve.* Via l'identification $S/\text{Fil}^1S \xrightarrow{\sim} \mathcal{O}_K$, $u \mapsto \pi$, on a une suite exacte de $\mathcal{O}_K$-modules libres de type fini: $0 \to \text{Fil}^1\mathcal{M}/\text{Fil}^1S\mathcal{M} \to \mathcal{M}/\text{Fil}^1S\mathcal{M} \to \mathcal{M}/\text{Fil}^1\mathcal{M} \to 0$. Soient $(\overline{e}_1, \ldots, \overline{e}_d)$ une base de $\mathcal{M}/\text{Fil}^1S\mathcal{M}$ sur $\mathcal{O}_K$ telle que $\text{Fil}^1\mathcal{M}/\text{Fil}^1S\mathcal{M} = \oplus_{i=d_1+1}^{d}\mathcal{O}_K\overline{e}_i$, $(e_1, \ldots, e_d)$ des relevés de $(\overline{e}_1, \ldots, \overline{e}_d)$ dans $\mathcal{M}$ tels que $(e_{d_1+1}, \ldots, e_d) \in \text{Fil}^1\mathcal{M}$, $(f_1, \ldots, f_d)$ une base quelconque de $\mathcal{M}$ sur $S$ et $(\overline{f}_1, \ldots, \overline{f}_d)$ la base image dans $\mathcal{M}/\text{Fil}^1S\mathcal{M}$. La matrice de $(\overline{e}_1, \ldots, \overline{e}_d)$ dans $(\overline{f}_1, \ldots, \overline{f}_d)$ a par définition son déterminant dans $\mathcal{O}_K^*$, donc le déterminant de la matrice de $(e_1, \ldots, e_d)$ dans la base $(f_1, \ldots, f_d)$ est dans $S^*$ et $(e_1, \ldots, e_d)$ est une base de $\mathcal{M}$ sur $S$. Si $x \in \text{Fil}^1\mathcal{M}$, son image $\overline{x}$ dans $\mathcal{M}/\text{Fil}^1S\mathcal{M}$ s'écrit $\overline{x} = \sum_{i=d_1+1}^{d}\overline{s}_i\overline{e}_i$, $\overline{s}_i \in \mathcal{O}_K$, donc $x - \sum_{i=d_1+1}^{d}s_ie_i \in \text{Fil}^1S\mathcal{M}$ où $s_i$ relève $\overline{s}_i$ dans $S$, donc $x \in \text{Fil}^1S\mathcal{M} + \sum_{i=d_1+1}^{d}Se_i$ d'où $\text{Fil}^1\mathcal{M} = \text{Fil}^1S\mathcal{M} + \sum_{i=d_1+1}^{d}Se_i$. $\square$

*Remarque* 2.1.1.10. La définition 2.1.1.8 se généralise à des crans supérieurs de filtration, voir (5.1.2), mais les modules fortement divisibles ne sont plus alors "filtered free" en général (cf. [Br3, 6]).

Si $\mathcal{M}$ est un module fortement divisible, notons que pour tout $n \in \mathbf{N}^*$, $\mathcal{M}/p^n\mathcal{M}$ est muni d'une structure d'objet de $(\mathrm{Mod\,FI}/S)$ en posant:

$$\mathrm{Fil}^1(\mathcal{M}/p^n\mathcal{M}) = \mathrm{Fil}^1\mathcal{M}/p^n\mathrm{Fil}^1\mathcal{M} \hookrightarrow \mathcal{M}/p^n\mathcal{M}$$

avec le $\phi_1$ induit.

2.1.2. *Les objets annulés par $p$.* Considérons d'un peu plus près la catégorie $(\mathrm{Mod}/S_1)$. Soient $F(u)$ l'unique polynôme dans $W[u]$ de degré $\leq e-1$ tel que $E(u) = u^e - pF(u)$, $\tilde{S}_1 = k[u]/u^{ep}$ et $\tilde{c}$ l'image de $-\phi(F(u)) \in W[u]$ dans $k[u]/u^{ep}$: c'est une unité de $\tilde{S}_1$. Notons $\mathrm{Fil}^1(\tilde{S}_1) = u^e\tilde{S}_1$ et $\tilde{\phi}_1 : \mathrm{Fil}^1(\tilde{S}_1) \to \tilde{S}_1$ l'unique opérateur semi-linéaire (par rapport au Frobenius naturel de $\tilde{S}_1$) tel que $\tilde{\phi}_1(u^e) = \tilde{c}$. Soit $'(\mathrm{Mod}/\tilde{S}_1)$ la catégorie dont les objets sont la donnée:

- d'un $\tilde{S}_1$-module $\tilde{\mathcal{M}}$,
- d'un sous-$\tilde{S}_1$-module $\mathrm{Fil}^1\tilde{\mathcal{M}}$ contenant $u^e\tilde{\mathcal{M}}$,
- d'une application $\tilde{S}_1$-semi-linéaire $\tilde{\phi}_1 : \mathrm{Fil}^1\tilde{\mathcal{M}} \to \tilde{\mathcal{M}}$,

et dont les flèches sont les morphismes $\tilde{S}_1$-linéaires qui préservent $\mathrm{Fil}^1\tilde{\mathcal{M}}$ et commutent à $\tilde{\phi}_1$. En procédant comme pour $'(\mathrm{Mod}/S)$, on munit cette catégorie d'une notion de suite exacte courte. On note $(\mathrm{Mod}/\tilde{S}_1)$ la sous-catégorie pleine de $'(\mathrm{Mod}/\tilde{S}_1)$ formée des objets $\mathcal{M}$ qui vérifient les deux conditions supplémentaires:

- le $\tilde{S}_1$-module $\tilde{\mathcal{M}}$ est libre de rang fini,
- $\tilde{\phi}_1(\mathrm{Fil}^1\tilde{\mathcal{M}})$ engendre $\tilde{\mathcal{M}}$ sur $\tilde{S}_1$.

LEMME 2.1.2.1. *Soit $\tilde{\mathcal{M}}$ un objet de $(\mathrm{Mod}/\tilde{S}_1)$, alors:*

(1) *La flèche $\mathrm{Id} \otimes \tilde{\phi}_1 : k[u]/u^{ep} \otimes_{(\phi), k[u]/u^e} \mathrm{Fil}^1\tilde{\mathcal{M}}/u^e\mathrm{Fil}^1\tilde{\mathcal{M}} \to \tilde{\mathcal{M}}$ est un isomorphisme de $k[u]/u^{ep}$-modules.*

(2) *La flèche canonique $k[u]/u^{ep} \otimes_{k[u^p]/u^{ep}} \tilde{\phi}_1(\mathrm{Fil}^1\tilde{\mathcal{M}}) \to \tilde{\mathcal{M}}$ est un isomorphisme de $k[u]/u^{ep}$-modules.*

*Preuve.* La preuve étant similaire à celle de [Br3, 2.2.1.2], nous y renvoyons le lecteur. □

On a un isomorphisme $k[u]/u^{ep}$-linéaire $S_1 \simeq \tilde{S}_1[X_i]/(X_i^p)_{i \in \mathbf{N}^*}$ où $\gamma_{p^i}(u^e) \mapsto X_i$. Soit $\sigma$ la surjection de $k[u]/u^{ep}$-algèbres $S_1 \to \tilde{S}_1$, $X_i \mapsto 0$ pour tout $i$.

Pour $\mathcal{M}$ dans $(\mathrm{Mod}/S_1)$, posons $T(\mathcal{M}) = \tilde{S}_1 \otimes_{(\sigma),S_1} \mathcal{M}$ et soit $s$ la surjection canonique $\mathcal{M} \to T(\mathcal{M})$. On définit:

- $\mathrm{Fil}^1 T(\mathcal{M}) = s(\mathrm{Fil}^1 \mathcal{M})$,

- si $x \in s(\mathrm{Fil}^1 \mathcal{M})$, $\tilde{\phi}_1(x) = s(\phi_1(\hat{x}))$ où $\hat{x}$ est un relevé quelconque de $x$ dans $\mathrm{Fil}^1 \mathcal{M}$, c'est indépendant du relevé car $p \geq 3$.

On obtient ainsi un foncteur $T : (\mathrm{Mod}/S_1) \to (\mathrm{Mod}/\tilde{S}_1)$.

PROPOSITION 2.1.2.2. *Le foncteur $T$ induit une équivalence de catégories entre $(\mathrm{Mod}/S_1)$ et $(\mathrm{Mod}/\tilde{S}_1)$.*

*Preuve.* La preuve est similaire à celle de [Br3, 2.2.2.1]. Un quasi-inverse à $T$ est donné par $\mathcal{M} = S_1 \otimes_{\tilde{S}_1} \tilde{\mathcal{M}}$, $\mathrm{Fil}^1 \mathcal{M} = s^{-1}(\mathrm{Fil}^1 \tilde{\mathcal{M}})$ où $s = $ "$\sigma \otimes \mathrm{Id}$": $\mathcal{M} \to \tilde{\mathcal{M}}$ et $\phi_1(x) = 1 \otimes \tilde{\phi}_1(s(x))$ si $x \in \mathrm{Fil}^1 \mathcal{M}$. Pour plus de détails, voir [Br3, 2.2.2.1]. □

*Remarque* 2.1.2.3. Comme dans [Br3, 2.2.2.2], on a une version de (2.1.2.1) directement dans $(\mathrm{Mod}/S_1)$: pour tout $\mathcal{M}$ de $(\mathrm{Mod}/S_1)$, on a des isomorphismes $S_1$-linéaires $\mathrm{Id} \otimes \phi_1 : S_1 \otimes_{(\phi),S_1} (\mathrm{Fil}^1 \mathcal{M}/(u^e \mathrm{Fil}^1 \mathcal{M} + \mathrm{Fil}^p S_1 \mathcal{M})) \xrightarrow{\sim} \mathcal{M}$ et $S_1 \otimes_{k[u^p]/u^{ep}} \phi_1(\mathrm{Fil}^1 \mathcal{M}) \xrightarrow{\sim} \mathcal{M}$.

LEMME 2.1.2.4. *Soient $\tilde{\mathcal{M}}$ un objet de $(\mathrm{Mod}/\tilde{S}_1)$ de rang $d$ sur $\tilde{S}_1$, $(f_1, \ldots, f_d)$ $d$ éléments de $\tilde{\mathcal{M}}$ et $(r_1, \ldots, r_d)$ $d$ entiers dans $\{0, \ldots, ep-1\}$ tels que $(u^{r_1} f_1, \ldots, u^{r_d} f_d)$ engendre $\mathrm{Fil}^1 \tilde{\mathcal{M}}$. Alors, $\forall i \in \{1, \ldots, d\}$, $r_i$ est le plus petit entier tel que $u^{r_i} f_i \in \mathrm{Fil}^1 \tilde{\mathcal{M}}$.*

*Preuve.* En effet, si tel n'était pas le cas pour un $i$, on aurait $u^{r_i - 1} f_i \in \mathrm{Fil}^1 \tilde{\mathcal{M}}$ soit $u^{r_i - 1} f_i = \sum_{j=1}^{d} \lambda_j u^{r_j} f_j$ d'où $u^{r_i - 1} f_i \in \sum_{j \neq i} \tilde{S}_1 u^{r_j} f_j$; i.e., les $u^{r_j} f_j$ pour $j \neq i$ engendrent $\mathrm{Fil}^1 \tilde{\mathcal{M}}$, ce qui est impossible car $\mathrm{Fil}^1 \tilde{\mathcal{M}}/u^e \mathrm{Fil}^1 \tilde{\mathcal{M}}$ est un $k[u]/u^e$-module libre de rang $d$ (2.1.2.1). □

PROPOSITION 2.1.2.5. *Soit $\tilde{\mathcal{M}}$ (resp. $\mathcal{M}$) un objet de $(\mathrm{Mod}/\tilde{S}_1)$ (resp. un objet de $(\mathrm{Mod}/S_1)$) de rang $d$ sur $\tilde{S}_1$ (resp. $S_1$). Il existe une base $(e_1, \ldots, e_d)$ de $\tilde{\mathcal{M}}$ (resp. $\mathcal{M}$) et des entiers $(r_1, \ldots, r_d)$ dans $\{0, \ldots, e\}$ tels que $(u^{r_1} e_1, \ldots, u^{r_d} e_d)$ (resp. $(u^{r_1} e_1, \ldots, u^{r_d} e_d) + \mathrm{Fil}^p S_1 \mathcal{M}$) engendre $\mathrm{Fil}^1 \tilde{\mathcal{M}}$ (resp. $\mathrm{Fil}^1 \mathcal{M}$).*

*Preuve.* Par (2.1.2.2), il suffit de le voir pour $\tilde{\mathcal{M}}$. Soit $(g_1, \ldots, g_d)$ dans $\mathrm{Fil}^1 \tilde{\mathcal{M}}$ dont l'image dans $\mathrm{Fil}^1 \tilde{\mathcal{M}}/u^e \mathrm{Fil}^1 \tilde{\mathcal{M}}$ forme une base de ce $k[u]/u^e$-module libre de rang $d$. On peut écrire $g_i = u^{r_i} f_i$ avec $f_i \notin u \tilde{\mathcal{M}}$ et permuter les $i$ pour avoir $r_i \leq r_{i+1}$. Pour un $x \in \tilde{\mathcal{M}}$, notons $\overline{x}$ son image dans $\tilde{\mathcal{M}}/u\tilde{\mathcal{M}}$. Si $(\overline{f}_1, \overline{f}_2)$ sont liés, il existe $k_1$ dans $k^*$ tels que $k_1 f_1 + f_2 = u^\epsilon f_2'$ avec $\epsilon \geq 1$ et $f_2' \notin u\tilde{\mathcal{M}}$.

La famille $(u^{r_1}f_1, u^{r_2+\epsilon}f'_2, u^{r_3}f_3, \ldots)$ engendre encore $\mathrm{Fil}^1\tilde{\mathcal{M}}$. Si $(\overline{f}_1, \overline{f}'_2)$ sont liés, on recommence le raisonnement ci-dessus. Comme l'exposant $r'_2 = r_2 + \epsilon$ ne peut augmenter indéfiniment à cause de $u^e\tilde{\mathcal{M}} \subset \mathrm{Fil}^1\tilde{\mathcal{M}}$ et de (2.1.2.4), on obtient à terme un nouvel élément $u^{r_2}f_2$ tel que $(\overline{f}_1, \overline{f}_2)$ est libre dans $\tilde{\mathcal{M}}/u\tilde{\mathcal{M}}$. Si le nouvel exposant $r_2$ ne vérifie plus $r_2 \leq r_3$, on recommence ce qui précède avec $(u^{r_1}f_1, u^{r_3}f_3)$ et ainsi de suite. En réordonnant, on obtient à la fin une nouvelle famille génératrice $(u^{r_i}f_i)_{1\leq i \leq d}$ avec $r_i \leq r_{i+1}$ et $(\overline{f}_1, \overline{f}_2)$ libre. On recommence: si $(\overline{f}_1, \overline{f}_2, \overline{f}_3)$ sont liés, on a $(k_1, k_2) \neq (0,0) \in k^2$ tels que $k_1 f_1 + k_2 f_2 + f_3 = u^\epsilon f'_3$ et la famille $(u^{r_1}f_1, u^{r_2}f_2, u^{r_3+\epsilon}f'_3, \ldots)$ engendre $\mathrm{Fil}^1\tilde{\mathcal{M}}$, etc. En procédant comme précédemment, on se ramène à $(u^{r_i}f_i)_{1\leq i \leq d}$ ($r_i \leq r_{i+1}$) qui engendre $\mathrm{Fil}^1\tilde{\mathcal{M}}$ et tel que $(\overline{f}_1, \overline{f}_2, \overline{f}_3)$ est libre. Par une récurrence, on obtient finalement une famille génératrice de $\mathrm{Fil}^1\tilde{\mathcal{M}}$, $(u^{r_i}f_i)_{1\leq i \leq d}$, telle que $(\overline{f}_i)_{1\leq i \leq d}$ est une base de $\tilde{\mathcal{M}}/u\tilde{\mathcal{M}}$. Donc $e_i = f_i$ est une base de $\tilde{\mathcal{M}}$ qui convient. □

*Définition* 2.1.2.6. On appelle base adaptée d'un objet $\tilde{\mathcal{M}}$ (resp. d'un objet $\mathcal{M}$) de $(\mathrm{Mod}/\tilde{S}_1)$ (resp. $(\mathrm{Mod}/S_1)$) toute base de $\tilde{\mathcal{M}}$ (resp. $\mathcal{M}$) satisfaisant la condition de (2.1.2.5).

Nous terminons avec la preuve du lemme 2.1.1.3 de la section précédente. Soient $\tilde{\mathcal{M}}$, $\tilde{\mathcal{M}}'$ deux objets de $(\mathrm{Mod}/\tilde{S}_1)$ tels qu'on ait un morphisme $\tilde{\mathcal{M}}' \to \tilde{\mathcal{M}}$ dans $'(\mathrm{Mod}/\tilde{S}_1)$ injectif sur les modules sous-jacents. Alors $\mathrm{Fil}^1\tilde{\mathcal{M}}'$ $= \tilde{\mathcal{M}}' \cap \mathrm{Fil}^1\tilde{\mathcal{M}}$. En effet, soit $x \in \tilde{\mathcal{M}}' \cap \mathrm{Fil}^1\tilde{\mathcal{M}}$ et $r$ le plus petit entier tel que $u^r x \in \mathrm{Fil}^1\tilde{\mathcal{M}}'$ ($0 \leq r \leq e$). On a $\tilde{\phi}_1(u^r x) = u^{pr}\tilde{\phi}_1(x)$ dans $\tilde{\mathcal{M}}$, donc $u^{p(e-r)}\tilde{\phi}_1(u^r x) = 0$ dans $\tilde{\mathcal{M}}'$, i.e. $\tilde{\phi}_1(u^r x) \in u^{pr}\tilde{\mathcal{M}}'$ puisque $\tilde{\mathcal{M}}'$ est libre sur $\tilde{S}_1$. De (2.1.2.1,(1)), on déduit facilement $\overline{u^r x} \in u^r\mathrm{Fil}^1\tilde{\mathcal{M}}'/u^e\mathrm{Fil}^1\tilde{\mathcal{M}}'$, i.e. $u^r x \in u^r\mathrm{Fil}^1\tilde{\mathcal{M}}'$ d'où $x \in \mathrm{Fil}^1\tilde{\mathcal{M}}' + u^{ep-r}\tilde{\mathcal{M}}' \subset \mathrm{Fil}^1\tilde{\mathcal{M}}'$ qui entraîne $r = 0$. De la construction en (2.1.2.2), on déduit le même résultat dans $(\mathrm{Mod}/S_1)$ pour une injection $\mathcal{M}' \hookrightarrow \mathcal{M}$. Soit maintenant $\mathcal{M}$ un objet de $(\mathrm{Mod}\,\mathrm{FI}/S)$, la même preuve qu'en ([Br3],2.3.1.2), qui utilise ce qui précède, montre alors que pour tout $r \in \mathbf{N}$: $p^r\mathrm{Fil}^1\mathcal{M} = \mathrm{Fil}^1\mathcal{M} \cap p^r\mathcal{M}$.

2.2. *Schémas en groupes et site syntomique.* On rappelle qu'un morphisme de schémas $X \to Y$ est dit syntomique s'il est plat, localement de présentation finie et s'il se factorise localement en une immersion fermée régulière dans un $Y$-schéma lisse ([Ma], [FM]). Les principales propriétés de ces morphismes sont les suivantes:

(1) Ils sont stables par composition et changement de base.

(2) Dans la définition, on peut remplacer immersion régulière par immersion transversalement régulière (E.G.A. IV.19.2).

(3) Si $X \to Y$ est syntomique, toute factorisation $X \stackrel{i}{\hookrightarrow} Z \stackrel{f}{\to} Y$ avec $i$ immersion fermée et $f$ lisse est telle que $i$ est régulière.

(4) Soit $Y' \hookrightarrow Y$ une immersion fermée et $X'$ syntomique sur $Y'$, il existe un recouvrement $X_i'$ de $X'$ pour la topologie de Zariski et des schémas $X_i$ syntomiques sur $Y$ tels que $Y' \times_Y X_i = X_i'$.

On appelle schéma formel $p$-adique, ou schéma formel tout court, sur $\mathrm{Sp}f(\mathcal{O}_K)$ un schéma formel sur $\mathrm{Sp}f(\mathcal{O}_K)$ qui localement s'identifie au schéma formel affine associé à une $\mathcal{O}_K$-algèbre séparée et complète pour la topologie $p$-adique. Un morphisme de schémas formels $p$-adiques $\mathfrak{X} \to \mathfrak{Y}$ sur $\mathrm{Sp}f(\mathcal{O}_K)$ est dit syntomique si pour tout $n \in \mathbf{N}^*$ le morphisme de schémas $\mathfrak{X}/p^n \to \mathfrak{Y}/p^n$ est syntomique ([FM, III.4.1]) ou, de manière équivalente, $\mathfrak{X}/\pi^n \to \mathfrak{Y}/\pi^n$ syntomique pour tout $n \in \mathbf{N}^*$: ces morphismes sont stables par composition et changement de base. On note $\mathrm{Sp}f(\mathcal{O}_K)_{\mathrm{syn}}$ la catégorie des schémas formels $p$-adiques syntomiques sur $\mathrm{Sp}f(\mathcal{O}_K)$ munie de la topologie de Grothendieck engendrée par les familles surjectives de morphismes syntomiques. Notons $\mathcal{O}_K\{X_1, \ldots, X_r\}$ le complété $p$-adique de l'anneau de polynômes $\mathcal{O}_K[X_1, \ldots, X_r]$.

LEMME 2.2.1. *Soit $\mathfrak{X}$ un schéma formel $p$-adique syntomique sur $\mathrm{Sp}f(\mathcal{O}_K)$, alors localement $\mathfrak{X}$ s'identifie au schéma formel $p$-adique affine associé au quotient d'un anneau $\mathcal{O}_K\{X_1, \ldots, X_r\}$ par une suite transversalement régulière relativement à $\mathcal{O}_K$.*

*Preuve.* Il résulte des définitions ci-dessus qu'on peut recouvrir $\mathfrak{X}$ par des schémas formels affines $\mathrm{Sp}f(B)$ où $B$ est une $\mathcal{O}_K$-algèbre séparée et complète pour la topologie $p$-adique telle que $B/\pi^n B$ est plat sur $\mathcal{O}_K/\pi^n \mathcal{O}_K$ pour tout $n \in \mathbf{N}^*$ et $B/\pi B \simeq k[X_1, \ldots, X_r]/(f_1, \ldots, f_s)$ avec $(f_1, \ldots, f_s)$ régulière. Puisque $\cap_{n \in \mathbf{N}^*} \pi^n B = \{0\}$, les platitudes entraînent que $B$ est sans $\pi$-torsion. Puisque $B$ est complet, si $b_1, \ldots, b_r$ désignent des relevés quelconques des $X_i$ dans $B$, on déduit que l'homomorphisme $\mathcal{O}_K$-linéaire $\mathcal{O}_K\{X_1, \ldots, X_r\} \to B$, $X_i \mapsto b_i$ est surjectif. Soit $I$ l'idéal noyau, une chasse au diagramme montre que l'isomorphisme

$$\mathcal{O}_K\{X_1, \ldots, X_r\}/\pi \mathcal{O}_K\{X_1, \ldots, X_r\} \stackrel{\sim}{\to} k[X_1, \ldots, X_r]$$

induit $I/\pi I \simeq (f_1, \ldots, f_s)$. Soient $\hat{f}_i$ des relevés quelconques des $f_i$ dans $I$, comme $I$ est complet (pour la topologie $p$-adique) et $B$ sans $\pi$-torsion, on déduit $I = (\hat{f}_1, \ldots, \hat{f}_s)$. Une récurrence facile et le critère classique de platitude par les suites régulières ([Mi, I.2.5]) permet de déduire la transverse régularité de la suite $(\hat{f}_1, \ldots, \hat{f}_s)$ par rapport à $\mathcal{O}_K$ de la régularité de $(f_1, \ldots, f_s)$. □

Soit $G$ un schéma en groupes (commutatifs) fini et plat sur $\mathrm{Spec}(\mathcal{O}_K)$:

PROPOSITION 2.2.2.  *$G$ est un objet de $\mathrm{Sp}f(\mathcal{O}_K)_{\mathrm{syn}}$.*

*Preuve* (voir [Ma]). Comme $G$ est fini sur $\mathrm{Spec}(\mathcal{O}_K)$, c'est un schéma formel $p$-adique sur $\mathrm{Sp}f(\mathcal{O}_K)$. Comme $G$ est plat sur $\mathrm{Spec}(\mathcal{O}_K)$, un examen de la preuve de (2.2.1) montre qu'il suffit de montrer $G_k \to \mathrm{Spec}(k)$ syntomique où $G_k = G \times_{\mathrm{Spec}(\mathcal{O}_K)} \mathrm{Spec}(k)$. Quitte à remplacer $k$ par une clôture algébrique, on peut supposer $k$ algébriquement clos (E.G.A. IV.19.3.4). Puisque $G_k$ est fini sur $\mathrm{Spec}(k)$, ses points sont tous fermés (E.G.A. I.6.4.4) et sa bigèbre s'identifie au produit de ses anneaux locaux (E.G.A. I.6.2.2). Comme $k$ est algébriquement clos, ces anneaux locaux sont tous isomorphes à l'anneau local en la section unité (faire une translation). Mais cet anneau local (fini) admet une présentation de la forme $k[X_1, \ldots, X_r]/(X_1^{p^{n_1}}, \ldots, X_r^{p^{n_r}})$ (S.G.A. III, $\mathrm{VII}_B$.5.4) donc est syntomique sur $\mathrm{Spec}(k)$, d'où le résultat. □

Tout schéma en groupe commutatifs $G$ sur $\mathrm{Spec}(\mathcal{O}_K)$ représente un faisceau abélien pour le gros site fppf sur $\mathrm{Spec}(\mathcal{O}_K)$ ([Mi, II.1.7]). En particulier, si $G$ est fini et plat, le préfaisceau $\mathfrak{X} \mapsto G(\mathfrak{X}) = \mathrm{Hom}_{\mathrm{Sp}f(\mathcal{O}_K)}(\mathfrak{X}, G)$ est un faisceau sur $\mathrm{Sp}f(\mathcal{O}_K)_{\mathrm{syn}}$ et le foncteur qui à $G$ associe le faisceau correspondant sur $\mathrm{Sp}f(\mathcal{O}_K)_{\mathrm{syn}}$ est pleinement fidèle par (2.2.2). Notons $(p\text{-}\mathrm{Gr}/\mathcal{O}_K)$ la catégorie des $p$-groupes finis et plats sur $\mathrm{Spec}(\mathcal{O}_K)$ et $(\mathrm{Ab}/\mathcal{O}_K)$ la catégorie (abélienne) des faisceaux de groupes commutatifs sur $\mathrm{Sp}f(\mathcal{O}_K)_{\mathrm{syn}}$. Si $G, G'$ et $G''$ sont des schémas en groupes (sur $\mathrm{Spec}(\mathcal{O}_K)$), rappelons que dire que $0 \to G' \to G \to G'' \to 0$ est une suite exacte c'est par définition dire que le diagramme correspondant de faisceaux abéliens sur le *gros site* fppf de $\mathrm{Spec}(\mathcal{O}_K)$ est une suite exacte.

PROPOSITION 2.2.3. *Soient $G, G'$ et $G''$ dans $(p\text{-}\mathrm{Gr}/\mathcal{O}_K)$, alors $0 \to G' \to G \to G'' \to 0$ est une suite exacte dans $(p\text{-}\mathrm{Gr}/\mathcal{O}_K)$ si et seulement si c'est une suite exacte dans $(\mathrm{Ab}/\mathcal{O}_K)$.*

*Preuve.* Soit $0 \to G' \to G \to G'' \to 0$ une suite exacte dans $(p\text{-}\mathrm{Gr}/\mathcal{O}_K)$, il faut montrer que la suite de faisceaux correspondante sur $\mathrm{Sp}f(\mathcal{O}_K)_{\mathrm{syn}}$ est encore exacte. Seule la surjectivité de droite n'est pas évidente. Pour la montrer, il suffit de montrer que $G \to G''$ est un morphisme syntomique et surjectif de schémas formels $p$-adiques i.e. que $G_n \to G_n''$ est un morphisme syntomique et surjectif de schémas pour tout $n \geq 1$ où

$$G_n = G \times_{\mathrm{Spec}(\mathcal{O}_K)} \mathrm{Spec}(\mathcal{O}_K/\pi^n)$$

(resp. $G_n'' = \cdots$). La surjectivité (sur les espaces topologiques) découle de la surjectivité des faisceaux pour la topologie fppf et la platitude est un résultat classique ([DG, III.3.2.5]). On a un isomorphisme canonique

$$G_n \times_{\mathrm{Spec}(\mathcal{O}_K/\pi^n)} G_n' \xrightarrow{\sim} G_n \times_{G_n''} G_n, \quad (g, g') \mapsto (g, gg').$$

Par (2.2.2), $G'_n$ est syntomique sur $\mathrm{Spec}(\mathcal{O}_K/\pi^n)$, par changement de base, $G_n \times_{\mathrm{Spec}(\mathcal{O}_K/\pi^n)} G'_n$ est donc syntomique sur $G_n$ d'où $\mathrm{pr}_1 : G_n \times_{G''_n} G_n \to G_n$ syntomique par l'isomorphisme précédent. Considérons le diagramme cartésien:

$$\begin{array}{ccc} G_n \times_{G''_n} G_n & \stackrel{\mathrm{pr}_2}{\to} & G_n \\ \mathrm{pr}_1 \downarrow & & \downarrow \\ G_n & \to & G''_n \end{array}$$

où $G_n \to G''_n$ est fidèlement plat, en utilisant (E.G.A. IV.16.9.10,(ii)), (E.G.A. IV.19.1.5,(ii)) et la propriété (3) des morphismes syntomiques, on déduit de la syntomicité de $\mathrm{pr}_1$ celle de $G_n \to G''_n$ (voir aussi [Ma]). Réciproquement, supposons qu'on a une suite exacte $0 \to G' \to G \to G'' \to 0$ dans $(\mathrm{Ab}/\mathcal{O}_K)$ avec $G'$, $G$ et $G''$ provenant de $(p\text{-}\mathrm{Gr}/\mathcal{O}_K)$. Il est facile de voir que le morphisme $G \to G''$ est encore un épimorphisme de faisceaux sur le gros site fppf de $\mathrm{Spec}(\mathcal{O}_K)$, et donc de noyau un schéma en groupes fini et plat $H$ tel qu'on a une suite exacte dans $(p\text{-}\mathrm{Gr}/\mathcal{O}_K)$: $0 \to H \to G \to G'' \to 0$. Par la première partie de la preuve, cette suite est exacte dans $(\mathrm{Ab}/\mathcal{O}_K)$, donc $H = G'$ et la suite $0 \to G' \to G \to G'' \to 0$ est exacte dans $(p\text{-}\mathrm{Gr}/\mathcal{O}_K)$. □

PROPOSITION 2.2.4. *Soit $0 \to G' \to G \to G'' \to 0$ une suite exacte dans $(\mathrm{Ab}/\mathcal{O}_K)$, si $G'$ et $G''$ proviennent de $(p\text{-}\mathrm{Gr}/\mathcal{O}_K)$, il en est de même pour $G$ et la suite est exacte dans $(p\text{-}\mathrm{Gr}/\mathcal{O}_K)$.*

*Preuve.* On a un isomorphisme de faisceaux sur $\mathrm{Sp}f(\mathcal{O}_K)_{\mathrm{syn}}$:

$$G \times_{G''} (G' \times G'') \xrightarrow{\sim} G \times_{G''} G, \ (g, (g', g'')) \mapsto (g, gg').$$

D'autre part, les hypothèses fournissent l'existence d'un recouvrement $Y$ de $G''$ dans $\mathrm{Sp}f(\mathcal{O}_K)_{\mathrm{syn}}$ et d'un morphisme de faisceaux $Y \to G$ (ici, on identifie $Y$ avec le faisceau qu'il représente) tels que le diagramme:

$$\begin{array}{ccc} & & Y \\ & \swarrow & \downarrow \\ G & \to & G'' \end{array}$$

commute. Par changement de base $Y \to G$, on en déduit un isomorphisme de faisceaux sur $\mathrm{Sp}f(\mathcal{O}_K)_{\mathrm{syn}}$:

$$Y \times_{G''} (G' \times G'') \xrightarrow{\sim} Y \times_{G''} G.$$

En particulier, $Y \times_{G''} G$ est représentable dans $\mathrm{Sp}f(\mathcal{O}_K)_{\mathrm{syn}}$ par le schéma formel $Y \times_{G''} (G' \times_{\mathrm{Sp}f(\mathcal{O}_K)} G'')$ qui est alors muni de données de descente provenant de $Y \times_{G''} G$. On en déduit l'existence d'un schéma formel $p$-adique plat sur $G''$ qui représente le faisceau $G$ dans $(\mathrm{Ab}/\mathcal{O}_K)$, et qui est aussi fini sur $G''$ par (E.G.A. IV.2.5.2,(iv)). En procédant comme dans la première partie de la preuve précédente, on voit qu'il est syntomique sur $G''$ parce qu'il l'est sur $Y$ après changement de base, donc c'est un objet de $\mathrm{Sp}f(\mathcal{O}_K)_{\mathrm{syn}}$. Tout ceci montre finalement qu'il est muni d'une structure d'objet de $(p\text{-}\mathrm{Gr}/\mathcal{O}_K)$ et on achève avec (2.2.3). □

*Remarque* 2.2.5. Un morphisme de schémas en groupes finis et plats (sur $\mathcal{O}_K$) qui induit un isomorphisme sur les fibres génériques est une injection de faisceaux dans $(\mathrm{Ab}/\mathcal{O}_K)$, mais *pas forcément* sur le gros site fppf. Considérer par exemple, lorsque $K$ contient les racines $p^{\mathrm{ièmes}}$ de 1, un morphisme non trivial $\mathbf{Z}/p\mathbf{Z} \to \mu_p$.

2.3. *Les faisceaux $\mathcal{O}^{\mathrm{cris}}_{\infty,\pi}$ et $\mathcal{J}^{\mathrm{cris}}_{\infty,\pi}$*. Pour $n \in \mathbf{N}$, on note $E_n = \mathrm{Spec}(S_n)$. Pour tout schéma $X$ sur $\mathrm{Spec}(\mathcal{O}_K/p^n)$, l'épaississement à puissances divisées de (2.1.1): $\mathrm{Spec}(\mathcal{O}_K/p^n) \hookrightarrow E_n$ permet de définir le petit site cristallin classique $(X/E_n)_{\mathrm{cris}}$ et les faisceaux $\mathcal{J}_{X/E_n} \subset \mathcal{O}_{X/E_n}$ dessus ([Be, III.1.1]). On définit des préfaisceaux $\mathcal{J}^{\mathrm{cris}}_{n,\pi} \subset \mathcal{O}^{\mathrm{cris}}_{n,\pi}$ sur $\mathrm{Sp}f(\mathcal{O}_K)_{\mathrm{syn}}$ en posant:

$$\mathcal{O}^{\mathrm{cris}}_{n,\pi}(\mathfrak{X}) = H^0_{\mathrm{cris}}((\mathfrak{X}_n/E_n)_{\mathrm{cris}}, \mathcal{O}_{\mathfrak{X}_n/E_n}),$$
$$\mathcal{J}^{\mathrm{cris}}_{n,\pi}(\mathfrak{X}) = H^0_{\mathrm{cris}}((\mathfrak{X}_n/E_n)_{\mathrm{cris}}, \mathcal{J}_{\mathfrak{X}_n/E_n})$$

où $\mathfrak{X}_n = \mathfrak{X} \times \mathbf{Z}/p^n\mathbf{Z}$. On vérifie comme dans [FM, 1.3] qu'on obtient ainsi des faisceaux sur $\mathrm{Sp}f(\mathcal{O}_K)_{\mathrm{syn}}$. On note $\mathcal{O}_n$ le faisceau structural "modulo $p^n$" sur $\mathrm{Sp}f(\mathcal{O}_K)_{\mathrm{syn}}$, i.e. le faisceau défini par $\mathcal{O}_n(\mathfrak{X}) = H^0(\mathfrak{X}_n, \mathcal{O}_{\mathfrak{X}_n})$.

*Remarque* 2.3.1. Il ne faut pas confondre les faisceaux $\mathcal{O}^{\mathrm{cris}}_{n,\pi}$ ci-dessus avec les faisceaux $\mathcal{O}^{\mathrm{cris}}_n$ introduits dans [FM] qui désignent, eux, simplement la cohomologie cristalline par rapport à la base $\mathrm{Spec}(W_n)$. L'indice $\pi$ dans la notation $\mathcal{O}^{\mathrm{cris}}_{n,\pi}$ est là pour rappeler que ces faisceaux dépendent du plongement $\mathrm{Spec}(\mathcal{O}_K/p^n) \hookrightarrow E_n$ que définit $\pi$.

Le lemme suivant est technique mais crucial:

LEMME 2.3.2. *Soient*

$$\mathfrak{A} = \mathcal{O}_K\{X_1, \ldots, X_r\}/(f_1, \ldots, f_s)$$

*où $(f_1, \ldots, f_s)$ est une suite régulière* (2.2.1),

$$\mathfrak{A}_i = \mathcal{O}_K\{X_0^{p^{-i}}, X_1^{p^{-i}}, \ldots, X_r^{p^{-i}}\}/(X_0 - \pi, f_1, \ldots, f_s),$$

$\mathfrak{A}_\infty = \varinjlim_i \mathfrak{A}_i$ *et*

$$A_\infty = \mathfrak{A}_\infty/p \simeq k[X_0^{p^{-\infty}}, X_1^{p^{-\infty}}, \ldots, X_r^{p^{-\infty}}]/(X_0^e, f_1, \ldots, f_s).$$

*Posons*

$$\mathcal{O}^{\mathrm{cris}}_{n,\pi}(\mathfrak{A}_\infty) = \varinjlim_i \mathcal{O}^{\mathrm{cris}}_{n,\pi}(\mathrm{Sp}f(\mathfrak{A}_i))$$

*et*

$$\mathcal{J}_{n,\pi}^{\mathrm{cris}}(\mathfrak{A}_\infty) = \varinjlim_{i} \mathcal{J}_{n,\pi}^{\mathrm{cris}}(\mathrm{Sp}f(\mathfrak{A}_i)),$$

*alors*:

(1) *Les $W_n$-modules $\mathcal{O}_{n,\pi}^{\mathrm{cris}}(\mathfrak{A}_\infty)$ et $\mathcal{J}_{n,\pi}^{\mathrm{cris}}(\mathfrak{A}_\infty)$ sont plats sur $W_n$ et on a des suites exactes*:

$$0 \to \mathcal{O}_{i,\pi}^{\mathrm{cris}}(\mathfrak{A}_\infty) \xrightarrow{p^n} \mathcal{O}_{n+i,\pi}^{\mathrm{cris}}(\mathfrak{A}_\infty) \to \mathcal{O}_{n,\pi}^{\mathrm{cris}}(\mathfrak{A}_\infty) \to 0$$

*et*

$$0 \to \mathcal{J}_{i,\pi}^{\mathrm{cris}}(\mathfrak{A}_\infty) \xrightarrow{p^n} \mathcal{J}_{n+i,\pi}^{\mathrm{cris}}(\mathfrak{A}_\infty) \to \mathcal{J}_{n,\pi}^{\mathrm{cris}}(\mathfrak{A}_\infty) \to 0.$$

(2) *Si $\psi_1, \ldots, \psi_s$ sont des éléments de $k[X_0^{p^{-\infty}}, X_1^{p^{-\infty}}, \ldots, X_r^{p^{-\infty}}]$ tels que $\psi_j^p = \overline{f}_j$, on a (à une torsion près par le Frobenius sur $k$)*:

$$\mathcal{O}_{1,\pi}^{\mathrm{cris}}(\mathfrak{A}_\infty) \simeq \bigoplus_{(m_0,\ldots,m_{s+1})\in \mathbf{N}^{s+1}} \frac{A_\infty[u]}{u^p - X_0} \gamma_{pm_0}((X_0^{p^{-1}})^e) \gamma_{pm_1}(\psi_1) \cdots$$
$$\cdots \gamma_{pm_s}(\psi_s) \gamma_{pm_{s+1}}(u - X_0^{p^{-1}}).$$

*Preuve*. Elle est très similaire à [Br4, 2.1.2.1]. Soit $(W_n(A_\infty) \otimes_{W_n,(\phi)^n} W_n[u])^{\mathrm{DP}}$ l'enveloppe aux puissances divisées (compatible avec les puissances divisées sur l'idéal $(p)$) par rapport au noyau de la surjection $W_n[u]$-linéaire canonique:

$$W_n(A_\infty) \otimes_{W_n,(\phi)^n} W_n[u] \to \mathfrak{A}_\infty/p^n\mathfrak{A}_\infty$$

qui envoie $(a_0, \ldots, a_{n-1}) \in W_n(A_\infty)$ sur $\hat{a}_o^{p^n} + p\hat{a}_1^{p^{n-1}} + \cdots + p^{n-1}\hat{a}_{n-1}^p$ où $\hat{a}_i$ désigne un relevé quelconque de $a_i$ dans $\mathfrak{A}_\infty/p^n\mathfrak{A}_\infty$ (et où $\mathfrak{A}_\infty/p^n\mathfrak{A}_\infty$ est vu sur $W_n[u]$ via $u \mapsto X_0 = \pi$). Puisque $(W_n(A_\infty) \otimes_{W_n,(\phi)^n} W_n[u])^{\mathrm{DP}}$ est un épaississement particulier de $\mathfrak{A}_\infty/p^n\mathfrak{A}_\infty$ sur $S/p^nS$, on a un morphisme canonique

$$\mathcal{O}_{n,\pi}^{\mathrm{cris}}(\mathfrak{A}_\infty) \to (W_n(A_\infty) \otimes_{W_n,(\phi)^n} W_n[u])^{\mathrm{DP}}.$$

Mais $W_n(A_\infty) \otimes_{W_n,(\phi)^n} W_n[u]$, donc $(W_n(A_\infty) \otimes_{W_n,(\phi)^n} W_n[u])^{\mathrm{DP}}$, s'envoie canoniquement dans tout DP-épaississement de $\mathfrak{A}_\infty/p^n\mathfrak{A}_\infty$ sur $S/p^nS$ par la même application que ci-dessus (en prenant les $\hat{a}_i$ dans l'épaississement), d'où une flèche canonique $(W_n(A_\infty) \otimes_{W_n,(\phi)^n} W_n[u])^{\mathrm{DP}} \to \mathcal{O}_{n,\pi}^{\mathrm{cris}}(\mathfrak{A}_\infty)$ dont il est formel de vérifier qu'elle est l'inverse de la précédente. On a finalement un isomorphisme canonique:

$$\mathcal{O}_{n,\pi}^{\mathrm{cris}}(\mathfrak{A}_\infty) \xrightarrow{\sim} (W_n(A_\infty) \otimes_{W_n,(\phi)^n} W_n[u])^{\mathrm{DP}}.$$

Puisque

$$\mathfrak{A}_\infty/p^n\mathfrak{A}_\infty \simeq W_n[X_0^{p^{-\infty}}, \ldots, X_r^{p^{-\infty}}]/(E(X_0), f_1, \ldots, f_s),$$

on a une surjection évidente

$$W_n[u][X_0^{p^{-\infty}}, \ldots, X_r^{p^{-\infty}}] \to \mathfrak{A}_\infty/p^n\mathfrak{A}_\infty$$

qui envoie $X_i^{p^{-\infty}}$ sur $X_i^{p^{-\infty}}$ et $u$ sur $X_0$, d'où un épaississement

$$(W_n[u][X_0^{p^{-\infty}}, \ldots, X_r^{p^{-\infty}}])^{\mathrm{DP}} \to \mathfrak{A}_\infty/p^n\mathfrak{A}_\infty$$

en prenant les puissances divisées par rapport au noyau, i.e. l'idéal $(u - X_0, E(X_0), f_1, \ldots, f_s)$. Par ce qui précède, on a donc un morphisme canonique:

$$(W_n(A_\infty) \otimes_{W_n,(\phi)^n} W_n[u])^{\mathrm{DP}} \to (W_n[u][X_0^{p^{-\infty}}, \ldots, X_r^{p^{-\infty}}])^{\mathrm{DP}}$$

dont on vérifie qu'il s'agit d'un isomorphisme en procédant comme dans la preuve de [Br4, 2.1.2.1]. La suite $(u - X_0, E(X_0), f_1, \ldots, f_s)$ étant transversalement régulière par rapport à $W_n$, le même argument que dans la preuve de *loc.cit.* utilisant les critères de platitude des enveloppes à puissances divisées s'applique à $(W_n[u][X_0^{p^{-\infty}}, \ldots, X_r^{p^{-\infty}}])^{\mathrm{DP}}$ et donne alors (1). La formule (2) n'est autre que l'explicitation de l'enveloppe aux puissances divisées $(A_\infty \otimes_{k,(\phi)} k[u])^{\mathrm{DP}}$ précédente (cf. [FM, II.1.7] ou la preuve de [Br4, 2.2.2.2]). □

COROLLAIRE 2.3.3. *Les faisceaux $\mathcal{O}_{n,\pi}^{\mathrm{cris}}$ et $\mathcal{J}_{n,\pi}^{\mathrm{cris}}$ sont plats sur $W_n$ et on a des suites exactes dans* $(\mathrm{Ab}/\mathcal{O}_K)$:

$$0 \to \mathcal{O}_{i,\pi}^{\mathrm{cris}} \xrightarrow{p^n} \mathcal{O}_{n+i,\pi}^{\mathrm{cris}} \to \mathcal{O}_{n,\pi}^{\mathrm{cris}} \to 0$$

*et*

$$0 \to \mathcal{J}_{i,\pi}^{\mathrm{cris}} \xrightarrow{p^n} \mathcal{J}_{n+i,\pi}^{\mathrm{cris}} \to \mathcal{J}_{n,\pi}^{\mathrm{cris}} \to 0.$$

De (2.3.2), on déduit également que le complexe $0 \to \mathcal{J}_{1,\pi}^{\mathrm{cris}} \to \mathcal{O}_{1,\pi}^{\mathrm{cris}} \to \mathcal{O}_1 \to 0$ est une suite exacte dans $(\mathrm{Ab}/\mathcal{O}_K)$. Les faisceaux $\mathcal{O}_n$ étant aussi trivialement plats sur $W_n$, en utilisant (2.3.3) et un dévissage élémentaire, on a pour tout $n \in \mathbf{N}^*$ des suites exactes $0 \to \mathcal{J}_{n,\pi}^{\mathrm{cris}} \to \mathcal{O}_{n,\pi}^{\mathrm{cris}} \to \mathcal{O}_n \to 0$ dans $(\mathrm{Ab}/\mathcal{O}_K)$.

Le Frobenius cristallin, qui existe car on a muni la base $E_n$ d'un (relèvement du) Frobenius (2.1.1), induit un opérateur semi-linéaire $\phi : \mathcal{O}_{n,\pi}^{\mathrm{cris}} \to \mathcal{O}_{n,\pi}^{\mathrm{cris}}$ tel que $\phi(\mathcal{J}_{n,\pi}^{\mathrm{cris}}) \subset p\mathcal{O}_{n,\pi}^{\mathrm{cris}}$ (on le vérifie à partir de (2.3.2)). Par le même argument de platitude qu'en [FM, II.2.3] ou [Br4, 2.1.2], on en déduit pour tout $n \in \mathbf{N}^*$ un opérateur $\phi_1 = $ "$\phi/p$": $\mathcal{J}_{n,\pi}^{\mathrm{cris}} \to \mathcal{O}_{n,\pi}^{\mathrm{cris}}$ dans $(\mathrm{Ab}/\mathcal{O}_K)$.

LEMME 2.3.4. *Pour tout $\mathfrak{X}$ de $\mathrm{Sp}f(\mathcal{O}_K)_{\mathrm{syn}}$, on a $\mathrm{Fil}^1 S \mathcal{O}_{n,\pi}^{\mathrm{cris}}(\mathfrak{X}) \subset \mathcal{J}_{n,\pi}^{\mathrm{cris}}(\mathfrak{X})$ et, si $s \in \mathrm{Fil}^1 S$ et $x \in \mathcal{O}_{n,\pi}^{\mathrm{cris}}(\mathfrak{X})$, $\phi_1(sx) = \phi_1(s)\phi(x)$. En particulier, la donnée $(\mathcal{O}_{n,\pi}^{\mathrm{cris}}(\mathfrak{X}), \mathcal{J}_{n,\pi}^{\mathrm{cris}}(\mathfrak{X}), \phi_1)$ est un objet de $'(\mathrm{Mod}/S)$ (2.1.1).*

*Preuve.* Si $s \in \text{Fil}^1 S$, le morphisme $\mathcal{O}_{n,\pi}^{\text{cris}} \xrightarrow{s\text{Id}} \mathcal{O}_{n,\pi}^{\text{cris}}$ se factorise par $\mathcal{J}_{n,\pi}^{\text{cris}}$ et le diagramme:

$$\begin{array}{ccc} \mathcal{O}_{n,\pi}^{\text{cris}} & \xrightarrow{s.Id} & \mathcal{J}_{n,\pi}^{\text{cris}} \\ \phi_1(s)Id \downarrow & & \downarrow \phi_1 \\ \mathcal{O}_{n,\pi}^{\text{cris}} & \xrightarrow{\phi} & \mathcal{O}_{n,\pi}^{\text{cris}} \end{array}$$

commute, d'où le résultat. $\square$

Enfin, on pose $\mathcal{O}_{\infty,\pi}^{\text{cris}} = \varinjlim_n \mathcal{O}_{n,\pi}^{\text{cris}}$ et $\mathcal{J}_{\infty,\pi}^{\text{cris}} = \varinjlim_n \mathcal{J}_{n,\pi}^{\text{cris}}$. Pour tout $n \in \mathbf{N}^*$, on a encore des suites exactes dans $(\text{Ab}/\mathcal{O}_K)$:

$$0 \to \mathcal{O}_{n,\pi}^{\text{cris}} \to \mathcal{O}_{\infty,\pi}^{\text{cris}} \xrightarrow{p^n} \mathcal{O}_{\infty,\pi}^{\text{cris}} \to 0,$$

$$0 \to \mathcal{J}_{n,\pi}^{\text{cris}} \to \mathcal{J}_{\infty,\pi}^{\text{cris}} \xrightarrow{p^n} \mathcal{J}_{\infty,\pi}^{\text{cris}} \to 0$$

et, pour tout objet $\mathfrak{X}$ de $\text{Sp}f(\mathcal{O}_K)_{\text{syn}}$, $(\mathcal{O}_{\infty,\pi}^{\text{cris}}(\mathfrak{X}), \mathcal{J}_{\infty,\pi}^{\text{cris}}(\mathfrak{X}), \phi_1)$ est un objet de $'(\text{Mod}/S)$.

## 3. Schémas en groupes de type $(p,\ldots,p)$

On suppose $p \neq 2$. On construit une anti-équivalence de catégories entre la catégorie des schémas en groupes finis et plats sur $\mathcal{O}_K$ tués par $p$ et $(\text{Mod}/S_1)$ ou $(\text{Mod}/\tilde{S}_1)$.

3.1. *Le schéma en groupes associé à un module annulé par $p$.* A tout objet $\mathcal{M}$ de $(\text{Mod}/S_1)$, on associe un faisceau $\text{Gr}(\mathcal{M})$ sur $\text{Sp}f(\mathcal{O}_K)_{\text{syn}}$ défini par (voir (2.3.4)):

$$\text{Gr}(\mathcal{M})(\mathfrak{X}) = \text{Hom}_{'(\text{Mod}/S_1)}(\mathcal{M}, \mathcal{O}_1^{\text{cris}}(\mathfrak{X})).$$

On obtient ainsi un foncteur contravariant:

$$\text{Gr} : (\text{Mod}/S_1) \to (\text{Ab}/\mathcal{O}_K), \ \mathcal{M} \mapsto \text{Gr}(\mathcal{M}).$$

Nous allons montrer que $\text{Gr}(\mathcal{M})$ est représentable par un schéma fini et syntomique sur $\mathcal{O}_K$.

Jusqu'à la fin de cette partie, on se fixe une racine $p^{\text{ième}}$ $\pi_1$ de $\pi$ dans $\mathcal{O}_{\overline{K}}$, un objet $\mathcal{M}$ de $(\text{Mod}/S_1)$ de rang $d$ et une base adaptée $(e_1,\ldots,e_d)$ de $\mathcal{M}$ (2.1.2.6). Soient $\mathcal{O}_{K_1} = \mathcal{O}_K[\pi_1]$, $(r_1,\ldots,r_d)$ les entiers de $\{0,\ldots,e\}$ tels que $u^{r_i} e_i \in \text{Fil}^1 \tilde{\mathcal{M}}$ et $\mathcal{G}$ l'unique matrice de $\text{GL}_d(S_1)$ telle que:

$$\begin{pmatrix} \phi_1(u^{r_1} e_1) \\ \vdots \\ \phi_1(u^{r_d} e_d) \end{pmatrix} = \mathcal{G} \begin{pmatrix} e_1 \\ \vdots \\ e_d \end{pmatrix}.$$

Soit $\tilde{\mathcal{G}}$ l'image de $\mathcal{G}$ dans $\mathrm{GL}_d(\tilde{S}_1) = \mathrm{GL}_d(k[u]/u^{ep})$, on a un isomorphisme $k$-linéaire $k[u]/u^{ep} \xrightarrow{\sim} \mathcal{O}_{K_1}/p \otimes_{k,(\phi)} k$, $\lambda u^i \mapsto \pi_1^i \otimes \lambda = \lambda^{p^{-1}} \pi_1^i \otimes 1$. On note $\mathcal{G}_\pi = (a_{ij})_{\substack{1 \leq i \leq d \\ 1 \leq j \leq d}}$ un relevé dans $\mathrm{GL}_d(\mathcal{O}_{K_1})$ de l'image de $\tilde{\mathcal{G}}$ dans $\mathrm{GL}_d(\mathcal{O}_{K_1}/p)$. On pose:

$$R_{1,\mathcal{M}} = \frac{\mathcal{O}_{K_1}[X_1, \ldots, X_d]}{X_1^p + \frac{\pi^{e-r_1}}{F(\pi)}(\sum_{j=1}^d a_{1j} X_j), \ldots, X_d^p + \frac{\pi^{e-r_d}}{F(\pi)}(\sum_{j=1}^d a_{dj} X_j)}$$

(rappelons que $u^e - pF(u)$ est le polynôme minimal de $\pi$). Il est clair que $R_{1,\mathcal{M}}$ est fini et syntomique sur $\mathcal{O}_K$ de rang $p^d$. On lui fait subir trois transformations successives:

- on note $''R_\mathcal{M}$ la restriction à la Weil de $\mathcal{O}_{K_1}$ à $\mathcal{O}_K$ de $R_{1,\mathcal{M}}$,
- on note $'R_\mathcal{M}$ le quotient de $''R_\mathcal{M}$ par ses éléments de $p$-torsion,
- on note $R_\mathcal{M}$ la complétion $p$-adique de $'R_\mathcal{M}$.

La $\mathcal{O}_K$-algèbre $R_\mathcal{M}$ est plate et $p$-adiquement complète.

PROPOSITION 3.1.1. *L'anneau $R_\mathcal{M}$ est fini et syntomique sur $\mathcal{O}_K$ de rang $p^d$.*

*Preuve.* La $\mathcal{O}_K$-algèbre $''R_\mathcal{M}$ se décrit explicitement comme suit: soient $S_{i,l} \in \mathcal{O}_K[X_{1,0}, \ldots, X_{1,p-1}, X_{2,0}, \ldots, X_{2,p-1}, \ldots, X_{d,0}, \ldots, X_{d,p-1}]$ tels que pour $i \in \{1, \ldots, d\}$:

$$\left(X_{i,0} + \pi_1 X_{i,1} + \cdots + \pi_1^{p-1} X_{i,p-1}\right)^p$$
$$+ \frac{\pi^{e-r_i}}{F(\pi)} \left(\sum_{j=1}^d a_{ij}(X_{j,0} + \pi_1 X_{j,1} + \cdots + \pi_1^{p-1} X_{j,p-1})\right) = \sum_{l=0}^{p-1} \pi_1^l S_{i,l}$$

dans $\mathcal{O}_{K_1}[X_{1,0}, \ldots, X_{1,p-1}, X_{2,0}, \ldots, X_{2,p-1}, \ldots, X_{d,0}, \ldots, X_{d,p-1}]$, on a $''R_\mathcal{M} = \mathcal{O}_K[X_{1,0}, \ldots, X_{d,p-1}]/(S_{i,l})_{\substack{1 \leq i \leq d \\ 0 \leq l \leq p-1}}$. Notons:

$$S_{\cdot,l} = \begin{pmatrix} S_{1,l} \\ \vdots \\ S_{d,l} \end{pmatrix}, \quad X_{\cdot,l} = \begin{pmatrix} X_{1,l} \\ \vdots \\ X_{d,l} \end{pmatrix}, \quad Y_{\cdot,l} = \begin{pmatrix} X_{1,0}^{p-1} X_{1,l} \\ \vdots \\ X_{d,0}^{p-1} X_{d,l} \end{pmatrix},$$

et $\mathcal{G}_\pi = \mathcal{G}_0 + \pi_1 \mathcal{G}_1 + \cdots + \pi_1^{p-1} \mathcal{G}_{p-1}$ où $\mathcal{G}_0 \in \mathrm{GL}_d(\mathcal{O}_K)$ et $\mathcal{G}_l \in M_d(\mathcal{O}_K)$ si $l \geq 1$, un examen des équations définissant les $S_{\cdot,l}$ donne:

$$S_{\cdot,0} = Y_{\cdot,0} + \mathrm{Diag}(\frac{\pi^{e-r_i}}{F(\pi)}) \cdot \mathcal{G}_0 \cdot X_{\cdot,0} + \pi f_0(X_{\cdot,0}, \ldots, X_{\cdot,p-1}),$$

$$S_{.,l} = pg_l(X_{.,0},\ldots,X_{.,p-1}) +$$
$$\operatorname{Diag}(\frac{\pi^{e-r_i}}{F(\pi)}) \cdot \Big(\mathcal{G}_0.X_{.,l} + \mathcal{G}_1.X_{.,l-1} + \cdots + \mathcal{G}_l.X_{.,0} + \pi f_l(X_{.,0},\ldots,X_{.,p-1})\Big)$$

où

$$g_l(X_{.,0},\ldots,X_{.,p-1}) = Y_{.,l} + \pi h_l(X_{.,0},\ldots,X_{.,p-1})$$

(si $l \geq 1$) avec $h_l(X_{.,0},\ldots,X_{.,p-1})$ et $f_l(X_{.,0},\ldots,X_{.,p-1})$ dans $(\mathcal{O}_K[X_{i,j}])^d$. Pour $l \geq 1$, cela donne (quitte à changer $h_l$ et $f_l$):

$$\begin{aligned} S_{.,l} &= \operatorname{Diag}(\frac{\pi^{e-r_i}}{F(\pi)}) \cdot \mathcal{G}_0 \cdot (Id + \mathcal{G}_0^{-1} \cdot \operatorname{Diag}(\pi^{r_i} X_{i,0}^{p-1})) \cdot X_{.,l} \\ &+ \operatorname{Diag}(\frac{\pi^{e-r_i}}{F(\pi)}) \cdot \mathcal{G}_0 \cdot (h_l(X_{.,0},\ldots,X_{.,l-1}) + \pi f_l(X_{.,0},\ldots,X_{.,p-1})). \end{aligned}$$

Soit $\mathcal{H} = \operatorname{Id} + \mathcal{G}_0^{-1} \cdot \operatorname{Diag}(\pi^{r_i} X_{i,0}^{p-1})$, je dis que l'image de $\mathcal{H}$ dans $M_d(R_\mathcal{M})$ tombe dans $\operatorname{GL}_d(R_\mathcal{M})$. Les coefficients de $\Big(\mathcal{G}_0^{-1} \cdot \operatorname{Diag}(\pi^{r_i} X_{i,0}^{p-1})\Big)^n$ sont des polynômes homogènes de degré $n$ à coefficients dans $\mathcal{O}_K$ en les variables $(Z_1,\ldots,Z_d)$ où $Z_i = \pi^{r_i} X_{i,0}^{p-1}$, donc, comme $R_\mathcal{M}$ est complet, il suffit de montrer que si $n_1 + \cdots + n_d = n$, $Z_1^{n_1} \cdots Z_d^{n_d}$ tend ($p$-adiquement) vers 0 dans $R_\mathcal{M}$ lorsque $n \to +\infty$, ou encore $Z_i^n \to 0$ (dans $R_\mathcal{M}$) si $n \to +\infty$. C'est clair si $r_i \geq 1$. Si $r_i = 0$, un examen de $S_{.,0}$ montre que l'image de $X_{i,0}^p$ dans $R_\mathcal{M}$ est dans $\pi R_\mathcal{M}$, donc $Z_i^n = X_{i,0}^{n(p-1)}$ tend vers 0 dans $R_\mathcal{M}$. On a donc pour $l \geq 1$ des équations dans $R_\mathcal{M}$ (car $R_\mathcal{M}$ est sans $\pi$-torsion et $\mathcal{G}_0$ est inversible):

$$X_{.,l} = -\mathcal{H}^{-1} \cdot (h_l(X_{.,0},\ldots,X_{.,l-1}) + \pi f_l(X_{.,0},\ldots,X_{.,p-1}))$$

qu'on peut récrire dans $R_\mathcal{M}$, en injectant successivement la formule donnant $X_{.,l-1}$ dans celle donnant $X_{.,l}$ et quitte à changer $f_l$: $X_{.,l} = -\mathcal{H}^{-1} \cdot (\tilde{h}_l(X_{.,0}) + \pi f_l(X_{.,0},\ldots,X_{.,p-1}))$, les $\tilde{h}_l$ étant des éléments de $(R_\mathcal{M})^d$ ne dépendant que de l'image des $X_{i,0}$. En réinjectant la formule des $X_{.,l}$ ($l \geq 1$) dans les $f_l$ et par un procédé inductif élémentaire, on voit que les $X_{.,l}$ sont déterminés dans $R_\mathcal{M}$ par les $X_{.,0}$. Plus précisément, si $I = \{i \in \{1,\ldots,d\}/r_i \geq 1\}$, $J = \{1,\ldots,d\} \setminus I$ et $R$ est l'image de $\mathcal{O}_K\{X_{i,0}, i \in I\}[[X_{j,0}, j \in J]]$ dans $R_\mathcal{M}$ (rappelons que l'image de $X_{j,0}^p$ est dans $\pi R_\mathcal{M}$ si $j \in J$), on voit que l'image de $X_{i,l}$ dans $R_\mathcal{M}$ tombe dans $R$. On en déduit des équations dans $R_\mathcal{M}$ de la forme:

$$X_{i,0}^p = -\frac{\pi^{e-r_i}}{F(\pi)}(\sum_{j=1}^d a_{ij}^0 X_{j,0}) + \pi \tilde{f}_i(X_{1,0},\ldots,X_{d,0})$$

où $\tilde{f}_i(X_{1,0},\ldots,X_{d,0}) \in R$ et $(a_{ij}^0)_{\substack{1 \leq i \leq d \\ 1 \leq j \leq d}} = \mathcal{G}_0 \in \operatorname{GL}_d(\mathcal{O}_K)$. Par un nouveau procédé inductif (en réinjectant la formule des $X_{j,0}$ dans $\tilde{f}_i$), on a finalement

des équations:

$$X_{i,0}^p = -\frac{\pi^{e-r_i}}{F(\pi)}(\sum_{j=1}^d a_{ij}^0 X_{j,0}) + \pi f_i(X_{1,0}, \ldots, X_{d,0})$$

pour certains $f_i(X_{1,0}, \ldots, X_{d,0}) \in \bigoplus_{n_i \leq p-1} \mathcal{O}_K X_{1,0}^{n_1} \cdots X_{d,0}^{n_d}$. Autrement dit, on a:

$$R_\mathcal{M} \simeq \frac{\mathcal{O}_K[X_{1,0}, \ldots, X_{d,0}]}{\left(X_{1,0}^p + \frac{\pi^{e-r_1}}{F(\pi)}(\sum_{j=1}^d a_{1j}^0 X_{j,0}) + \pi f_1, \ldots, X_{d,0}^p + \frac{\pi^{e-r_d}}{F(\pi)}(\sum_{j=1}^d a_{dj}^0 X_{j,0}) + \pi f_d\right)}$$

où les $f_i$ sont comme ci-dessus. Il est alors clair que $R_\mathcal{M}$ est fini et plat de rang $p^d$ sur $\mathcal{O}_K$ et que les $d$ équations du quotient forment une suite régulière (regarder modulo $\pi$). □

La preuve précédente permet une description calculatoire un peu fastidieuse de $R_\mathcal{M}$. Mais il existe un cas où $R_\mathcal{M}$ se décrit facilement:

PROPOSITION 3.1.2. *Supposons qu'on puisse choisir $\mathcal{G}_\pi$ dans $\mathrm{GL}_d(\mathcal{O}_K)$ (par exemple si $\tilde{\mathcal{G}} \in \mathrm{GL}_d(k[u^p]/u^{ep})$), alors:*

$$R_\mathcal{M} \simeq \frac{\mathcal{O}_K[X_1, \ldots, X_d]}{\left(X_1^p + \frac{\pi^{e-r_1}}{F(\pi)}(\sum_{j=1}^d a_{1j} X_j), \ldots, X_d^p + \frac{\pi^{e-r_d}}{F(\pi)}(\sum_{j=1}^d a_{dj} X_j)\right)}$$

*Preuve.* Avec les notations de la preuve précédente, on a $\mathcal{G}_l = 0$ si $l \geq 1$ (ou $\mathcal{G}_\pi = \mathcal{G}_0$), et on vérifie qu'on aboutit à des équations dans $R_\mathcal{M}$ pour $l \geq 1$ de la forme $X_{\cdot,l} = -\mathcal{H}^{-1} \cdot (\pi f_l(X_{\cdot,0}, \ldots, X_{\cdot,p-1}))$ où $f_l \in (\mathcal{O}_K[X_{i,j}])^d$ sans termes constants et sans monômes en les $X_{i,0}$ seulement. Un procédé inductif montre alors que les $X_{i,l}$ sont tous nuls dans $R_\mathcal{M}$ dès que $l \geq 1$. On en déduit facilement le résultat en posant $X_i = X_{i,0}$. □

*Remarque* 3.1.3. Si $e = 1$, on sait que $\mathcal{M}$ provient toujours de la catégorie de Fontaine-Laffaille $\underline{MF}_k^{f,1}$ (voir [Br3, 2.4.1] et [Br4, 4.4.1]). On en déduit que dans le cas $e = 1$, on peut toujours se ramener à $\tilde{\mathcal{G}} \in \mathrm{GL}_d(k)$ et (3.1.2) donne dans ce cas une description directe de $R_\mathcal{M}$.

PROPOSITION 3.1.4. *Pour tout $\mathfrak{X}$ dans $\mathrm{Spf}(\mathcal{O}_K)_{\mathrm{syn}}$, on a:*

$$\mathrm{Hom}_{\mathrm{Spf}(\mathcal{O}_K)}\Big(\mathfrak{X}, \mathrm{Spf}(R_\mathcal{M})\Big)$$
$$= \mathrm{Hom}_{\mathrm{Spf}(\mathcal{O}_{K_1})}\Big(\mathrm{Spf}(\mathcal{O}_{K_1}) \times_{\mathrm{Spf}(\mathcal{O}_K)} \mathfrak{X}, \mathrm{Spf}(R_{1,\mathcal{M}})\Big).$$

*Preuve.* On a:

$$\mathrm{Hom}_{\mathrm{Sp}f(\mathcal{O}_K)}\Big(\mathfrak{X}, \mathrm{Sp}f(R_\mathcal{M})\Big) = \mathrm{Hom}_{\mathcal{O}_K}\Big(R_\mathcal{M}, \Gamma(\mathfrak{X}, \mathcal{O}_\mathfrak{X})\Big),$$

$$\mathrm{Hom}_{\mathrm{Sp}f(\mathcal{O}_{K_1})}\Big(\mathrm{Sp}f(\mathcal{O}_{K_1}) \times \mathfrak{X}, \mathrm{Sp}f(R_{1,\mathcal{M}})\Big)$$
$$= \mathrm{Hom}_{\mathcal{O}_{K_1}}\Big(R_{1,\mathcal{M}}, \mathcal{O}_{K_1} \otimes_{\mathcal{O}_K} \Gamma(\mathfrak{X}, \mathcal{O}_\mathfrak{X})\Big).$$

Mais:

$$\begin{aligned}
\mathrm{Hom}_{\mathcal{O}_{K_1}}\Big(R_{1,\mathcal{M}}, \mathcal{O}_{K_1} \otimes_{\mathcal{O}_K} \Gamma(\mathfrak{X}, \mathcal{O}_\mathfrak{X})\Big) &= \mathrm{Hom}_{\mathcal{O}_K}\Big({''}R_\mathcal{M}, \Gamma(\mathfrak{X}, \mathcal{O}_\mathfrak{X})\Big) \\
&= \mathrm{Hom}_{\mathcal{O}_K}\Big({'}R_\mathcal{M}, \Gamma(\mathfrak{X}, \mathcal{O}_\mathfrak{X})\Big) \\
&= \mathrm{Hom}_{\mathcal{O}_K}\Big(R_\mathcal{M}, \Gamma(\mathfrak{X}, \mathcal{O}_\mathfrak{X})\Big)
\end{aligned}$$

en utilisant la définition de ${''}R_\mathcal{M}$ et le fait que $\Gamma(\mathfrak{X}, \mathcal{O}_\mathfrak{X})$ est plat sur $\mathcal{O}_K$ et $p$-adiquement complet. $\square$

PROPOSITION 3.1.5. *Soit $\mathrm{Sp}f(R_\mathcal{M})^\sim$ le faisceau sur $\mathrm{Sp}f(\mathcal{O}_K)_\mathrm{syn}$ représenté par $\mathrm{Sp}f(R_\mathcal{M})$, il existe un morphisme canonique non nul de faisceaux sur $\mathrm{Sp}f(\mathcal{O}_K)_\mathrm{syn}$: $\mathrm{Sp}f(R_\mathcal{M})^\sim \to \mathrm{Gr}(\mathcal{M})$.*

*Preuve.* Notons pour simplifier $\mathcal{O}_{n,\pi}^\mathrm{cris}(\mathfrak{A})$ au lieu de $\mathcal{O}_{n,\pi}^\mathrm{cris}(\mathrm{Sp}f(\mathfrak{A}))$ si $\mathfrak{A}$ est une $\mathcal{O}_K$-algèbre formelle syntomique. Rappelons (voir 2.3.2) que le faisceau $\mathcal{O}_{n,\pi}^\mathrm{cris}$ s'identifie au faisceau associé au préfaisceau

$$\mathfrak{A} \mapsto (W_n(\mathfrak{A}/p) \otimes_{W_n,(\phi)^n} W_n[u])^\mathrm{DP}$$

où les puissances divisées sont prises par rapport au noyau de l'unique flèche

$$W_n(\mathfrak{A}/p) \otimes_{W_n,(\phi)^n} W_n[u] \to \mathfrak{A}/p^n$$

qui envoie $u$ sur l'image de $\pi$ et $(\overline{\mathfrak{a}}_0, \ldots, \overline{\mathfrak{a}}_{n-1})$ sur $\mathfrak{a}_0^{p^n} + p\mathfrak{a}_1^{p^{n-1}} + \cdots + p^{n-1}\mathfrak{a}_{n-1}^p$ (les $\mathfrak{a}_i$ relevant les $\overline{\mathfrak{a}}_i$ dans $\mathfrak{A}/p^n$), et on a donc une flèche canonique

$$W_n(\mathfrak{A}/p) \otimes_{W_n,(\phi)^n} W_n[u] \to \mathcal{O}_{n,\pi}^\mathrm{cris}(\mathfrak{A}).$$

Par recollement, il suffit de construire pour tout $\mathfrak{A}$ comme précédemment une application canonique et fonctorielle

$$\mathrm{Hom}_{\mathcal{O}_K}(R_\mathcal{M}, \mathfrak{A}) \to \mathrm{Hom}'_{(\mathrm{Mod}/S_1)}(\mathcal{M}, \mathcal{O}_{1,\pi}^\mathrm{cris}(\mathfrak{A}))$$

ou encore

$$\mathrm{Hom}_{\mathcal{O}_K}({''}R_\mathcal{M}, \mathfrak{A}) \to \mathrm{Hom}'_{(\mathrm{Mod}/S_1)}(\mathcal{M}, \mathcal{O}_{1,\pi}^\mathrm{cris}(\mathfrak{A})).$$

Soit $f \in \mathrm{Hom}_{\mathcal{O}_K}({''R}_{\mathcal{M}}, \mathfrak{A})$ et $\mathfrak{a}_{i,j} = f(X_{i,j})$ (cf. preuve de 3.1.1), à $f$ on associe dans un premier temps $\tilde{g} : \mathcal{M} \to \mathcal{O}^{\mathrm{cris}}_{1,\pi}(\mathfrak{A})$ définie par

$$\tilde{g}(e_i) = \overline{\mathfrak{a}}_{i,0} \otimes 1 + \overline{\mathfrak{a}}_{i,1} \otimes u + \cdots + \overline{\mathfrak{a}}_{i,p-1} \otimes u^{p-1}$$

(avec des notations évidentes) qu'on écrit pour simplifier $\overline{\mathfrak{a}}_{i,0} + u\overline{\mathfrak{a}}_{i,1} + \cdots + u^{p-1}\overline{\mathfrak{a}}_{i,p-1}$. On vérifie que

$$u^{r_i}(\overline{\mathfrak{a}}_{i,0} + u\overline{\mathfrak{a}}_{i,1} + \cdots + u^{p-1}\overline{\mathfrak{a}}_{i,p-1}) \in \mathcal{J}^{\mathrm{cris}}_{1,\pi}(\mathfrak{A})$$

(car son image est nulle dans $\mathfrak{A}/p$): calculons son image par $\phi_1$. Soient $T_i \in \mathbf{Z}[u^p][X_0, \ldots, X_{p-1}]$ les polynômes tels que

$$(X_0 + uX_1 + \cdots + u^{p-1}X_{p-1})^p = T_0 + puT_1 + \cdots + pu^{p-1}T_{p-1},$$

$\tilde{T}_0 \in \mathbf{Z}[u][X_0, \ldots, X_{p-1}]$ le polynôme obtenu en remplaçant les $u^p$ par des $u$ dans $T_0$ et $C_{i,j}$ les polynômes sur $\mathcal{O}_K$ (de degré $\leq 1$) en les $X_{i,j}$ définis par (cf. preuve de 3.1.1):

$$(*) \qquad (\mathcal{G}_0 + \pi_1 \mathcal{G}_1 + \cdots + \pi_1^{p-1} \mathcal{G}_{p-1}) \cdot (X_{.,0} + \pi_1 X_{.,1} + \cdots + \pi_1^{p-1} X_{.,p-1})$$

$$= \begin{pmatrix} C_{1,0} \\ \vdots \\ C_{d,0} \end{pmatrix} + \pi_1 \begin{pmatrix} C_{1,1} \\ \vdots \\ C_{d,1} \end{pmatrix} + \cdots + \pi_1^{p-1} \begin{pmatrix} C_{1,p-1} \\ \vdots \\ C_{d,p-1} \end{pmatrix}.$$

Quitte à prendre un recouvrement syntomique de $\mathfrak{A}$, on vérifie, en utilisant que $\tilde{T}_0 \equiv X_0^p + uX_1^p + \cdots + u^{p-1}X_{p-1}^p \bmod p$, que pour tout $i \in \{1, \ldots, d\}$ l'élément:

$$u^{r_i}\tilde{T}_0([\overline{\mathfrak{a}}_{i,0}^{p^{-2}}], \ldots, [\overline{\mathfrak{a}}_{i,p-1}^{p^{-2}}]) + p[C_{i,0}(\overline{\mathfrak{a}}_{m,m'})^{p^{-2}}]$$

de $\mathcal{O}^{\mathrm{cris}}_{2,\pi}(\mathfrak{A})$ relève

$$u^{r_i}(\overline{\mathfrak{a}}_{i,0} + u\overline{\mathfrak{a}}_{i,1} + \cdots + u^{p-1}\overline{\mathfrak{a}}_{i,p-1})$$

et a pour image 0 dans $\mathfrak{A}/p^2\mathfrak{A}$ ($[.]$ désignant le représentant de Teichmüller). Son image par $\phi_1 = ``\frac{\phi}{p}$'' dans $\mathcal{O}^{\mathrm{cris}}_{1,\pi}(\mathfrak{A})$ donne:

$$\frac{(u^{r_i}(\overline{\mathfrak{a}}_{i,0} + u\overline{\mathfrak{a}}_{i,1} + \cdots + u^{p-1}\overline{\mathfrak{a}}_{i,p-1}))^p}{p} - u^{pr_i+1}T_1(\overline{\mathfrak{a}}_{i,0}, \ldots, \overline{\mathfrak{a}}_{i,p-1}) - \cdots$$

$$\cdots - u^{pr_i+p-1}T_{p-1}(\overline{\mathfrak{a}}_{i,0}, \ldots, \overline{\mathfrak{a}}_{i,p-1}) + C_{i,0}(\overline{\mathfrak{a}}_{m,m'}).$$

Notons ${'T}_i \in \mathcal{O}_K[X_0, \ldots, X_{p-1}]$ les polynômes obtenus en remplaçant $u^p$ par $\pi$ dans les $T_i$, par définition les $\mathfrak{a}_{i,j}$ vérifient les égalités pour $1 \leq l \leq p-1$: $\pi^{r_i}{'T}_l(\mathfrak{a}_{i,0}, \ldots, \mathfrak{a}_{i,p-1}) = -C_{i,l}(\mathfrak{a}_{m,m'})$. Mais, dans $\mathcal{O}^{\mathrm{cris}}_{1,\pi}(\mathfrak{A})$, on a $u^p = \overline{\pi}$ (car, avec les notations de (2.3.2), $u - X_0^{p^{-1}}$ est muni de puissances divisées dans $\varinjlim_i \mathcal{O}^{\mathrm{cris}}_{1,\pi}(\mathfrak{A}_i))$ de sorte que:

$$u^{pr_i+l}T_l(\overline{\mathfrak{a}}_{i,0}, \ldots, \overline{\mathfrak{a}}_{i,p-1}) = u^l(\overline{\pi}^{r_i}{'T}_l(\overline{\mathfrak{a}}_{i,0}, \ldots, \overline{\mathfrak{a}}_{i,p-1})) = -u^l C_{i,l}(\overline{\mathfrak{a}}_{m,m'}).$$

Comme on a une égalité analogue à $(*)$ $(\mathrm{mod}\, p)$ en remplaçant $\mathcal{G}_\pi$ par $\tilde{\mathcal{G}} \in \mathrm{GL}_d(\tilde{S}_1)$, $\pi_1$ par $u$ et $\pi$ par $u^p$ dans les coefficients des $C_{i,j}$ $\mathrm{mod}\, p$, on déduit de ce qui précède que $\phi_1(u^{r_i}\tilde{g}(e_i)) = \gamma_p(u^{r_i}\tilde{g}(e_i)) + \sum_{j=1}^{d}\tilde{a}_{ij}\tilde{g}(e_j)$ où $(\tilde{a}_{ij}) = \tilde{\mathcal{G}}$.

Finalement, il existe $c_1,\ldots,c_d$ dans $\mathcal{J}^{cris,[p]}_{1,\pi}(\mathfrak{A})$ (en fait dans l'idéal de $\mathcal{O}^{\mathrm{cris}}_{1,\pi}(\mathfrak{A})$ engendré par les $\gamma_p(x)$, $x \in \mathcal{J}^{\mathrm{cris}}_{1,\pi}(\mathfrak{A})$) tels que:

$$\begin{pmatrix} \phi_1(u^{r_1}\tilde{g}(e_1)) \\ \vdots \\ \phi_1(u^{r_d}\tilde{g}(e_d)) \end{pmatrix} = \mathcal{G} \begin{pmatrix} \tilde{g}(e_1) \\ \vdots \\ \tilde{g}(e_d) \end{pmatrix} + \begin{pmatrix} c_1 \\ \vdots \\ c_d \end{pmatrix}.$$

Et il n'y a qu'une façon d'en déduire un élément $g$ de $\mathrm{Hom}_{'(\mathrm{Mod}/S_1)}(\mathcal{M},\mathcal{O}^{\mathrm{cris}}_{1,\pi}(\mathfrak{A}))$, c'est de poser:

$$\begin{pmatrix} g(e_1) \\ \vdots \\ g(e_d) \end{pmatrix} = \begin{pmatrix} \tilde{g}(e_1) \\ \vdots \\ \tilde{g}(e_d) \end{pmatrix} + \mathcal{G}^{-1} \begin{pmatrix} c_1 \\ \vdots \\ c_d \end{pmatrix}. \qquad \square$$

Soient une $\mathcal{O}_K$-algèbre syntomique formelle

$$\mathfrak{A} = \mathcal{O}_K\{X_1,\ldots,X_r\}/(f_1,\ldots,f_s)$$

avec $(f_1,\ldots,f_s)$ transversalement régulière,

$$\mathfrak{A}_i = \mathcal{O}_K\{X_0^{p^{-i}}, X_1^{p^{-i}},\ldots,X_r^{p^{-i}}\}/(X_0 - \pi, f_1,\ldots,f_s)$$

comme en (2.3.2),

$$\mathfrak{B}_i = \mathcal{O}_{K_1} \otimes_{\mathcal{O}_K} \mathfrak{A}_i, \qquad \mathfrak{A}_\infty = \varinjlim_i \mathfrak{A}_i$$

et

$$\mathfrak{B}_\infty = \mathcal{O}_{K_1} \otimes_{\mathcal{O}_K} \mathfrak{A}_\infty.$$

On pose

$$\mathrm{Fil}^1\mathfrak{B}_i \;=\; \{x \in \mathfrak{B}_i \;/\; x^p \in p\mathfrak{B}_i\}$$

et

$$\mathrm{Fil}^1\mathfrak{B}_\infty \;=\; \{x \in \mathfrak{B}_\infty \;/\; x^p \in p\mathfrak{B}_\infty\};$$

ce sont des idéaux de $\mathfrak{B}_i$ et $\mathfrak{B}_\infty$ respectivement et on note $\mathrm{Fil}^m\mathfrak{B}_i = (\mathrm{Fil}^1\mathfrak{B}_i)^m$, $\mathrm{Fil}^m\mathfrak{B}_\infty = (\mathrm{Fil}^1\mathfrak{B}_\infty)^m$ $(m \geq 1)$. On définit $\phi_1 : \mathrm{Fil}^1\mathfrak{B}_i \to \mathfrak{B}_i$, $x \mapsto x^p/(-p) = \phi_1(x)$ (resp. avec $\mathfrak{B}_\infty$). Si $x \in \mathrm{Fil}^1\mathfrak{B}_i$ et $\delta \in \mathrm{Fil}^m\mathfrak{B}_i$, un calcul simple montre que pour $m \geq 2$, $\phi_1(x+\delta) - \phi_1(x) \in \mathrm{Fil}^m\mathfrak{B}_i$, de sorte qu'on peut encore définir $\phi_1 : \mathrm{Fil}^1\mathfrak{B}_i/\mathrm{Fil}^m\mathfrak{B}_i \to \mathfrak{B}_i/\mathrm{Fil}^m\mathfrak{B}_i$ pour $m \geq 2$. En particulier, $\mathfrak{B}_i/\mathrm{Fil}^p\mathfrak{B}_i$ et $\mathfrak{B}_\infty/\mathrm{Fil}^p\mathfrak{B}_\infty$ étant tués par $p$, on peut les munir d'une structure d'objet de $'(\text{Mod}/\tilde{S}_1)$ en posant $\text{Fil}^1(\mathfrak{B}_i/\text{Fil}^p\mathfrak{B}_i) = \text{Fil}^1\mathfrak{B}_i/\text{Fil}^p\mathfrak{B}_i$, $\text{Fil}^1(\mathfrak{B}_\infty/\text{Fil}^p\mathfrak{B}_\infty) = \text{Fil}^1\mathfrak{B}_\infty/\text{Fil}^p\mathfrak{B}_\infty$, $\phi_1$ comme ci-dessus et en définissant la structure de $\tilde{S}_1$-module par $(\lambda u^i)x = \lambda^{p^{-1}}(\pi_1^i \otimes 1)x$, $\lambda \in k$, $x \in \mathfrak{B}_i/\text{Fil}^p\mathfrak{B}_i$ ou $\mathfrak{B}_\infty/\text{Fil}^p\mathfrak{B}_\infty$ (les conditions en (2.1.2) sont bien vérifiées).

LEMME 3.1.6. *On a un isomorphisme*:
$$\text{Hom}_{'(\text{Mod}/S_1)}(\mathcal{M}, \mathcal{O}_{1,\pi}^{\text{cris}}(\mathfrak{A}_\infty)) \xrightarrow{\sim} \text{Hom}_{'(\text{Mod}/\tilde{S}_1)}(\tilde{\mathcal{M}}, \mathfrak{B}_\infty/\text{Fil}^p\mathfrak{B}_\infty)$$
*où $\tilde{\mathcal{M}}$ est l'objet associé à $\mathcal{M}$ en (2.1.2.2) et $\mathcal{O}_{1,\pi}^{\text{cris}}(\mathfrak{A}_\infty)$ est comme en (2.3.2).*

*Preuve.* Soit $A_\infty = \mathfrak{A}_\infty/p\mathfrak{A}_\infty$, on a un isomorphisme $A_\infty$-linéaire
$$A_\infty[u]/(u^p - \overline{\pi}) \xrightarrow{\sim} \mathfrak{B}_\infty/p\mathfrak{B}_\infty, \quad u \mapsto \overline{\pi}_1 \otimes 1.$$
De (2.3.2), on en déduit un isomorphisme de $\tilde{S}_1$-modules:
$$\mathcal{O}_{1,\pi}^{\text{cris}}(\mathfrak{A}_\infty)/\mathcal{J}_{1,\pi}^{cris,[p]}(\mathfrak{A}_\infty) \xrightarrow{\sim} \mathfrak{B}_\infty/\text{Fil}^p\mathfrak{B}_\infty$$
et c'est un exercice de voir qu'il s'agit d'un isomorphisme dans $'(\text{Mod}/\tilde{S}_1)$. Il suffit donc de montrer que la surjection
$$\mathcal{O}_{1,\pi}^{\text{cris}}(\mathfrak{A}_\infty) \to \mathcal{O}_{1,\pi}^{\text{cris}}(\mathfrak{A}_\infty)/\mathcal{J}_{1,\pi}^{cris,[p]}(\mathfrak{A}_\infty)$$
induit un isomorphisme
$$\text{Hom}_{'(\text{Mod}/S_1)}(\mathcal{M}, \mathcal{O}_{1,\pi}^{\text{cris}}(\mathfrak{A}_\infty)) \xrightarrow{\sim} \text{Hom}_{'(\text{Mod}/\tilde{S}_1)}(\tilde{\mathcal{M}}, \mathcal{O}_{1,\pi}^{\text{cris}}(\mathfrak{A}_\infty)/\mathcal{J}_{1,\pi}^{cris,[p]}(\mathfrak{A}_\infty)).$$
Mais soit
$$f \in \text{Hom}_{'(\text{Mod}/\tilde{S}_1)}(\tilde{\mathcal{M}}, \mathcal{O}_{1,\pi}^{\text{cris}}(\mathfrak{A}_\infty)/\mathcal{J}_{1,\pi}^{cris,[p]}(\mathfrak{A}_\infty))$$
et $\widehat{f(e_i)}$ des relevés de $f(e_i)$ dans $\mathcal{O}_{1,\pi}^{\text{cris}}(\mathfrak{A}_\infty)$, il est clair que:
$$\begin{pmatrix} \phi_1(u^{r_1}\widehat{f(e_1)}) \\ \vdots \\ \phi_1(u^{r_d}\widehat{f(e_d)}) \end{pmatrix} - \mathcal{G} \begin{pmatrix} \widehat{f(e_1)} \\ \vdots \\ \widehat{f(e_d)} \end{pmatrix} \in \left(\mathcal{J}_{1,\pi}^{cris,[p]}(\mathfrak{A}_\infty)\right)^d$$
d'où le résultat en raisonnant comme à la fin de la preuve de (3.1.5). $\square$

LEMME 3.1.7. *On a un isomorphisme:*
$$\text{Hom}_{\mathcal{O}_{K_1}}(R_{1,\mathcal{M}}, \mathfrak{B}_\infty) \xrightarrow{\sim} \text{Hom}_{'(\text{Mod}/\tilde{S}_1)}(\tilde{\mathcal{M}}, \mathfrak{B}_\infty/\text{Fil}^p\mathfrak{B}_\infty).$$

*Preuve.* Soit $f \in \text{Hom}_{\mathcal{O}_{K_1}}(R_{1,\mathcal{M}}, \mathfrak{B}_\infty)$ et $\mathfrak{b}_i = f(X_i)$ (cf. description de $R_{1,\mathcal{M}}$), vu la structure d'objet de $'(\text{Mod}/\tilde{S}_1)$ de $\mathfrak{B}_\infty/\text{Fil}^p\mathfrak{B}_\infty$, il est presque évident que $g : \tilde{\mathcal{M}} \to \mathfrak{B}_\infty/\text{Fil}^p\mathfrak{B}_\infty$ définie par $g(e_i) = \overline{\mathfrak{b}}_i$ (avec des notations évidentes) est un morphisme de $'(\text{Mod}/\tilde{S}_1)$. Ceci définit une flèche:
$$\text{Hom}_{\mathcal{O}_{K_1}}(R_{1,\mathcal{M}}, \mathfrak{B}_\infty) \to \text{Hom}_{'(\text{Mod}/\tilde{S}_1)}(\tilde{\mathcal{M}}, \mathfrak{B}_\infty/\text{Fil}^p\mathfrak{B}_\infty).$$

Montrons qu'il s'agit d'un isomorphisme. Comme $R_{1,\mathcal{M}}$ est de type fini sur $\mathcal{O}_{K_1}$ et $\tilde{\mathcal{M}}$ de type fini sur $\tilde{S}_1$, on a

$$\operatorname{Hom}_{\mathcal{O}_{K_1}}(R_{1,\mathcal{M}}, \mathfrak{B}_\infty) \simeq \varinjlim_i \operatorname{Hom}_{\mathcal{O}_{K_1}}(R_{1,\mathcal{M}}, \mathfrak{B}_i),$$

$$\operatorname{Hom}'_{(\operatorname{Mod}/\tilde{S}_1)}(\tilde{\mathcal{M}}, \mathfrak{B}_\infty/\operatorname{Fil}^p\mathfrak{B}_\infty) \simeq \varinjlim_i \operatorname{Hom}'_{(\operatorname{Mod}/\tilde{S}_1)}(\tilde{\mathcal{M}}, \mathfrak{B}_i/\operatorname{Fil}^p\mathfrak{B}_i)$$

et il suffit de montrer que la flèche (construite de la même manière)

$$\operatorname{Hom}_{\mathcal{O}_{K_1}}(R_{1,\mathcal{M}}, \mathfrak{B}_i) \to \operatorname{Hom}'_{(\operatorname{Mod}/\tilde{S}_1)}(\tilde{\mathcal{M}}, \mathfrak{B}_i/\operatorname{Fil}^p\mathfrak{B}_i)$$

est un isomorphisme pour $i \geq 0$. Pour cela, il suffit de montrer que toute solution dans $(\mathfrak{B}_i/\operatorname{Fil}^p\mathfrak{B}_i)^d$ du système:

(S) $$\begin{pmatrix} \phi_1(\pi_1^{r_1}\mathfrak{b}_1) \\ \vdots \\ \phi_1(\pi_1^{r_d}\mathfrak{b}_d) \end{pmatrix} = \mathcal{G}_\pi \begin{pmatrix} \mathfrak{b}_1 \\ \vdots \\ \mathfrak{b}_d \end{pmatrix}$$

où $\pi_1^{r_j}\mathfrak{b}_j \in \operatorname{Fil}^1\mathfrak{B}_i/\operatorname{Fil}^p\mathfrak{B}_i$ se remonte de façon unique en une solution (du même système (S)) dans $(\mathfrak{B}_i)^d$. Comme $\mathfrak{B}_i$ est noethérien, $\operatorname{Fil}^1\mathfrak{B}_i$ est un idéal de type fini et la topologie $\operatorname{Fil}^1\mathfrak{B}_i$-adique sur $\mathfrak{B}_i$ s'identifie à la topologie $p$-adique, de sorte que $\mathfrak{B}_i \xrightarrow{\sim} \varprojlim_m(\mathfrak{B}_i/\operatorname{Fil}^m\mathfrak{B}_i)$. De proche en proche, il suffit donc de montrer que toute solution de (S) dans $(\mathfrak{B}_i/\operatorname{Fil}^m\mathfrak{B}_i)^d$ se remonte de façon unique dans $(\mathfrak{B}_i/\operatorname{Fil}^{m+1}\mathfrak{B}_i)^d$ pour $m \geq p$. Soient donc $(\mathfrak{b}_1, \ldots, \mathfrak{b}_d)$ une solution de (S) dans $(\mathfrak{B}_i/\operatorname{Fil}^m\mathfrak{B}_i)^d$ et $\hat{\mathfrak{b}}_j$ des relevés dans $\mathfrak{B}_i/\operatorname{Fil}^{m+1}\mathfrak{B}_i$, on a $\pi_1^{r_j}\hat{\mathfrak{b}}_j \in \operatorname{Fil}^1\mathfrak{B}_i/\operatorname{Fil}^{m+1}\mathfrak{B}_i$ et on cherche des $\delta_j \in \operatorname{Fil}^m\mathfrak{B}_i/\operatorname{Fil}^{m+1}\mathfrak{B}_i$ tels que $(\hat{\mathfrak{b}}_j + \delta_j)$ soit une solution de (S) dans $(\mathfrak{B}_i/\operatorname{Fil}^{m+1}\mathfrak{B}_i)^d$. Soit $\mathfrak{b}$ un élément de $\operatorname{Fil}^m\mathfrak{B}_i$, $\mathfrak{b} = \sum \mathfrak{x}_1 \cdots \mathfrak{x}_m$ avec $\mathfrak{x}_j \in \operatorname{Fil}^1\mathfrak{B}_i$. Un calcul donne $\phi_1(\mathfrak{b}) = \sum \frac{\mathfrak{x}_1^p \cdots \mathfrak{x}_m^p}{-p}$ + un élément de $\operatorname{Fil}^{m+1}\mathfrak{B}_i$ et $\frac{\mathfrak{x}_1^p \cdots \mathfrak{x}_m^p}{-p} \in \operatorname{Fil}^{p(m-1)}\mathfrak{B}_i$. Puisque $p \geq 3$, on a $p(m-1) \geq m+1$ quand $m \geq p$ et $\phi_1(\mathfrak{b}) \in \operatorname{Fil}^{m+1}\mathfrak{B}_i$. De cela on déduit $\phi_1(\pi_1^{r_j}(\hat{\mathfrak{b}}_j + \delta_j)) = \phi_1(\pi_1^{r_j}\hat{\mathfrak{b}}_j)$ dans $\mathfrak{B}_i/\operatorname{Fil}^{m+1}\mathfrak{B}_i$ et les $(\hat{\mathfrak{b}}_j + \delta_j)$ satisfont (S) si et seulement si:

$$\begin{pmatrix} \delta_1 \\ \vdots \\ \delta_d \end{pmatrix} = \mathcal{G}_\pi^{-1}\left(\begin{pmatrix} \phi_1(\pi_1^{r_1}\hat{\mathfrak{b}}_1) \\ \vdots \\ \phi_1(\pi_1^{r_d}\hat{\mathfrak{b}}_d) \end{pmatrix} - \mathcal{G}_\pi\begin{pmatrix} \hat{\mathfrak{b}}_1 \\ \vdots \\ \hat{\mathfrak{b}}_d \end{pmatrix}\right)$$

dans $(\operatorname{Fil}^m\mathfrak{B}_i/\operatorname{Fil}^{m+1}\mathfrak{B}_i)^d$, d'où le résultat. $\square$

COROLLAIRE 3.1.8. *Le morphisme de faisceaux sur* $\operatorname{Spf}(\mathcal{O}_K)_{\operatorname{syn}}$:

$$\operatorname{Spf}(R_\mathcal{M})^\sim \to \operatorname{Gr}(\mathcal{M})$$

*défini en* (3.1.5) *est un isomorphisme.*

*Preuve.* Il suffit de le vérifier localement pour la topologie syntomique. On reprend les notations précédant (3.1.6). Pour tout $i$, on a défini en (3.1.5) un morphisme $\mathrm{Hom}_{\mathcal{O}_K}(R_{\mathcal{M}}, \mathfrak{A}_i) \to \mathrm{Hom}_{'(\mathrm{Mod}/S_1)}(\mathcal{M}, \mathcal{O}^{\mathrm{cris}}_{1,\pi}(\mathfrak{A}_i))$. Par passage à la limite inductive, on obtient un diagramme dont on vérifie facilement la commutativité à partir de la construction des diverses flèches:

$$\begin{array}{ccc} \mathrm{Hom}_{\mathcal{O}_K}(R_{\mathcal{M}}, \mathfrak{A}_\infty) & \longrightarrow & \mathrm{Hom}_{'(\mathrm{Mod}/S_1)}(\mathcal{M}, \mathcal{O}^{\mathrm{cris}}_{1,\pi}(\mathfrak{A}_\infty)) \\ \wr \downarrow (déf.) & & (3.1.6) \downarrow \wr \\ \mathrm{Hom}_{\mathcal{O}_{K_1}}(R_{1,\mathcal{M}}, \mathfrak{B}_\infty) & \stackrel{(3.1.7)}{\xrightarrow{\sim}} & \mathrm{Hom}_{'(\mathrm{Mod}/\tilde{S}_1)}(\tilde{\mathcal{M}}, \mathfrak{B}_\infty/\mathrm{Fil}^p\mathfrak{B}_\infty) \end{array}$$

et la flèche horizontale supérieure est bien un isomorphisme puisque les trois autres le sont. □

3.2. *Le module associé à un schéma en groupes annulés par $p$.* A tout groupe $G$ de $(p\text{-Gr}/\mathcal{O}_K)$ annulé par $p$, on associe un objet $\mathrm{Mod}(G)$ de $'(\mathrm{Mod}/S_1)$ défini par:

- $\mathrm{Mod}(G) = \mathrm{Hom}_{(\mathrm{Ab}/\mathcal{O}_K)}(G, \mathcal{O}^{\mathrm{cris}}_{1,\pi})$, la structure de $S_1$-module provenant de celle de $\mathcal{O}^{\mathrm{cris}}_{1,\pi}$;

- $\mathrm{Fil}^1\mathrm{Mod}(G) = \mathrm{Hom}_{(\mathrm{Ab}/\mathcal{O}_K)}(G, \mathcal{J}^{\mathrm{cris}}_{1,\pi})$;

- $\phi_1 : \mathrm{Fil}^1\mathrm{Mod}(G) \to \mathrm{Mod}(G)$ est induit par $\phi_1 : \mathcal{J}^{\mathrm{cris}}_{1,\pi} \to \mathcal{O}^{\mathrm{cris}}_{1,\pi}$ (2.3).

Nous allons montrer, en utilisant de façon essentielle plusieurs résultats de [BBM], que $\mathrm{Mod}(G)$ est un objet de $(\mathrm{Mod}/S_1)$, i.e. est un $S_1$-module libre de type fini tel que $\phi_1(\mathrm{Fil}^1\mathrm{Mod}(G))$ engendre $\mathrm{Mod}(G)$ sur $S_1$.

Pour suivre les notations de [BBM], nous posons, uniquement dans cette partie, $S = \mathrm{Spec}(\mathcal{O}_K/p)$, qu'on prendra garde à ne pas confondre avec l'anneau $S$ de (2.1.1) (nous utiliserons seulement les anneaux $S_n$ ici). Soit $n \in \mathbf{N}^*$, on rappelle qu'on a fixé un épaississement à puissances divisées $S \hookrightarrow E_n$ où $E_n = \mathrm{Spec}(S_n)$ en envoyant $\gamma_i(u^e)$ sur $(F(\pi))^i p^i/i! = 0$ dans $\mathcal{O}_K/p$ (cf. 2.3, les puissances divisées sur $\mathrm{Ker}(S_n \to \mathcal{O}_K/p)$ étant bien sûr prises compatibles aux puissances divisées canoniques sur $pS_n$). On note $(S/E_n)_{CRIS}$ la catégorie dont les objets sont les quadruples $(U, T, i, \delta)$ où $U$ est un $S$-schéma, $T$ un $E_n$-schéma, $i : U \hookrightarrow T$ une immersion fermée avec un diagramme commutatif:

$$\begin{array}{ccc} U & \stackrel{i}{\hookrightarrow} & T \\ \downarrow & & \downarrow \\ S & \hookrightarrow & E_n \end{array}$$

et $\delta$ une structure d'idéal à puissances divisées sur l'idéal définissant $i$, compatible aux puissances divisées sur $E_n$. On abrègera $(U, T, i, \delta)$ en $(U, T)$. On définit les morphismes de la façon usuelle ([BBM, 1.1.1]) et on dit qu'un morphisme:

$$\begin{array}{ccc} U' & \hookrightarrow & T' \\ \downarrow & & \downarrow \\ U & \hookrightarrow & T \end{array}$$

est syntomique s'il est cartésien (i.e. $U' = U \times_T T'$) et si $T' \to T$ est syntomique. On munit $(S/E_n)_{\mathrm{CRIS}}$ de la topologie engendrée par les familles surjectives de morphismes syntomiques. Pour tout objet $(U,T)$ de $(S/E_n)_{\mathrm{CRIS}}$, on définit de manière analogue le site $(U/T)_{\mathrm{CRIS}}$ en remplaçant $S$ par $U$ et $E_n$ par $T$ et, pour tout $S$-schéma $X$, le site $(X/E_n)_{\mathrm{CRIS}}$ en remplaçant $S$ par $X$. On note $\mathrm{Ab}_{S/E_n}$ et $\mathrm{Ab}_{U/T}$ les catégories de faisceaux abéliens sur $(S/E_n)_{\mathrm{CRIS}}$ et $(U/T)_{\mathrm{CRIS}}$ (définies comme d'habitude, cf. [BBM, 1.1.3]). Soit $G$ un schéma en groupes fini et plat sur $\mathcal{O}_K$, $G_1 = G \times_{\mathrm{Spec}(\mathcal{O}_K)} \mathrm{Spec}(\mathcal{O}_K/p)$, on note avec ([BBM, 1.3.1]) $\underline{G}_1$ le faisceau sur $(S/E_n)_{\mathrm{CRIS}}$ défini par $\underline{G}_1(U,T) = G_1(U) = \mathrm{Hom}_S(U,G_1)$. Pour tout faisceau $F$ de $\mathrm{Ab}_{S/E_n}$, soit $\mathcal{H}\mathrm{om}_{S/E_n}(\underline{G}_1,F)$ le faisceau défini par $\mathcal{H}\mathrm{om}_{S/E_n}(\underline{G}_1,F)(U,T) = \mathrm{Hom}_{\mathrm{Ab}_{U/T}}(\underline{G}_{1,U},F_U)$, où $\underline{G}_{1,U}$ et $F_U$ désignent les restrictions de $\underline{G}_1$ et $F$ à $(U/T)_{\mathrm{CRIS}}$, et:

$$\mathrm{Hom}_{S/E_n}(\underline{G}_1,F) = \mathrm{Hom}_{\mathrm{Ab}_{S/E_n}}(\underline{G}_1,F) = H^0((S/E_n)_{\mathrm{CRIS}}, \mathcal{H}\mathrm{om}_{S/E_n}(\underline{G}_1,F)).$$

Pour $i \geq 1$, on définit également les faisceaux $\mathcal{E}\mathrm{xt}^i_{S/E_n}(\underline{G}_1,F)$ comme les faisceaux associés aux préfaisceaux $(U,T) \mapsto \mathrm{Ext}^i_{\mathrm{Ab}_{U/T}}(\underline{G}_{1,U},F_U)$ (cf. [BBM],1.3.3) et $\mathrm{Ext}^i_{S/E_n}(\underline{G}_1,F) = H^0((S/E_n)_{\mathrm{CRIS}}, \mathcal{E}\mathrm{xt}^i_{S/E_n}(\underline{G}_1,F))$. Soient $S_{\mathrm{SYN}}$ la catégorie des $S$-schémas munie de la topologie syntomique et $(\mathrm{Ab}/S)$ la catégorie des faisceaux abéliens sur $S_{\mathrm{SYN}}$, on définit enfin un foncteur $w_* : \mathrm{Ab}_{S/E_n} \to (\mathrm{Ab}/S)$ qui à un faisceau $F$ associe:

$$\begin{array}{rcl} X \mapsto w_*(F)(X) & = & H^0((X/E_n)_{\mathrm{CRIS}}, F|_{(X/E_n)_{\mathrm{CRIS}}}) \\ & = & H^0((X/E_n)_{\mathrm{cris}}, F|_{(X/E_n)_{\mathrm{cris}}}) \end{array}$$

qui est bien un faisceau sur $S_{\mathrm{SYN}}$ ($(X/E_n)_{\mathrm{cris}}$ est le petit site cristallin, cf. 2.3).

LEMME 3.2.1. *Avec les notations précédentes, soit $F$ un faisceau de $\mathrm{Ab}_{S/E_n}$, on a:*

$$\mathrm{Hom}_{S/E_n}(\underline{G}_1,F) = \mathrm{Hom}_{(\mathrm{Ab}/S)}(G_1, w_*F).$$

*Preuve.* Si $H \in (\mathrm{Ab}/S)$, on note $w^*H$ le faisceau sur $(S/E_n)_{\mathrm{CRIS}}$ défini par $(w^*H)(U,T) = H(U)$: $(w^*, w_*)$ forme un couple de foncteurs adjoints ([Be, III.4.4]). Comme $\underline{G}_1 = w^*G_1$, on a:

$$\mathrm{Hom}_{S/E_n}(\underline{G}_1, F) = \mathrm{Hom}_{S/E_n}(w^*G_1, F) = \mathrm{Hom}_{(\mathrm{Ab}/S)}(G_1, w_*F). \qquad \square$$

On note $\mathbf{G}_a$ le $S$-schéma en groupes $\mathrm{Spec}((\mathcal{O}_K/p\mathcal{O}_K)[X])$ et $\underline{\mathbf{G}}_a$ le faisceau de $\mathrm{Ab}_{S/E_n}$ défini par $\underline{\mathbf{G}}_a(U,T) = \mathrm{Hom}_S(U, \mathbf{G}_a)$.

COROLLAIRE 3.2.2. *On a des isomorphismes*:

$$\begin{aligned}
\mathrm{Hom}_{S/E_n}(\underline{G}_1, \mathcal{O}_{S/E_n}) &\simeq \mathrm{Hom}_{(\mathrm{Ab}/\mathcal{O}_K)}(G, \mathcal{O}_{n,\pi}^{\mathrm{cris}}), \\
\mathrm{Hom}_{S/E_n}(\underline{G}_1, \underline{\mathbf{G}}_a) &\simeq \mathrm{Hom}_{(\mathrm{Ab}/\mathcal{O}_K)}(G, \mathcal{O}_1).
\end{aligned}$$

*Preuve.* Soient $\Delta_{G_1}$ (resp. $\Delta_G$) la multiplication $G_1 \times_S G_1 \to G_1$ (resp. $G \times_{\mathrm{Spf}(\mathcal{O}_K)} G \to G$) et $\mathrm{pr}_i$, $i \in \{1,2\}$, les deux projections $G_1 \times_S G_1 \to G_1$ (resp. $G \times_{\mathrm{Spf}(\mathcal{O}_K)} G \to G$), on a (avec des notations évidentes):

$$\begin{aligned}
\mathrm{Hom}_{(\mathrm{Ab}/S)}(G_1, w_*\mathcal{O}_{S/E_n}) &= \{x \in w_*\mathcal{O}_{S/E_n}(G_1) / (\Delta_{G_1} - \mathrm{pr}_1 - \mathrm{pr}_2)(x) = 0\} \\
&= \{x \in \mathcal{O}_{n,\pi}^{\mathrm{cris}}(G) / (\Delta_G - \mathrm{pr}_1 - \mathrm{pr}_2)(x) = 0\} \\
&= \mathrm{Hom}_{(\mathrm{Ab}/\mathcal{O}_K)}(G, \mathcal{O}_{n,\pi}^{\mathrm{cris}}) \text{ par } (2.2.2),
\end{aligned}$$

d'où le résultat dans ce cas par (3.2.1). La preuve avec $\underline{\mathbf{G}}_a$ à la place de $\mathcal{O}_{S/E_n}$ est la même puisque $w_*\underline{\mathbf{G}}_a$ est le faisceau structural sur $S_{\mathrm{SYN}}$ ([Be, III.4.4.8]). □

Dans [BBM, 3.1.5] il est associé à $G_1$ un cristal $\mathbf{D}(G_1)$ (dit de Dieudonné) sur $(S/E_n)_{\mathrm{CRIS}}$ défini par $\mathbf{D}(G_1) = \mathcal{E}\mathrm{xt}^1_{S/E_n}(\underline{G}_1, \mathcal{O}_{S/E_n})$ ([BBM, 3.1.3]). En fait, le cristal construit dans (loc. cit.) vit sur un site plus gros, noté $(S/\Sigma)_{\mathrm{CRIS}}$, et ce que nous appelons $\mathbf{D}(G_1)$ ici en est la restriction au site $(S/E_n)_{\mathrm{CRIS}}$ ([BBM, 1.3.3.8]). La topologie syntomique n'est pas considérée dans [BBM], mais il s'agit bien du même cristal par [BBM, 2.3.11]. En particulier:

PROPOSITION 3.2.3 ([BBM, 4.3.1]). *Supposons $G$ tué par $p$ de rang $p^d$, alors le cristal $\mathbf{D}(G_1)$ est un faisceau de $\mathcal{O}_{S/E_n}/p\mathcal{O}_{S/E_n}$-modules localement libres de rang $d$.*

PROPOSITION 3.2.4. *Soit $G$ un schéma en groupes fini et plat sur $\mathcal{O}_K$ tué par $p^n$, on a un isomorphisme canonique de $S_n$-modules*:

$$H^0((S/E_n)_{\mathrm{CRIS}}, \mathbf{D}(G_1)) = \mathbf{D}(G_1)(S, E_n) \xrightarrow{\sim} Hom_{S/E_n}(\underline{G}_1, \mathcal{O}_{S/E_n}).$$

*Preuve.* On applique [BBM, 4.2.9] avec $m = n$ et $T_m = E_n$. □

Comme $E_n = \mathrm{Spec}(S_n)$ où $S_n$ est un anneau local, on déduit de (3.2.2), (3.2.4) et (3.2.3):

COROLLAIRE 3.2.5. *Soit $G$ un schéma en groupes fini et plat sur $\mathcal{O}_K$ tué par $p$ de rang $p^d$, le $S_1$-module $\mathrm{Mod}(G) = \mathrm{Hom}_{(\mathrm{Ab}/\mathcal{O}_K)}(G, \mathcal{O}_{1,\pi}^{\mathrm{cris}})$ est libre de rang $d$.*

Avant de passer au $\phi_1$, notons le lemme suivant, qui sera utile:

LEMME 3.2.6. *Soit $0 \to G' \to G \to G'' \to 0$ une suite exacte de schémas en groupes finis et plats sur $\mathcal{O}_K$ tués par $p^n$, alors on a une suite exacte de $S_n$-modules*:

$$0 \to \mathrm{Hom}_{(\mathrm{Ab}/\mathcal{O}_K)}(G'', \mathcal{O}_{n,\pi}^{\mathrm{cris}}) \to \mathrm{Hom}_{(\mathrm{Ab}/\mathcal{O}_K)}(G, \mathcal{O}_{n,\pi}^{\mathrm{cris}})$$
$$\to \mathrm{Hom}_{(\mathrm{Ab}/\mathcal{O}_K)}(G', \mathcal{O}_{n,\pi}^{\mathrm{cris}}) \to 0.$$

*Preuve.* Par [BBM, 4.2.7], on a une suite exacte courte de $\mathcal{O}_{E_n}$-modules dans $(E_n)_{\mathrm{zar}}^\sim$:

$$0 \to \mathcal{E}\mathrm{xt}^1_{S/E_n}(\underline{G}_1'', \mathcal{O}_{S/E_n})|_{E_{n,\mathrm{zar}}} \to \mathcal{E}\mathrm{xt}^1_{S/E_n}(\underline{G}_1, \mathcal{O}_{S/E_n})|_{E_{n,\mathrm{zar}}}$$
$$\to \mathcal{E}\mathrm{xt}^1_{S/E_n}(\underline{G}_1', \mathcal{O}_{S/E_n})|_{E_{n,\mathrm{zar}}} \to 0.$$

Comme ces faisceaux sont quasi-cohérents ([BBM, 2.3.1]) et que $E_n$ est affine, on a encore une suite exacte avec les $H^0((S/E_n)_{\mathrm{CRIS}}, -)$. Par (3.2.4) et (3.2.2), on en déduit une suite exacte comme dans l'énoncé. □

Nous montrons maintenant que $\phi_1(\mathrm{Fil}^1\mathrm{Mod}(G))$ engendre $\mathrm{Mod}(G)$ ($G$ tué par $p$). Bien que certains résultats soient valables sur $(S/E_n)_{\mathrm{CRIS}}$, $n \geq 1$, et avec un schéma en groupes fini et plat $G$ annulé par une puissance quelconque de $p$, on se restreint la plupart du temps au cas $n = 1$ et $pG = 0$ pour simplifier. Pour tout $(U,T)$ de $(S/E_1)_{\mathrm{CRIS}}$, le morphisme $\mathcal{O}_T \to \mathcal{O}_U$ induit une surjection $\mathcal{O}_{S/E_1} \to \underline{\mathbf{G}}_a$ dans $\mathrm{Ab}_{S/E_1}$. Comme $\mathcal{E}\mathrm{xt}^1_{S/E_1}(\underline{G}_1, \mathcal{O}_{S/E_1})$, le faisceau $\mathcal{E}\mathrm{xt}^1_{S/E_1}(\underline{G}_1, \underline{\mathbf{G}}_a)$ est un $\mathcal{O}_{S/E_1}$-module (quasi-cohérent au sens de [BBM, 1.1.18]) et son calcul pour la topologie syntomique s'identifie à son calcul pour la topologie de Zariski ([BBM, 2.3.12]).

LEMME 3.2.7. *Soit $G$ un schéma en groupes fini et plat sur $\mathcal{O}_K$ tué par $p$ et $G_1 = G \times_{\mathrm{Spec}(\mathcal{O}_K)} \mathrm{Spec}(\mathcal{O}_K/p)$, la surjection $\mathcal{O}_{S/E_1} \to \underline{\mathbf{G}}_a$ induit une surjection de $\mathcal{O}_{S/E_1}$-modules $\mathcal{E}\mathrm{xt}^1_{S/E_1}(\underline{G}_1, \mathcal{O}_{S/E_1}) \to \mathcal{E}\mathrm{xt}^1_{S/E_1}(\underline{G}_1, \underline{\mathbf{G}}_a)$ dans $\mathrm{Ab}_{S/E_1}$.*

*Preuve.* Pour tout groupe $p$-divisible $H$ sur $S$, on note $\underline{H}$ le faisceau $(U,T) \mapsto H(U)$ sur $(S/E_n)_{\mathrm{CRIS}}$ ($n \geq 1$). Rappelons que $\mathcal{J}_{S/E_1}$ est le faisceau sur $(S/E_1)_{\mathrm{CRIS}}$ défini par $\mathcal{J}_{S/E_1}(U,T) = \mathrm{Ker}(\Gamma(T, \mathcal{O}_T) \to \Gamma(U, \mathcal{O}_U))$ et qu'on a une suite exacte dans $\mathrm{Ab}_{S/E_1}$: $0 \to \mathcal{J}_{S/E_1} \to \mathcal{O}_{S/E_1} \to \underline{\mathbf{G}}_a \to 0$. Puisque $S$ est le spectre d'un anneau local, par [BBM, 3.1.1] il existe un $S$-schéma abélien $A$ et une $S$-immersion fermée $G_1 \hookrightarrow A$ (c'est un cas particulier d'un théorème de Raynaud). Soit $H$ le groupe $p$-divisible associé à $A$, alors $H' = H/G_1$ est encore un groupe $p$-divisible ([BBM, 3.3.12]), de sorte qu'on a une suite exacte pour la topologie fppf: $0 \to G_1 \to H \to H' \to 0$. Pour $i \leq 2$ et $* \in \{G_1, H, H'\}$, les faisceaux $\mathcal{E}\mathrm{xt}^i_{S/E_1}(\underline{*}, \mathcal{J}_{S/E_1})$, $\mathcal{E}\mathrm{xt}^i_{S/E_1}(\underline{*}, \mathcal{O}_{S/E_1})$ et $\mathcal{E}\mathrm{xt}^i_{S/E_1}(\underline{*}, \mathbf{G}_a)$ sont égaux aux $\mathcal{E}\mathrm{xt}^i$ correspondants calculés pour la topologie fppf ([BBM, 2.3.10, 2.4.5, 2.3.12 et 2.4.6]). On a donc un diagramme commutatif dans $\mathrm{Ab}_{S/E_1}$:

$$\begin{array}{ccccc}
\mathcal{E}\mathrm{xt}^1_{S/E_1}(\underline{H}, \mathcal{O}_{S/E_1}) & \to & \mathcal{E}\mathrm{xt}^1_{S/E_1}(\underline{G}_1, \mathcal{O}_{S/E_1}) & & \\
\downarrow & & \downarrow & & \\
\mathcal{E}\mathrm{xt}^1_{S/E_1}(\underline{H}, \mathbf{G}_a) & \to & \mathcal{E}\mathrm{xt}^1_{S/E_1}(\underline{G}_1, \mathbf{G}_a) & \to & \mathcal{E}\mathrm{xt}^2_{S/E_1}(\underline{H}', \mathbf{G}_a) \\
\downarrow & & & & \\
\mathcal{E}\mathrm{xt}^2_{S/E_1}(\underline{H}, \mathcal{J}_{S/E_1}) & & & &
\end{array}$$

où les deux suites de trois termes sont exactes pour la topologie fppf. Par [BBM, 3.3.2(iii), 3.3.4(iii)) et 1.3.8], on a

$$\mathcal{E}\mathrm{xt}^2_{S/E_1}(\underline{H}', \mathbf{G}_a) = \mathcal{E}\mathrm{xt}^2_{S/E_1}(\underline{H}, \mathcal{J}_{S/E_1}) = 0.$$

Une chasse au diagramme donne alors la surjection pour la topologie fppf, donc pour la topologie syntomique puisque les deux faisceaux sont quasi-cohérents au sens de [BBM, 1.1.18] (voir [BBM, 2.3.1, 2.3.2, 2.3.9 et 2.3.12]). □

COROLLAIRE 3.2.8. *Avec les notations de* (3.2.7), *on a une surjection de $S_1$-modules*:
$$\mathrm{Ext}^1_{S/E_1}(\underline{G}_1, \mathcal{O}_{S/E_1}) \to \mathrm{Ext}^1_{S/E_1}(\underline{G}_1, \mathbf{G}_a).$$

*Preuve*. Pour tout faisceau $F$ sur $(S/E_1)_{\mathrm{CRIS}}$, on a $H^0((S/E_1)_{\mathrm{CRIS}}, F) = F(S, E_1)$ puisque $(S, E_1)$ est l'objet final du site. Les faisceaux

$$\mathcal{E}\mathrm{xt}^1_{S/E_1}(\underline{G}_1, \mathcal{O}_{S/E_1}) \quad \text{et} \quad \mathcal{E}\mathrm{xt}^1_{S/E_1}(\underline{G}_1, \mathbf{G}_a)$$

sont quasi-cohérents lorsque restreints à $(E_1)_{\mathrm{zar}}$ (voir preuve précédente) et comme $E_1$ est affine, le résultat est clair en vertu de (3.2.7). □

Pour tout groupe $p$-divisible $H$ sur $S$, on note $H(n)$ le schéma en groupes fini et plat $\mathrm{Ker}(p^n : H \to H)$ et $\underline{H}(n)$ le faisceau $(U, T) \mapsto \mathrm{Hom}_S(U, H(n))$ sur $(S/E_n)_{\mathrm{CRIS}}$ ($n \geq 1$).

PROPOSITION 3.2.9. *Soit $H$ un groupe $p$-divisible sur $S$, on a une surjection de $S_1$-modules*: $\mathrm{Hom}_{S/E_1}(\underline{H}(1), \mathcal{O}_{S/E_1}) \to \mathrm{Hom}_{S/E_1}(\underline{H}(1), \mathbf{G}_a)$.

*Preuve*. Si $F$ est un faisceau sur $(S/E_1)_{\mathrm{CRIS}}$ annulé par $p$, on a une application naturelle $\mathrm{Ext}^1_{\mathrm{Ab}_{S/E_1}}(\underline{H}(1), F) \to \mathrm{Hom}_{S/E_1}(\underline{H}(1), F)$ qui à une classe d'extension $0 \to F \to \mathcal{E} \to \underline{H}(1) \to 0$ dans $\mathrm{Ab}_{S/E_1}$ associe l'homomorphisme $\underline{H}(1) \to F$ donné par le diagramme du serpent relatif à la multiplication par $p$ sur cette extension. Cette application étant clairement fonctorielle en $F$, on a un diagramme commutatif:

$$\begin{array}{ccc} \mathrm{Ext}^1_{\mathrm{Ab}_{S/E_1}}(\underline{H}(1), \mathcal{O}_{S/E_1}) & \to & \mathrm{Hom}_{S/E_1}(\underline{H}(1), \mathcal{O}_{S/E_1}) \\ \downarrow & & \downarrow \\ \mathrm{Ext}^1_{\mathrm{Ab}_{S/E_1}}(\underline{H}(1), \underline{\mathbf{G}}_a) & \to & \mathrm{Hom}_{S/E_1}(\underline{H}(1), \underline{\mathbf{G}}_a) \end{array}$$

qui se factorise par:

$$\begin{array}{ccc} \mathrm{Ext}^1_{S/E_1}(\underline{H}(1), \mathcal{O}_{S/E_1}) & \to & \mathrm{Hom}_{S/E_1}(\underline{H}(1), \mathcal{O}_{S/E_1}) \\ \downarrow & & \downarrow \\ \mathrm{Ext}^1_{S/E_1}(\underline{H}(1), \underline{\mathbf{G}}_a) & \to & \mathrm{Hom}_{S/E_1}(\underline{H}(1), \underline{\mathbf{G}}_a) \end{array}$$

puisque le préfaisceau des homomorphismes locaux est un faisceau. Par (3.2.8) et une chasse au diagramme triviale, il suffit de montrer la surjectivité de la flèche horizontale inférieure, ou encore la surjectivité de la flèche:

$$(*) \qquad \mathrm{Ext}^1_{\mathrm{Ab}_{S/E_1}}(\underline{H}(1), \underline{\mathbf{G}}_a) \to \mathrm{Hom}_{S/E_1}(\underline{H}(1), \underline{\mathbf{G}}_a).$$

Mais, puisque $H$ est un groupe $p$-divisible, on a une extension dans $\mathrm{Ab}_{S/E_1}$: $0 \to \underline{H}(1) \to \underline{H}(2) \xrightarrow{p} \underline{H}(1) \to 0$ (en effet, $0 \to H(1) \to H(2) \xrightarrow{p} H(1) \to 0$ est exact pour la topologie syntomique sur $S$ en utilisant la première partie de la preuve de (2.2.3) et on applique la propriété 4 des morphismes syntomiques pour voir que la suite reste exacte dans $(S/E_1)_{\mathrm{CRIS}}$), elle fournit donc un homomorphisme cobord $\mathrm{Hom}_{S/E_1}(\underline{H}(1), \underline{\mathbf{G}}_a) \to \mathrm{Ext}^1_{\mathrm{Ab}_{S/E_1}}(\underline{H}(1), \underline{\mathbf{G}}_a)$ qui envoie un homomorphisme $f$ sur la classe d'extensions $0 \to \underline{\mathbf{G}}_a \to \underline{H}(2) \oplus_{\underline{H}(1),f} \underline{\mathbf{G}}_a \to \underline{H}(1) \to 0$. On vérifie en regardant le diagramme du serpent de la multiplication par $p$ sur cette extension qu'on a construit une section de la flèche $(*)$, d'où le résultat. □

De (3.2.2), on déduit:

COROLLAIRE 3.2.10. *Soit $H$ un groupe $p$-divisible sur $\mathrm{Spec}(\mathcal{O}_K)$, on a une surjection de $S_1$-modules*:

$$\mathrm{Hom}_{(\mathrm{Ab}/\mathcal{O}_K)}(H(1), \mathcal{O}^{\mathrm{cris}}_{1,\pi}) \to \mathrm{Hom}_{(\mathrm{Ab}/\mathcal{O}_K)}(H(1), \mathcal{O}_1).$$

LEMME 3.2.11. *Soit $H$ un groupe $p$-divisible sur $\mathrm{Spec}(\mathcal{O}_K)$, on a des surjections de $S_2$-modules*:

$$\begin{array}{ccc} \mathrm{Hom}_{(\mathrm{Ab}/\mathcal{O}_K)}(H(2), \mathcal{O}^{\mathrm{cris}}_{2,\pi}) & \to & \mathrm{Hom}_{(\mathrm{Ab}/\mathcal{O}_K)}(H(1), \mathcal{O}^{\mathrm{cris}}_{1,\pi}), \\ \mathrm{Hom}_{(\mathrm{Ab}/\mathcal{O}_K)}(H(2), \mathcal{J}^{\mathrm{cris}}_{2,\pi}) & \to & \mathrm{Hom}_{(\mathrm{Ab}/\mathcal{O}_K)}(H(1), \mathcal{J}^{\mathrm{cris}}_{1,\pi}). \end{array}$$

*Preuve.* Par (2.3.3), on a un isomorphisme:

$$\mathrm{Hom}_{(\mathrm{Ab}/\mathcal{O}_K)}(H(1), \mathcal{O}^{\mathrm{cris}}_{1,\pi}) \simeq \mathrm{Hom}_{(\mathrm{Ab}/\mathcal{O}_K)}(H(1), \mathcal{O}^{\mathrm{cris}}_{2,\pi})$$

et, par (2.2.3), une suite exacte dans $(\mathrm{Ab}/\mathcal{O}_K)$: $0 \to H(1) \to H(2) \xrightarrow{p} H(1) \to 0$, d'où la première surjectivité par (3.2.6) appliqué avec $n = 2$. On a aussi un isomorphisme:

$$\mathrm{Hom}_{(\mathrm{Ab}/\mathcal{O}_K)}(H(1), \mathcal{O}^{\mathrm{cris}}_{1,\pi}) \simeq \mathrm{Hom}_{(\mathrm{Ab}/\mathcal{O}_K)}(H(2), \mathcal{O}^{\mathrm{cris}}_{1,\pi})$$

(resp. avec $\mathcal{J}^{\mathrm{cris}}_{1,\pi}$ et $\mathcal{O}_1$). Le diagramme commutatif de suites exactes dans $(\mathrm{Ab}/\mathcal{O}_K)$ (2.3.3):

$$\begin{array}{ccccccccc}
& & 0 & & 0 & & 0 & & \\
& & \downarrow & & \downarrow & & \downarrow & & \\
0 & \to & \mathcal{J}^{\mathrm{cris}}_{1,\pi} & \to & \mathcal{O}^{\mathrm{cris}}_{1,\pi} & \to & \mathcal{O}_1 & \to & 0 \\
& & \downarrow & & \downarrow & & \downarrow & & \\
0 & \to & \mathcal{J}^{\mathrm{cris}}_{2,\pi} & \to & \mathcal{O}^{\mathrm{cris}}_{2,\pi} & \to & \mathcal{O}_2 & \to & 0 \\
& & \downarrow & & \downarrow & & \downarrow & & \\
0 & \to & \mathcal{J}^{\mathrm{cris}}_{1,\pi} & \to & \mathcal{O}^{\mathrm{cris}}_{1,\pi} & \to & \mathcal{O}_1 & \to & 0 \\
& & \downarrow & & \downarrow & & \downarrow & & \\
& & 0 & & 0 & & 0 & &
\end{array}$$

fournit un diagramme commutatif de suites exactes:

$$\begin{array}{ccccccccc}
& & 0 & & 0 & & 0 & & \\
& & \downarrow & & \downarrow & & \downarrow & & \\
0 \to & \mathrm{Hom}_{(\mathrm{Ab}/\mathcal{O}_K)}(H(2),\mathcal{J}^{\mathrm{cris}}_{1,\pi}) & \to & \mathrm{Hom}_{(\mathrm{Ab}/\mathcal{O}_K)}(H(2),\mathcal{O}^{\mathrm{cris}}_{1,\pi}) & \to & \mathrm{Hom}_{(\mathrm{Ab}/\mathcal{O}_K)}(H(2),\mathcal{O}_1) & \to & 0 \\
& & \downarrow & & \downarrow & & \downarrow & & \\
0 \to & \mathrm{Hom}_{(\mathrm{Ab}/\mathcal{O}_K)}(H(2),\mathcal{J}^{\mathrm{cris}}_{2,\pi}) & \to & \mathrm{Hom}_{(\mathrm{Ab}/\mathcal{O}_K)}(H(2),\mathcal{O}^{\mathrm{cris}}_{2,\pi}) & \to & \mathrm{Hom}_{(\mathrm{Ab}/\mathcal{O}_K)}(H(2),\mathcal{O}_2) & & \\
& & \downarrow & & \downarrow & & \downarrow & & \\
0 \to & \mathrm{Hom}_{(\mathrm{Ab}/\mathcal{O}_K)}(H(2),\mathcal{J}^{\mathrm{cris}}_{1,\pi}) & \to & \mathrm{Hom}_{(\mathrm{Ab}/\mathcal{O}_K)}(H(2),\mathcal{O}^{\mathrm{cris}}_{1,\pi}) & \to & \mathrm{Hom}_{(\mathrm{Ab}/\mathcal{O}_K)}(H(2),\mathcal{O}_1) & & \\
& & & & \downarrow & & & & \\
& & & & 0 & & & &
\end{array}$$

où les surjections proviennent de (3.2.10) et de ce qui précède. Une chasse au diagramme donne la deuxième surjection. $\square$

Pour tout schéma $X$ sur $\mathbf{F}_p$, on note $F_X$ le Frobenius absolu. Soit $G$ un schéma en groupes fini et plat sur $\mathcal{O}_K$ et $G_1 = G \times_{\mathrm{Spec}(\mathcal{O}_K)} \mathrm{Spec}(\mathcal{O}_K/p)$, on note $G_1^{(p)}$ le produit fibré: $\begin{array}{ccc} G_1^{(p)} & \to & G_1 \\ \downarrow & & \downarrow \\ S & \xrightarrow{F_S} & S \end{array}$. On a un Frobenius (relatif)

$F: G_1 \to G_1^{(p)}$ (car $G_1$ est un schéma sur $\mathbf{F}_p$) et un Verschiebung $V: G_1^{(p)} \to G_1$ (car $G_1$ est un schéma en groupes fini et plat) tels que $F \circ V = p$ et $V \circ F = p$.

LEMME 3.2.12. *Supposons de plus $G$ tué par $p^n$, on a un isomorphisme $S_n$-linéaire*:

$$S_n \otimes_{(\phi), S_n} \mathrm{Hom}_{(\mathrm{Ab}/\mathcal{O}_K)}(G, \mathcal{O}^{\mathrm{cris}}_{n,\pi}) \xrightarrow{\sim} \mathrm{Hom}_{(\mathrm{Ab}/S)}(G_1^{(p)}, w_* \mathcal{O}_{S/E_n})$$

*tel que le diagramme*:

$$\begin{array}{ccccc}
S_n \otimes_{(\phi), S_n} \mathrm{Hom}_{(\mathrm{Ab}/\mathcal{O}_K)}(G, \mathcal{O}_{n,\pi}^{\mathrm{cris}}) & \xrightarrow{\sim} & \mathrm{Hom}_{(\mathrm{Ab}/S)}(G_1^{(p)}, w_*\mathcal{O}_{S/E_n}) & \xrightarrow{F} & \mathrm{Hom}_{(\mathrm{Ab}/S)}(G_1, w_*\mathcal{O}_{S/E_n}) \\
\| & & & & \| \\
S_n \otimes_{(\phi), S_n} \mathrm{Hom}_{(\mathrm{Ab}/\mathcal{O}_K)}(G, \mathcal{O}_{n,\pi}^{\mathrm{cris}}) & \xrightarrow{\mathrm{Id} \otimes \phi} & \mathrm{Hom}_{(\mathrm{Ab}/\mathcal{O}_K)}(G, \mathcal{O}_{n,\pi}^{\mathrm{cris}}) & \simeq & \mathrm{Hom}_{(\mathrm{Ab}/S)}(G_1, w_*\mathcal{O}_{S/E_n})
\end{array}$$

*commute, où $\phi$ est induit par le Frobenius $\mathcal{O}_{n,\pi}^{\mathrm{cris}} \to \mathcal{O}_{n,\pi}^{\mathrm{cris}}$ et $F$ par le Frobenius $G_1 \to G_1^{(p)}$.*

*Preuve.* Notons, avec [BBM], $\sigma : E_n \to E_n$ le relèvement du Frobenius de $S$ que définit $\phi : S_n \to S_n$ dans (2.1.1), par [BBM, 1.3.5, 4.2.9] et des fonctorialités faciles, on a un diagramme commutatif de $\mathcal{O}_{E_n}$-modules:

$$\begin{array}{ccc}
\sigma^*(\mathcal{H}\mathrm{om}_{S/E_n}(\underline{G}_1, \mathcal{O}_{S/E_n})|_{E_{n,\mathrm{zar}}}) & \longrightarrow & \mathcal{H}\mathrm{om}_{S/E_n}(\underline{G}_1^{(p)}, \mathcal{O}_{S/E_n})|_{E_{n,\mathrm{zar}}} \\
{}^{\text{"}p^n\text{"}} \uparrow \wr & & \wr \uparrow {}^{\text{"}p^n\text{"}} \\
\sigma^*(\mathcal{E}\mathrm{xt}^1_{S/E_n}(\underline{G}_1, \mathcal{O}_{S/E_n})|_{E_{n,\mathrm{zar}}}) & \xrightarrow{\sim} & \mathcal{E}\mathrm{xt}^1_{S/E_n}(\underline{G}_1^{(p)}, \mathcal{O}_{S/E_n})|_{E_{n,\mathrm{zar}}}
\end{array}$$

où la flèche horizontale inférieure est un isomorphisme car $\mathcal{E}\mathrm{xt}^1_{S/E_n}(\underline{G}_1, \mathcal{O}_{S/E_n})$ est un cristal (cf. [BBM, 3.1.3 et 1.3.5]). Compte tenu de [BBM, 2.3.1] et du fait que $E_n$ est un schéma affine, on en déduit un isomorphisme $S_n$-linéaire:

$$S_n \otimes_{(\phi), S_n} \mathrm{Hom}_{S/E_n}(\underline{G}_1, \mathcal{O}_{S/E_n}) \xrightarrow{\sim} \mathrm{Hom}_{S/E_n}(\underline{G}_1^{(p)}, \mathcal{O}_{S/E_n})$$

d'où l'isomorphisme de l'énoncé par (3.2.1) et (3.2.2). La commutativité du diagramme provient du fait que la composée $\mathcal{O}_{n,\pi}^{\mathrm{cris}}(G_1) \to \mathcal{O}_{n,\pi}^{\mathrm{cris}}(G_1^{(p)}) \xrightarrow{F} \mathcal{O}_{n,\pi}^{\mathrm{cris}}(G_1)$ est précisément la définition du Frobenius cristallin. $\square$

LEMME 3.2.13. *Soit $G$ un schéma en groupes fini et plat sur $\mathcal{O}_K$ tué par $p$, on a*:

$$\mathrm{Ker}(\mathrm{Id} \otimes \phi : S_1 \otimes_{(\phi), S_1} \mathrm{Mod}(G) \to \mathrm{Mod}(G)) = \mathrm{Im}(S_1 \otimes_{(\phi), S_1} \mathrm{Fil}^1 \mathrm{Mod}(G)).$$

*Preuve.* On le montre par un calcul local, indépendant de [BBM]. Soit $f \in S_1 \otimes_{(\phi), S_1} \mathrm{Mod}(G)$, il existe des $f_{i,m} \in \mathrm{Mod}(G)$ tels que $f$ s'écrive:

$$f = \sum_{\substack{m \in \mathbf{N} \\ i \in \{0, \ldots, p-1\}}} u^i \gamma_{pm}(u^e) \otimes f_{i,m}.$$

Si $(\mathrm{Id} \otimes \phi)(f) = 0$, c'est que pour tout $\mathfrak{A}$ de $\mathrm{Sp}f(\mathcal{O}_K)_{\mathrm{syn}}$ et tout $x \in G(\mathfrak{A})$:

$$\sum_{\substack{m \in \mathbf{N} \\ i \in \{0, \ldots, p-1\}}} u^i \gamma_{pm}(u^e) \phi(f_{i,m}(x)) = 0$$

dans $\mathcal{O}_{1,\pi}^{\mathrm{cris}}(\mathfrak{A})$. Il suffit donc de montrer $f_{i,m}(x) \in \mathcal{J}_{1,\pi}^{\mathrm{cris}}(\mathfrak{A})$. C'est un problème local; on peut donc supposer $\mathfrak{A} = \mathfrak{A}_\infty$ (2.3.2). La formule (2.3.2,(2)) montre que pour tout $(i,m)$, on a $\phi(f_{i,m}(x)) \in A_\infty \subset \mathcal{O}_{1,\pi}^{\mathrm{cris}}(\mathfrak{A}_\infty)$ et qu'on a une injection:
$$\bigoplus_{\substack{m \in \mathbf{N} \\ i \in \{0,\ldots,p-1\}}} u^i \gamma_{pm}(u^e) A_\infty \hookrightarrow \mathcal{O}_{1,\pi}^{\mathrm{cris}}(\mathfrak{A}_\infty)$$
d'où $\phi(f_{i,m}(x)) = 0$ dans $\mathcal{O}_{1,\pi}^{\mathrm{cris}}(\mathfrak{A}_\infty)$. On conclut alors par la suite exacte: $0 \to \mathcal{J}_{1,\pi}^{\mathrm{cris}}(\mathfrak{A}_\infty) \to \mathcal{O}_{1,\pi}^{\mathrm{cris}}(\mathfrak{A}_\infty) \xrightarrow{\phi} \mathcal{O}_{1,\pi}^{\mathrm{cris}}(\mathfrak{A}_\infty)$ qu'on vérifie à partir de (2.3.2) ($p \geq 3$). □

PROPOSITION 3.2.14. *Soit $H$ un groupe $p$-divisible sur $\mathrm{Spec}(\mathcal{O}_K)$, alors $\phi_1(\mathrm{Fil}^1 \mathrm{Mod}(H(1)))$ engendre $\mathrm{Mod}(H(1))$ sur $S_1$.*

*Preuve.* Soit $H_1 = H \times_{\mathrm{Spec}(\mathcal{O}_K)} S$, le Verschiebung $V : H_1(2)^{(p)} \to H_1(2)$ induit un homomorphisme:
$$\mathrm{Hom}_{(\mathrm{Ab}/S)}(H_1(2), w_* \mathcal{O}_{S/E_2}) \longrightarrow \mathrm{Hom}_{(\mathrm{Ab}/S)}(H_1(2)^{(p)}, w_* \mathcal{O}_{S/E_2}).$$
Par (3.2.1), (3.2.2), (3.2.12) et comme $F \circ V = p$, on déduit un diagramme commutatif:

$$\begin{array}{ccc}
\mathrm{Hom}_{(\mathrm{Ab}/\mathcal{O}_K)}(H(2), \mathcal{O}_{2,\pi}^{\mathrm{cris}}) & \xrightarrow{V} & S_2 \otimes_{(\phi), S_2} \mathrm{Hom}_{(\mathrm{Ab}/\mathcal{O}_K)}(H(2), \mathcal{O}_{2,\pi}^{\mathrm{cris}}) \\
\| & & \downarrow \mathrm{Id} \otimes \phi \\
\mathrm{Hom}_{(\mathrm{Ab}/\mathcal{O}_K)}(H(2), \mathcal{O}_{2,\pi}^{\mathrm{cris}}) & \xrightarrow{p} & \mathrm{Hom}_{(\mathrm{Ab}/\mathcal{O}_K)}(H(2), \mathcal{O}_{2,\pi}^{\mathrm{cris}})
\end{array}$$

et un diagramme analogue avec $\mathcal{O}_{1,\pi}^{\mathrm{cris}}$ (et $H(1)$). Soit
$$f \in \mathrm{Hom}_{(\mathrm{Ab}/\mathcal{O}_K)}(H(1), \mathcal{O}_{1,\pi}^{\mathrm{cris}}),$$
par (3.2.13), $V(f)$ s'écrit: $V(f) = \sum s_i \otimes f_i$ où $s_i \in S_2$ et $f_i$ provient de
$$\mathrm{Hom}_{(\mathrm{Ab}/\mathcal{O}_K)}(H(1), \mathcal{J}_{1,\pi}^{\mathrm{cris}}) = \mathrm{Hom}_{(\mathrm{Ab}/\mathcal{O}_K)}(H(2), \mathcal{J}_{1,\pi}^{\mathrm{cris}}).$$
Par (3.2.11), il existe des relevés $\hat{f}_i$ de $f_i$ dans $\mathrm{Hom}_{(\mathrm{Ab}/\mathcal{O}_K)}(H(2), \mathcal{J}_{2,\pi}^{\mathrm{cris}})$ et l'image de
$$\sum s_i \otimes \hat{f}_i \in S_2 \otimes_{(\phi), S_2} \mathrm{Hom}_{(\mathrm{Ab}/\mathcal{O}_K)}(H(2), \mathcal{O}_{2,\pi}^{\mathrm{cris}})$$
dans $S_1 \otimes_{(\phi), S_1} \mathrm{Hom}_{(\mathrm{Ab}/\mathcal{O}_K)}(H(1), \mathcal{O}_{1,\pi}^{\mathrm{cris}})$ est par construction $V(f)$. Soit $\hat{f}$ un relevé quelconque de $f$ dans $\mathrm{Hom}_{(\mathrm{Ab}/\mathcal{O}_K)}(H(2), \mathcal{O}_{2,\pi}^{\mathrm{cris}})$ (3.2.11), $V(\hat{f}) - \sum s_i \otimes \hat{f}_i$ a donc pour image $0$ dans $S_2 \otimes_{(\phi), S_2} \mathrm{Hom}_{(\mathrm{Ab}/\mathcal{O}_K)}(H(2), \mathcal{O}_{1,\pi}^{\mathrm{cris}})$. Mais on a une suite exacte:
$$\begin{array}{rcl}
0 \to \mathrm{Hom}_{(\mathrm{Ab}/\mathcal{O}_K)}(H(2), \mathcal{O}_{1,\pi}^{\mathrm{cris}}) & \xrightarrow{p} & \mathrm{Hom}_{(\mathrm{Ab}/\mathcal{O}_K)}(H(2), \mathcal{O}_{2,\pi}^{\mathrm{cris}}) \\
& \to & \mathrm{Hom}_{(\mathrm{Ab}/\mathcal{O}_K)}(H(2), \mathcal{O}_{1,\pi}^{\mathrm{cris}}) \to 0
\end{array}$$

(preuve de (3.2.11)), d'où une suite exacte:

$$S_2 \otimes_{(\phi),S_2} \mathrm{Hom}_{(\mathrm{Ab}/\mathcal{O}_K)}(H(2),\mathcal{O}^{\mathrm{cris}}_{1,\pi}) \xrightarrow{p} S_2 \otimes_{(\phi),S_2} \mathrm{Hom}_{(\mathrm{Ab}/\mathcal{O}_K)}(H(2),\mathcal{O}^{\mathrm{cris}}_{2,\pi})$$
$$\to S_2 \otimes_{(\phi),S_2} \mathrm{Hom}_{(\mathrm{Ab}/\mathcal{O}_K)}(H(2),\mathcal{O}^{\mathrm{cris}}_{1,\pi}) \to 0.$$

Donc il existe
$$g \in S_2 \otimes_{(\phi),S_2} \mathrm{Hom}_{(\mathrm{Ab}/\mathcal{O}_K)}(H(2),\mathcal{O}^{\mathrm{cris}}_{1,\pi})$$

tel que $pg = V(\hat{f}) - \sum s_i \otimes \hat{f}_i$ d'où, dans $\mathrm{Hom}_{(\mathrm{Ab}/\mathcal{O}_K)}(H(2),\mathcal{O}^{\mathrm{cris}}_{2,\pi}) \subset \mathcal{O}^{\mathrm{cris}}_{2,\pi}(H(2))$:

$$\begin{aligned}(\mathrm{Id}\otimes\phi)(\sum s_i \otimes \hat{f}_i) &= p(\mathrm{Id}\otimes\phi_1)(\sum s_i \otimes \hat{f}_i)\\ &= (\mathrm{Id}\otimes\phi)(V(\hat{f})) - p(\mathrm{Id}\otimes\phi)(g)\\ &= p\hat{f} - p(\mathrm{Id}\otimes\phi)(g)\end{aligned}$$

qui entraîne, dans $\mathrm{Hom}_{(\mathrm{Ab}/\mathcal{O}_K)}(H(2),\mathcal{O}^{\mathrm{cris}}_{1,\pi}) \subset \mathcal{O}^{\mathrm{cris}}_{1,\pi}(H(2))$ (2.3.3):

$$(\mathrm{Id}\otimes\phi_1)(\sum s_i \otimes f_i) = f - (\mathrm{Id}\otimes\phi)(g).$$

Mais $g$ s'écrit $\sum t_i \otimes g_i$ dans $S_1 \otimes_{(\phi),S_1} \mathrm{Mod}(H(1))$, d'où $f = (\mathrm{Id}\otimes\phi_1)(\sum s_i \otimes f_i + \sum c^{-1} t_i \otimes u^e g_i)$: tout $f$ est dans l'image de $\mathrm{Id}\otimes\phi_1$. □

COROLLAIRE 3.2.15. *Soit $G$ un schéma en groupes fini et plat sur $\mathcal{O}_K$ tué par $p$, alors $\phi_1(\mathrm{Fil}^1\mathrm{Mod}(G))$ engendre $\mathrm{Mod}(G)$ sur $S_1$.*

*Preuve.* En procédant comme en (3.2.7), on trouve un groupe $p$-divisible $H$ sur $\mathrm{Spec}(\mathcal{O}_K)$ et une $\mathcal{O}_K$-immersion fermée ($\mathcal{O}_K$ est local): $G \hookrightarrow H(1)$. En notant $G'$ le conoyau, on a une suite exacte dans ($p$-$\mathrm{Gr}/\mathcal{O}_K$): $0 \to G \to H(1) \to G' \to 0$, d'où un diagramme commutatif où la flèche du haut est surjective par (3.2.6):

$$\begin{array}{ccccc}\mathrm{Hom}_{(\mathrm{Ab}/\mathcal{O}_K)}(H(1),\mathcal{O}^{\mathrm{cris}}_{1,\pi}) & \to & \mathrm{Hom}_{(\mathrm{Ab}/\mathcal{O}_K)}(G,\mathcal{O}^{\mathrm{cris}}_{1,\pi}) & \to & 0\\ \mathrm{Id}\otimes\phi_1 \uparrow & & \uparrow \mathrm{Id}\otimes\phi_1 & & \\ S_1\otimes_{(\phi),S_1}\mathrm{Hom}_{(\mathrm{Ab}/\mathcal{O}_K)}(H(1),\mathcal{J}^{\mathrm{cris}}_{1,\pi}) & \to & S_1\otimes_{(\phi),S_1}\mathrm{Hom}_{(\mathrm{Ab}/\mathcal{O}_K)}(G,\mathcal{J}^{\mathrm{cris}}_{1,\pi}) & & .\end{array}$$

Par (3.2.14), la flèche verticale de gauche est aussi surjective, d'où le résultat. □

*Remarque* 3.2.16. Soit $G$ comme dans (3.2.15) et posons $\mathcal{M} = \mathrm{Mod}(G)$ et $\mathrm{Fil}^1\mathcal{M} = \mathrm{Fil}^1\mathrm{Mod}(G)$, des preuves de (3.2.14) et (3.2.15), on déduit aisément que la composée:

$$\mathcal{M} \xrightarrow[\sim]{(\mathrm{Id}\otimes\phi_1)^{-1}} S_1 \otimes_{(\phi),S_1} \frac{\mathrm{Fil}^1\mathcal{M}}{(u^e\mathrm{Fil}^1\mathcal{M} + \mathrm{Fil}^p S_1\mathcal{M})} \longrightarrow S_1 \otimes_{(\phi),S_1} \mathcal{M}$$

(voir 2.1.2.3) n'est autre que le Verschiebung. Rappelons que le Frobenius $F = \text{Id} \otimes \phi$ (cf. 3.2.12): $S_1 \otimes_{(\phi),S_1} \mathcal{M} \to \mathcal{M}$ s'exprime aussi à partir de $\phi_1$ par $\phi(x) = \frac{1}{\phi_1(v)}\phi_1(vx)$ (2.3.4).

3.3. *L'équivalence de catégories.* En (3.1), on a associé à tout objet $\mathcal{M}$ de $(\text{Mod}/S_1)$ un schéma en groupes $\text{Gr}(\mathcal{M})$ de $(p\text{-Gr}/\mathcal{O}_K)$ annulé par $p$. En (3.2), on a associé à tout objet $G$ de $(p\text{-Gr}/\mathcal{O}_K)$ annulé par $p$ un module $\text{Mod}(G)$ de $(\text{Mod}/S_1)$. On montre ici qu'on a des identifications naturelles (et fonctorielles): $\mathcal{M} \xrightarrow{\sim} \text{Mod}(\text{Gr}(\mathcal{M}))$ et $G \xrightarrow{\sim} \text{Gr}(\text{Mod}(G))$.

LEMME 3.3.1. *Soit $\mathcal{M}' \to \mathcal{M}$ un morphisme d'objets de $(\text{Mod}/S_1)$ tel que le noyau (du morphisme de modules sous-jacent) est contenu dans $\text{Fil}^p S_1 \mathcal{M}'$, alors ce noyau est nul.*

*Preuve.* Soit $\mathcal{K} = \text{Ker}(\phi_1(\text{Fil}^1 \mathcal{M}') \to \phi_1(\text{Fil}^1 \mathcal{M}))$ et $x = \phi_1(y) \in \mathcal{K}$, on a $x \in \text{Fil}^p S_1 \mathcal{M}'$ donc $\overline{x} = 0$ dans $\tilde{\mathcal{M}}'$ (2.1.2.2) donc $\overline{y} \in u^e \text{Fil}^1 \tilde{\mathcal{M}}'$ par (2.1.2.1) donc $y \in u^e \text{Fil}^1 \mathcal{M}' + \text{Fil}^p S_1 \mathcal{M}'$ donc $\phi_1(y) = 0 = x$; i.e. $\mathcal{K} = 0$. Comme $S_1$ est plat sur $k[u^p]/u^{ep}$, on a par (2.1.2.3): $\text{Ker}(\mathcal{M}' \to \mathcal{M}) = S_1 \otimes_{k[u^p]/u^{ep}} \mathcal{K}$ d'où le résultat. □

LEMME 3.3.2. *Soit $f : \mathcal{M}' \to \mathcal{M}$ un morphisme d'objets de $(\text{Mod}/S_1)$, on suppose:* (1) $\text{rg}_{S_1} \mathcal{M}' = \text{rg}_{S_1} \mathcal{M}$ *et* (2) $\text{Ker}(f) \subset \text{Fil}^p S_1 \mathcal{M}'$, *alors $f$ est un isomorphisme.*

*Preuve.* Par (3.3.1) on a en particulier $\phi_1(\text{Fil}^1 \mathcal{M}') \hookrightarrow \phi_1(\text{Fil}^1 \mathcal{M})$. Soit $\tilde{f} : \tilde{\mathcal{M}}' \to \tilde{\mathcal{M}}$ le morphisme dans $(\text{Mod}/\tilde{S}_1)$ déduit de $f$ (2.1.2.2), comme $\phi_1(\text{Fil}^1 \mathcal{M}') \xrightarrow{\sim} \tilde{\phi}_1(\text{Fil}^1 \tilde{\mathcal{M}}')$ et $\phi_1(\text{Fil}^1 \mathcal{M}) \xrightarrow{\sim} \tilde{\phi}_1(\text{Fil}^1 \tilde{\mathcal{M}})$ (même raisonnement que dans la preuve de (3.3.1)), on a $\tilde{\phi}_1(\text{Fil}^1 \tilde{\mathcal{M}}') \hookrightarrow \tilde{\phi}_1(\text{Fil}^1 \tilde{\mathcal{M}})$ d'où par (2.1.2.1) $\tilde{\mathcal{M}}' \hookrightarrow \tilde{\mathcal{M}}$. Comme il s'agit par ailleurs de deux $k$-espaces vectoriels de même dimension, on a $\tilde{\mathcal{M}}' \xrightarrow{\sim} \tilde{\mathcal{M}}$ d'où $\mathcal{M}' \xrightarrow{\sim} \mathcal{M}$ par (2.1.2.2). □

Soit $\mathcal{M}$ un objet de $(\text{Mod}/S_1)$, il est clair qu'on a un morphisme $S_1$-linéaire canonique:

$$\mathcal{M} \longrightarrow \text{Hom}_{(\text{Ab}/\mathcal{O}_K)}\Big(\text{Hom}'_{(\text{Mod}/S_1)}(\mathcal{M}, \mathcal{O}^{\text{cris}}_{1,\pi}), \mathcal{O}^{\text{cris}}_{1,\pi}\Big) \simeq \text{Mod}(\text{Gr}(\mathcal{M}))$$

qui préserve $\text{Fil}^1$ et commute à $\phi_1$: c'est donc un morphisme d'objets de $(\text{Mod}/S_1)$.

LEMME 3.3.3. *Avec les notations ci-dessus, le noyau (sur les modules sous-jacents) du morphisme $\mathcal{M} \to \text{Mod}(\text{Gr}(\mathcal{M}))$ est contenu dans $\text{Fil}^p S_1 \mathcal{M}$.*

*Preuve.* Soit $x \in \mathcal{M}$ tel que, pour tout $\mathfrak{A}$ de $\text{Sp}f(\mathcal{O}_K)_{\text{syn}}$ et tout $f$ de $\text{Hom}'_{(\text{Mod}/S_1)}(\mathcal{M}, \mathcal{O}^{\text{cris}}_{1,\pi}(\mathfrak{A}))$, on ait $f(x) = 0$ dans $\mathcal{O}^{\text{cris}}_{1,\pi}(\mathfrak{A})$. Notons, comme en (3.1), $R_{\mathcal{M}}$ la bigèbre de $\text{Gr}(\mathcal{M})$ et soit $f$ l'image de $Id_{R_{\mathcal{M}}}$ dans l'isomorphisme $\mathrm{Hom}_{\mathcal{O}_K}(R_{\mathcal{M}}, R_{\mathcal{M}}) \xrightarrow{\sim} \mathrm{Hom}'_{(\mathrm{Mod}/S_1)}(\mathcal{M}, \mathcal{O}_{1,\pi}^{\mathrm{cris}}(R_{\mathcal{M}}))$ (3.1.8), on a en particulier $f(x) = 0$ dans $\mathcal{O}_{1,\pi}^{\mathrm{cris}}(R_{\mathcal{M}})$. Soient $(e_1, \ldots, e_d)$ une base adaptée de $\mathcal{M}$ (2.1.2.6), $(s_i)_{1 \leq i \leq d} \in S_1^d$ tels que $x = \sum_{i=1}^{d} s_i e_i$ et $J_1^{[p]}(R_{\mathcal{M}})$ l'idéal de $\mathcal{O}_{1,\pi}^{\mathrm{cris}}(R_{\mathcal{M}})$ engendré par les $\gamma_p(x)$, $x \in \mathcal{J}_{1,\pi}^{\mathrm{cris}}(R_{\mathcal{M}})$, un examen de la preuve de (3.1.5) montre que, dans $\mathcal{O}_{1,\pi}^{\mathrm{cris}}(R_{\mathcal{M}})/J_1^{[p]}(R_{\mathcal{M}})$, on a $f(e_i) = \overline{X}_{i,0} \otimes 1 + \overline{X}_{i,1} \otimes u + \cdots + \overline{X}_{i,p-1} \otimes u^{p-1}$ où les $X_{i,j}$ interviennent dans la présentation de $R_{\mathcal{M}}$ associée à la base $(e_1, \ldots, e_d)$ (cf. 3.1.1). On a donc dans $\mathcal{O}_{1,\pi}^{\mathrm{cris}}(R_{\mathcal{M}})/J_1^{[p]}(R_{\mathcal{M}})$ (et en allégeant les notations): $\sum_{i=1}^{d} \overline{s}_i(\overline{X}_{i,0} + u\overline{X}_{i,1} + \cdots + u^{p-1}\overline{X}_{i,p-1}) = 0$. En procédant comme en (2.3.2), on peut extraire des racines $p^{n^{\mathrm{i\grave{e}mes}}}$ sur $R_{\mathcal{M}}$ et fabriquer un $R_{\mathcal{M},\infty}$ tel que

$$(R_{\mathcal{M},\infty}/pR_{\mathcal{M},\infty})[u]/(u^p - \pi) \hookrightarrow \mathcal{O}_{1,\pi}^{\mathrm{cris}}(R_{\mathcal{M},\infty})/J_1^{[p]}(R_{\mathcal{M},\infty})$$

où $J_1^{[p]}(R_{\mathcal{M},\infty})$ est défini comme pour $R_{\mathcal{M}}$ (l'injection résulte de (2.3.2,(2))). Comme

$$(R_{\mathcal{M}}/pR_{\mathcal{M}})[u]/(u^p - \pi) \hookrightarrow (R_{\mathcal{M},\infty}/pR_{\mathcal{M},\infty})[u]/(u^p - \pi)$$

puisque $R_{\mathcal{M},\infty}/pR_{\mathcal{M},\infty}$ est une limite inductive de recouvrements syntomiques de $R_{\mathcal{M}}/pR_{\mathcal{M}}$, on a:

$$(R_{\mathcal{M}}/pR_{\mathcal{M}})[u]/(u^p - \pi) \hookrightarrow \mathcal{O}_{1,\pi}^{\mathrm{cris}}(R_{\mathcal{M}})/J_1^{[p]}(R_{\mathcal{M}})$$

d'où

$$\sum_{i=1}^{d} \overline{s}_i(\overline{X}_{i,0} + u\overline{X}_{i,1} + \cdots + u^{p-1}\overline{X}_{i,p-1}) = 0$$

dans $(R_{\mathcal{M}}/pR_{\mathcal{M}})[u]/(u^p - \pi)$, ce qui entraîne modulo $u$: $\sum_{i=1}^{d} \overline{s}_i \overline{X}_{i,0} = 0$ dans $R_{\mathcal{M}}/\pi R_{\mathcal{M}}$. En vertu de l'expression de $R_{\mathcal{M}}$ à la fin de la preuve de (3.1.1), on en déduit $\overline{s}_i = 0$ dans $k$ pour tout $i \in \{1, \ldots, d\}$; i.e. $s_i \in uS_1 + \mathrm{Fil}^p S_1$. Soit $s_i' \in S_1$ tel que $s_i - us_i' \in \mathrm{Fil}^p S_1$, on a donc

$$\sum_{i=1}^{d} u\overline{s}_i'(\overline{X}_{i,0} + u\overline{X}_{i,1} + \cdots + u^{p-1}\overline{X}_{i,p-1}) = 0$$

dans $(R_{\mathcal{M}}/pR_{\mathcal{M}})[u]/(u^p - \pi)$ soit

$$\sum_{i=1}^{d} \overline{s}_i'(\overline{X}_{i,0} + u\overline{X}_{i,1} + \cdots + u^{p-1}\overline{X}_{i,p-1}) \in u^{ep-1}(R_{\mathcal{M}}/pR_{\mathcal{M}})[u]/(u^p - \pi)$$

puisque $(R_{\mathcal{M}}/pR_{\mathcal{M}})[u]/(u^p - \pi)$ est libre sur $k[u]/u^{ep}$. A nouveau, on en déduit modulo $u$: $\sum_{i=1}^{d} \overline{s}'_i \overline{X}_{i,0} = 0$ dans $R_{\mathcal{M}}/\pi R_{\mathcal{M}}$; i.e. $\overline{s}'_i = 0$ dans $k$; i.e. $s_i \in u^2 S_1 + \mathrm{Fil}^p S_1$. Une récurrence évidente donne finalement pour tout $i \in \{1, \ldots, d\}$, $s_i \in u^{ep} S_1 + \mathrm{Fil}^p S_1 = \mathrm{Fil}^p S_1$. $\square$

COROLLAIRE 3.3.4. *Pour tout objet $\mathcal{M}$ de $(\mathrm{Mod}/S_1)$, on a $\mathcal{M} \xrightarrow{\sim} \mathrm{Mod}(\mathrm{Gr}(\mathcal{M}))$.*

*Preuve.* Par (3.1.1) et (3.2.5), on a $\mathrm{rg}_{S_1} \mathcal{M} = \mathrm{rg}_{S_1} \mathrm{Mod}(\mathrm{Gr}(\mathcal{M}))$, d'où le résultat par (3.3.3) et (3.3.2). $\square$

LEMME 3.3.5. *Soient $G', G$ deux schémas en groupes finis et plats sur $\mathcal{O}_K$ tués par $p$ et $f : G' \to G$ un morphisme dans $(p\text{-}\mathrm{Gr}/\mathcal{O}_K)$ qui induit un isomorphisme $\mathrm{Mod}(f) : \mathrm{Mod}(G) \xrightarrow{\sim} \mathrm{Mod}(G')$ dans $(\mathrm{Mod}/S_1)$, alors $f$ est un isomorphisme dans $(p\text{-}\mathrm{Gr}/\mathcal{O}_K)$.*

*Preuve.* Soient $G_k = G \times_{\mathrm{Spec}(\mathcal{O}_K)} \mathrm{Spec}(k)$ et $G'_k = G' \times_{\mathrm{Spec}(\mathcal{O}_K)} \mathrm{Spec}(k)$, puisque les bigèbres de $G$ et $G'$ sont complètes pour la topologie $\pi$-adique et sans $\pi$-torsion, il suffit de montrer que le morphisme de schémas en groupes $G'_k \to G_k$ est un isomorphisme. A partir de $\phi_1$, on reconstruit sur $\mathrm{Mod}(G)$ le Frobenius $F = \mathrm{Id} \otimes \phi : S_1 \otimes_{(\phi), S_1} \mathrm{Mod}(G) \to \mathrm{Mod}(G)$ et le Verschiebung $V : \mathrm{Mod}(G) \to S_1 \otimes_{(\phi), S_1} \mathrm{Mod}(G)$ (cf. 3.2.16) et on pose $M(G_k) = \mathrm{Mod}(G) \otimes_{S_1} k$ muni des $F$ et $V$ induits où $S_1 \to k$ est la surjection $\gamma_i(u) \mapsto 0$, $i \geq 1$ (resp. $M(G'_k) = \mathrm{Mod}(G') \otimes_{S_1} k$). L'isomorphisme dans $(\mathrm{Mod}/S_1)$ induit donc un isomorphisme de $k$-espaces vectoriels compatible à $F$ et $V$: $M(G_k) \xrightarrow{\sim} M(G'_k)$. Par ailleurs, $\mathrm{Mod}(G)$ s'identifie à $\mathbf{D}(G_1)(S, E_1) = \mathrm{Ext}^1_{S/E_1}(\underline{G}_1, \mathcal{O}_{S/E_1})$ où $G_1 = G \times_{\mathrm{Spec}(\mathcal{O}_K)} S$ (3.2.4 et 3.2.2), cette identification étant compatible aux opérateurs $F$ et $V$ (resp. avec $G'$). Par [BBM, preuve de 4.3.1], $(M(G_k), F, V)$ et $(M(G'_k), F, V)$ s'identifient alors aux modules de Dieudonné classiques associés à $G_k$ et $G'_k$. Comme le foncteur de Dieudonné classique est pleinement fidèle ([BM, Th. 1]), on en déduit $G'_k \xrightarrow{\sim} G_k$. $\square$

COROLLAIRE 3.3.6. *Pour tout groupe $G$ de $(p\text{-}\mathrm{Gr}/\mathcal{O}_K)$ annulé par $p$, on a $G \xrightarrow{\sim} \mathrm{Gr}(\mathrm{Mod}(G))$.*

*Preuve.* Soit $G$ un objet de $(p\text{-}\mathrm{Gr}/\mathcal{O}_K)$ tué par $p$, il est clair qu'on a un morphisme canonique dans $(p\text{-}\mathrm{Gr}/\mathcal{O}_K)$: $G \to \mathrm{Gr}(\mathrm{Mod}(G))$, d'où un morphisme dans $(\mathrm{Mod}/S_1)$: $\mathrm{Mod}(\mathrm{Gr}(\mathrm{Mod}(G))) \to \mathrm{Mod}(G)$ tel que la composée avec le morphisme canonique en (3.3.3): $\mathrm{Mod}(G) \to \mathrm{Mod}(\mathrm{Gr}(\mathrm{Mod}(G))) \to \mathrm{Mod}(G)$ est l'identité. Comme $\mathrm{Mod}(G) \to \mathrm{Mod}(\mathrm{Gr}(\mathrm{Mod}(G)))$ est un isomorphisme par (3.3.4), la flèche $\mathrm{Mod}(\mathrm{Gr}(\mathrm{Mod}(G))) \to \mathrm{Mod}(G)$ est aussi un isomorphisme dans $(\mathrm{Mod}/S_1)$, d'où le résultat par (3.3.5). $\square$

On a donc démontré:

THÉORÈME 3.3.7. *Supposons $p \neq 2$, le foncteur* Mod *(3.2) établit une antiéquivalence de catégories entre la catégorie des schémas en groupes finis et plats sur $\mathcal{O}_K$ tués par $p$ et la catégorie de modules $(\mathrm{Mod}/S_1)$ (2.1.1). Le foncteur* Gr *(3.1) en est un quasi-inverse.*

## 4. Schémas en groupes annulés par une puissance de $p$

On suppose toujours $p \neq 2$ et on montre qu'un certain faisceau d'extensions est "nul", ce qui permet de faire les dévissages nécessaires pour étendre l'équivalence (3.3.7) au cas des $p$-groupes (finis, plats), puis au cas des groupes $p$-divisibles.

4.1. *Nullité de certains faisceaux d'extensions.* Rappelons qu'une suite $0 \to \mathcal{M}' \to \mathcal{M} \to \mathcal{M}'' \to 0$ dans $'(\mathrm{Mod}/S)$ est dite exacte si la suite de $S$-modules sous-jacente est exacte ainsi que la suite $0 \to \mathrm{Fil}^1\mathcal{M}' \to \mathrm{Fil}^1\mathcal{M} \to \mathrm{Fil}^1\mathcal{M}'' \to 0$ (2.1.1). La catégorie $'(\mathrm{Mod}/S)$ n'est pas abélienne mais, si $\mathcal{M}$ et $\mathcal{N}$ sont dans $'(\mathrm{Mod}/S)$, on définit l'ensemble $\mathrm{Ext}^1_{'(\mathrm{Mod}/S)}(\mathcal{M}, \mathcal{N})$ comme les classes d'extensions $0 \to \mathcal{N} \to \mathcal{E} \to \mathcal{M} \to 0$ dans $'(\mathrm{Mod}/S)$. Si on a un morphisme $\mathcal{M}' \to \mathcal{M}$ dans $'(\mathrm{Mod}/S)$, le produit fibré de $S$-modules $\mathcal{E} \times_{\mathcal{M}} \mathcal{M}'$ muni de $\mathrm{Fil}^1\mathcal{E} \times_{\mathrm{Fil}^1\mathcal{M}} \mathrm{Fil}^1\mathcal{M}'$ et du $\phi_1$ induit est un objet de $'(\mathrm{Mod}/S)$ et on a une suite exacte dans $'(\mathrm{Mod}/S)$: $0 \to \mathcal{N} \to \mathcal{E} \times_{\mathcal{M}} \mathcal{M}' \to \mathcal{M}' \to 0$. Si on a un morphisme $\mathcal{N} \to \mathcal{N}''$ dans $'(\mathrm{Mod}/S)$, la somme amalgamée de $S$-modules $\mathcal{N}'' \oplus_{\mathcal{N}} \mathcal{E}$ munie de $\mathrm{Fil}^1\mathcal{N}'' \oplus_{\mathrm{Fil}^1\mathcal{N}} \mathrm{Fil}^1\mathcal{E}$ (qui est bien un sous-module) et du $\phi_1$ induit est un objet de $'(\mathrm{Mod}/S)$ et on a une suite exacte dans $'(\mathrm{Mod}/S)$: $0 \to \mathcal{N}'' \to \mathcal{N}'' \oplus_{\mathcal{N}} \mathcal{E} \to \mathcal{M} \to 0$. Donc $\mathrm{Ext}^1_{'(\mathrm{Mod}/S)}(\mathcal{M}, \mathcal{N})$ est un foncteur contravariant en $\mathcal{M}$ et covariant en $\mathcal{N}$. On dira qu'une extension est nulle si elle est scindée dans $'(\mathrm{Mod}/S)$. Le lemme suivant, dont la démonstration est laissée en exercice au lecteur, nous suffira:

LEMME 4.1.1. *Soient $0 \to \mathcal{M}' \to \mathcal{M} \to \mathcal{M}'' \to 0$ une suite exacte courte dans $'(\mathrm{Mod}/S)$, $\mathcal{N}$ un objet quelconque de $'(\mathrm{Mod}/S)$, $f$ un élément de $\mathrm{Hom}_{'(\mathrm{Mod}/S)}(\mathcal{M}', \mathcal{N})$ et $(\mathcal{E})$ un élément de $\mathrm{Ext}^1_{'(\mathrm{Mod}/S)}(\mathcal{M}, \mathcal{N})$, alors:*

(1) *l'image de $f$ dans $\mathrm{Ext}^1_{'(\mathrm{Mod}/S)}(\mathcal{M}'', \mathcal{N})$ est nulle si et seulement si $f$ provient de $\mathrm{Hom}_{'(\mathrm{Mod}/S)}(\mathcal{M}, \mathcal{N})$;*

(2) *l'image de $(\mathcal{E})$ dans $\mathrm{Ext}^1_{'(\mathrm{Mod}/S)}(\mathcal{M}', \mathcal{N})$ est nulle si et seulement si $(\mathcal{E})$ provient de $\mathrm{Ext}^1_{'(\mathrm{Mod}/S)}(\mathcal{M}'', \mathcal{N})$.*

On rappelle que pour tout schéma formel $\mathfrak{X}$ de $\mathrm{Sp}f(\mathcal{O}_K)_{\mathrm{syn}}$ la donnée $(\mathcal{O}^{\mathrm{cris}}_{\infty,\pi}(\mathfrak{X}), \mathcal{J}^{\mathrm{cris}}_{\infty,\pi}(\mathfrak{X}), \phi_1)$ est un objet de $'(\mathrm{Mod}/S)$ et qu'on a des suites exactes dans $'(\mathrm{Mod}/S)$: $0 \to \mathcal{O}_{n,\pi}^{\mathrm{cris}}(\mathfrak{X}) \to \mathcal{O}_{\infty,\pi}^{\mathrm{cris}}(\mathfrak{X}) \xrightarrow{p^n} \mathcal{O}_{\infty,\pi}^{\mathrm{cris}}(\mathfrak{X})$ (2.3). Soit $\mathcal{M}$ un objet de $'(\mathrm{Mod}/S)$, on définit $\mathcal{E}\mathrm{xt}^1_{'(\mathrm{Mod}/S)}(\mathcal{M}, \mathcal{O}_{\infty,\pi}^{\mathrm{cris}})$ comme le faisceau sur $\mathrm{Sp}f(\mathcal{O}_K)_{\mathrm{syn}}$ associé au préfaisceau $\mathfrak{X} \mapsto \mathrm{Ext}^1_{'(\mathrm{Mod}/S)}(\mathcal{M}, \mathcal{O}_{\infty,\pi}^{\mathrm{cris}}(\mathfrak{X}))$. Si $\mathcal{M}$ est dans $'(\mathrm{Mod}/S_1)$, on définit de même $\mathcal{E}\mathrm{xt}^1_{'(\mathrm{Mod}/S_1)}(\mathcal{M}, \mathcal{O}_{1,\pi}^{\mathrm{cris}})$ comme le faisceau sur $\mathrm{Sp}f(\mathcal{O}_K)_{\mathrm{syn}}$ associé au préfaisceau $\mathfrak{X} \mapsto \mathrm{Ext}^1_{'(\mathrm{Mod}/S_1)}(\mathcal{M}, \mathcal{O}_{1,\pi}^{\mathrm{cris}}(\mathfrak{X}))$ (il s'agit donc dans ce cas d'extensions tuées par $p$). On écrit $\mathcal{E}\mathrm{xt}^1_{'(\mathrm{Mod}/S)}(\mathcal{M}, \mathcal{O}_{\infty,\pi}^{\mathrm{cris}}) = 0$ (resp. $\mathcal{E}\mathrm{xt}^1_{'(\mathrm{Mod}/S_1)}(\mathcal{M}, \mathcal{O}_{1,\pi}^{\mathrm{cris}}) = 0$) si, quel que soit $\mathfrak{X}$ dans $\mathrm{Sp}f(\mathcal{O}_K)_{\mathrm{syn}}$, toute extension de $\mathrm{Ext}^1_{'(\mathrm{Mod}/S)}(\mathcal{M}, \mathcal{O}_{\infty,\pi}^{\mathrm{cris}}(\mathfrak{X}))$ (resp. $\mathrm{Ext}^1_{'(\mathrm{Mod}/S_1)}(\mathcal{M}, \mathcal{O}_{1,\pi}^{\mathrm{cris}}(\mathfrak{X}))$) se scinde sur un recouvrement de $\mathfrak{X}$ dans $\mathrm{Sp}f(\mathcal{O}_K)_{\mathrm{syn}}$.

LEMME 4.1.2. *Soit $\mathcal{M}$ un objet de* $(\mathrm{Mod}/S_1)$, *si* $\mathcal{E}\mathrm{xt}^1_{'(\mathrm{Mod}/S_1)}(\mathcal{M}, \mathcal{O}_{1,\pi}^{\mathrm{cris}}) = 0$ *alors* $\mathcal{E}\mathrm{xt}^1_{'(\mathrm{Mod}/S)}(\mathcal{M}, \mathcal{O}_{\infty,\pi}^{\mathrm{cris}}) = 0$.

*Preuve*. Il suffit de montrer que pour toute $\mathcal{O}_K$-algèbre $\mathfrak{A}$ comme en (2.2.1) et toute extension $0 \to \mathcal{O}_{\infty,\pi}^{\mathrm{cris}}(\mathfrak{A}) \to \mathcal{E} \to \mathcal{M} \to 0$ dans $'(\mathrm{Mod}/S)$, il existe un recouvrement syntomique de $\mathfrak{A}$ sur lequel l'extension se scinde. Notons $p_\mathcal{E}$ la multiplication par $p$ sur $\mathcal{E}$, les diagrammes du serpent relatifs à la multiplication par $p$ sur les suites exactes de $S$-modules $0 \to \mathcal{O}_{\infty,\pi}^{\mathrm{cris}}(\mathfrak{A}) \to \mathcal{E} \to \mathcal{M} \to 0$ et $0 \to \mathcal{J}_{\infty,\pi}^{\mathrm{cris}}(\mathfrak{A}) \to \mathrm{Fil}^1 \mathcal{E} \to \mathrm{Fil}^1 \mathcal{M} \to 0$ fournissent des suites exactes de $S$-modules:

$$0 \to \mathcal{O}_{1,\pi}^{\mathrm{cris}}(\mathfrak{A}) \to \mathrm{Ker}(p_\mathcal{E}) \to \mathcal{M} \xrightarrow{\tau} \mathcal{O}_{\infty,\pi}^{\mathrm{cris}}(\mathfrak{A})/p\mathcal{O}_{\infty,\pi}^{\mathrm{cris}}(\mathfrak{A}),$$

$$0 \to \mathcal{J}_{1,\pi}^{\mathrm{cris}}(\mathfrak{A}) \to \mathrm{Fil}^1\mathrm{Ker}(p_\mathcal{E}) \to \mathrm{Fil}^1 \mathcal{M} \xrightarrow{\tau'} \mathcal{J}_{\infty,\pi}^{\mathrm{cris}}(\mathfrak{A})/p\mathcal{J}_{\infty,\pi}^{\mathrm{cris}}(\mathfrak{A}).$$

Soient $(x_1, \ldots, x_d)$ une base de $\mathcal{M}$ sur $S_1$ et $(y_1, \ldots, y_d)$ une famille de $\mathrm{Fil}^1\mathcal{M}$ telle que $\mathrm{Fil}^1\mathcal{M} = \sum S_1 y_i + \sum \mathrm{Fil}^1 S_1 x_i$. A cause des suites exactes de faisceaux $0 \to \mathcal{O}_{1,\pi}^{\mathrm{cris}} \to \mathcal{O}_{\infty,\pi}^{\mathrm{cris}} \xrightarrow{p} \mathcal{O}_{\infty,\pi}^{\mathrm{cris}} \to 0$ et $0 \to \mathcal{J}_{1,\pi}^{\mathrm{cris}} \to \mathcal{J}_{\infty,\pi}^{\mathrm{cris}} \xrightarrow{p} \mathcal{J}_{\infty,\pi}^{\mathrm{cris}} \to 0$, il existe un recouvrement $(\mathfrak{A}_j)_{j \in J}$ de $\mathfrak{A}$ dans $\mathrm{Sp}f(\mathcal{O}_K)_{\mathrm{syn}}$ tel que pour tout $i, j$: $\mathrm{res}_{\mathfrak{A}_j}(\tau(x_i)) = 0$ et $\mathrm{res}_{\mathfrak{A}_j}(\tau'(y_i)) = 0$. Soit $\mathcal{E}_j = \mathcal{O}_{\infty,\pi}^{\mathrm{cris}}(\mathfrak{A}_j) \oplus_{\mathcal{O}_{\infty,\pi}^{\mathrm{cris}}(\mathfrak{A})} \mathcal{E}$ (somme amalgamée dans $'(\mathrm{Mod}/S)$), on a donc pour tout $j$ un diagramme commutatif de suites exactes dans $'(\mathrm{Mod}/S)$:

$$\begin{array}{ccccccccc} 0 & \to & \mathcal{O}_{1,\pi}^{\mathrm{cris}}(\mathfrak{A}_j) & \to & \mathrm{Ker}(p_{\mathcal{E}_j}) & \to & \mathcal{M} & \to & 0 \\ & & \downarrow & & \downarrow & & \downarrow & & \\ 0 & \to & \mathcal{O}_{\infty,\pi}^{\mathrm{cris}}(\mathfrak{A}_j) & \to & \mathcal{E}_j & \to & \mathcal{M} & \to & 0. \end{array}$$

Mais par hypothèse pour tout $j$ il existe un recouvrement syntomique $(\mathfrak{A}_{jj'})_{j' \in J'}$ de $\mathfrak{A}_j$ tel que la suite exacte:

$$0 \to \mathcal{O}_{1,\pi}^{\mathrm{cris}}(\mathfrak{A}_{jj'}) \to \mathcal{O}_{1,\pi}^{\mathrm{cris}}(\mathfrak{A}_{jj'}) \oplus_{\mathcal{O}_{1,\pi}^{\mathrm{cris}}(\mathfrak{A}_j)} \mathrm{Ker}(p_{\mathcal{E}_j}) \to \mathcal{M} \to 0$$

est scindée. Comme on a un diagramme commutatif:

$$
\begin{array}{ccccccccc}
0 & \to & \mathcal{O}_{1,\pi}^{\mathrm{cris}}(\mathfrak{A}_{jj'}) & \to & \mathcal{O}_{1,\pi}^{\mathrm{cris}}(\mathfrak{A}_{jj'}) \oplus_{\mathcal{O}_{1,\pi}^{\mathrm{cris}}(\mathfrak{A}_j)} \mathrm{Ker}(p_{\mathcal{E}_j}) & \to & \mathcal{M} & \to & 0 \\
 & & \downarrow & & \downarrow & & \| & & \\
0 & \to & \mathcal{O}_{\infty,\pi}^{\mathrm{cris}}(\mathfrak{A}_{jj'}) & \to & \mathcal{O}_{\infty,\pi}^{\mathrm{cris}}(\mathfrak{A}_{jj'}) \oplus_{\mathcal{O}_{\infty,\pi}^{\mathrm{cris}}(\mathfrak{A}_j)} \mathcal{E}_j & \to & \mathcal{M} & \to & 0,
\end{array}
$$

la suite exacte inférieure est également scindée, d'où le résultat. □

PROPOSITION 4.1.3. *Soit $\mathcal{M}$ un objet de $(\mathrm{Mod}/S_1)$, alors*:

$$\mathcal{E}\mathrm{xt}^1_{'(\mathrm{Mod}/S_1)}(\mathcal{M}, \mathcal{O}_{1,\pi}^{\mathrm{cris}}) = 0.$$

*Preuve.* Soit $\mathfrak{A}$ dans $\mathrm{Sp}f(\mathcal{O}_K)_{\mathrm{syn}}$ et $\mathfrak{A}_\infty$, $\mathcal{O}_{1,\pi}^{\mathrm{cris}}(\mathfrak{A}_\infty)$ comme en (2.3.2), il suffit de montrer que toute extension $0 \to \mathcal{O}_{1,\pi}^{\mathrm{cris}}(\mathfrak{A}_\infty) \to \mathcal{E} \to \mathcal{M} \to 0$ dans $'(\mathrm{Mod}/S_1)$ se scinde sur un recouvrement syntomique (fini) de $\mathfrak{A}_\infty$. Soit $(e_1, \ldots, e_d)$ une base adaptée de $\mathcal{M}$ (2.1.2.6), $(r_1, \ldots, r_d)$ les entiers de $\{0, \ldots, e\}$ tels que $u^{r_i} e_i \in \mathrm{Fil}^1 \mathcal{M}$ et $\mathcal{G}$ l'unique matrice de $\mathrm{GL}_d(S_1)$ telle que:

$$\begin{pmatrix} \phi_1(u^{r_1} e_1) \\ \vdots \\ \phi_1(u^{r_d} e_d) \end{pmatrix} = \mathcal{G} \begin{pmatrix} e_1 \\ \vdots \\ e_d \end{pmatrix}.$$

On va donc chercher à construire une section $s : \mathcal{M} \to \mathcal{E}$ dans $'(\mathrm{Mod}/S_1)$ quitte à autoriser un recouvrement syntomique de $\mathfrak{A}_\infty$. Soient $\hat{e}_j$ des relevés de $e_j$ dans $\mathcal{E}$, on va les corriger pour que $s(e_j) = \hat{e}_j$ soit un morphisme de $'(\mathrm{Mod}/S_1)$.

*Compatibilité aux* $\mathrm{Fil}^1$. Pour tout $j$, il existe $\delta_j \in \mathcal{O}_{1,\pi}^{\mathrm{cris}}(\mathfrak{A}_\infty)$ tel que $u^{r_j} \hat{e}_j + \delta_j \in \mathrm{Fil}^1 \mathcal{E}$. On rappelle la formule (2.3.2):

$$\mathcal{O}_{1,\pi}^{\mathrm{cris}}(\mathfrak{A}_\infty) \simeq \bigoplus_{(m_0,\ldots,m_{s+1}) \in \mathbf{N}^{s+1}} \frac{A_\infty[u - Y_0]}{(u - Y_0)^p} \gamma_{pm_0}(Y_0^e) \gamma_{pm_1}(\psi_1) \cdots$$
$$\cdots \gamma_{pm_s}(\psi_s) \gamma_{pm_{s+1}}(u - Y_0)$$

où, pour alléger, on a noté $Y_0 = X_0^{p^{-1}}$. Donc $\delta_j$ s'écrit $\delta_j = \delta_{j,0} + \delta_{j,1}$ avec $\delta_{j,0} \in A_\infty = \mathfrak{A}_\infty/p$ et $\delta_{j,1} \in \mathcal{J}_{1,\pi}^{\mathrm{cris}}(\mathfrak{A}_\infty)$. Comme $u^{e-r_j} \delta_j \in \mathrm{Fil}^1 \mathcal{E}$ (car $u^e \mathcal{E} \subset \mathrm{Fil}^1 \mathcal{E}$) on a $u^{e-r_j} \delta_{j,0} \in \mathcal{J}_{1,\pi}^{\mathrm{cris}}(\mathfrak{A}_\infty)$; i.e. $X_0^{e-r_j} \delta_{j,0}^p = 0$ dans $A_\infty$ ou encore $\pi^{e-r_j} \hat{\delta}_{j,0}^p = \pi^e x_j$ dans $\mathfrak{A}_\infty$, où $\hat{\delta}_{j,0}$ désigne un relevé quelconque de $\delta_{j,0}$. Comme $\mathfrak{A}_\infty$ est plat sur $\mathcal{O}_K$, donc sans $\pi$-torsion, on déduit $\delta_{j,0}^p = \pi^{r_j} x_j$ dans $A_\infty$. Soit $y_j \in A_\infty$ tel que $y_j^p = x_j$, on a $(\delta_{j,0} - Y_0^{r_j} y_j)^p = 0$ dans $A_\infty$, i.e. $\delta_{j,0} - Y_0^{r_j} y_j \in \mathcal{J}_{1,\pi}^{\mathrm{cris}}(\mathfrak{A}_\infty)$; i.e. $\delta_{j,0} - u^{r_j} y_j \in \mathcal{J}_{1,\pi}^{\mathrm{cris}}(\mathfrak{A}_\infty)$ et, quitte à remplacer $\hat{e}_j$ par $\hat{e}_j + y_j$, on peut supposer $u^{r_j} \hat{e}_j \in \mathrm{Fil}^1 \mathcal{E}$.

*Compatibilité aux $\phi_1$.* Posons:

$$\begin{pmatrix} c_1 \\ \vdots \\ c_d \end{pmatrix} = \begin{pmatrix} \phi_1(u^{r_1}\hat{e}_1) \\ \vdots \\ \phi_1(u^{r_d}\hat{e}_d) \end{pmatrix} - \mathcal{G}\begin{pmatrix} \hat{e}_1 \\ \vdots \\ \hat{e}_d \end{pmatrix}$$

et cherchons des $\delta_j$ dans $\mathcal{O}_{1,\pi}^{\mathrm{cris}}(\mathfrak{A}_\infty)$ tels que $u^{r_j}\delta_j \in \mathcal{J}_{1,\pi}^{\mathrm{cris}}(\mathfrak{A}_\infty)$ et:

$$\begin{pmatrix} \phi_1(u^{r_1}(\hat{e}_1 + \delta_1)) \\ \vdots \\ \phi_1(u^{r_d}(\hat{e}_d + \delta_d)) \end{pmatrix} = \mathcal{G}\begin{pmatrix} \hat{e}_1 + \delta_1 \\ \vdots \\ \hat{e}_d + \delta_d \end{pmatrix},$$

c'est-à-dire solution du système:

$$(*)\qquad \begin{pmatrix} \phi_1(u^{r_1}\delta_1) \\ \vdots \\ \phi_1(u^{r_d}\delta_d) \end{pmatrix} = \mathcal{G}\begin{pmatrix} \delta_1 \\ \vdots \\ \delta_d \end{pmatrix} - \begin{pmatrix} c_1 \\ \vdots \\ c_d \end{pmatrix}.$$

Ecrivons dans $\mathcal{O}_{1,\pi}^{\mathrm{cris}}(\mathfrak{A}_\infty)$ $c_j = c_{j,0} + c_{j,1} + c_{j,2}$ où $c_{j,0} \in A_\infty$, $c_{j,1} \in (u-Y_0)A_\infty$, $c_{j,2} \in \mathcal{J}_{1,\pi}^{cris,[2]}(\mathfrak{A}_\infty) = (\mathcal{J}_{1,\pi}^{\mathrm{cris}}(\mathfrak{A}_\infty))^{[2]}$, il suffit de résoudre séparément les trois systèmes $(*_0)$, $(*_1)$, $(*_2)$ obtenus en remplaçant les $c_j$ respectivement par les $c_{j,0}$, $c_{j,1}$, $c_{j,2}$ dans le système $(*)$ et d'additionner les trois solutions pour obtenir une solution de $(*)$.

*Résolution de $(*_2)$.* C'est le plus facile. Puisque $\phi_1(\mathcal{J}_{1,\pi}^{cris,[2]}(\mathfrak{A}_\infty)) = 0$ (car $p \geq 3$), il est clair que:

$$\begin{pmatrix} \delta_{1,2} \\ \vdots \\ \delta_{d,2} \end{pmatrix} = \mathcal{G}^{-1}\begin{pmatrix} c_{1,2} \\ \vdots \\ c_{d,2} \end{pmatrix}$$

est solution.

*Résolution de $(*_1)$.* Un calcul simple montre que

$$\phi_1(u - Y_0) - Y_0^{p-1}(u - Y_0) \in \mathcal{J}_{1,\pi}^{cris,[2]}(\mathfrak{A}_\infty);$$

on note $\Delta_{1,2}$ cet élément et on cherche des solutions $\delta_{j,1}$ de $(*_1)$ sous la forme:

$$\begin{pmatrix} \delta_{1,1} \\ \vdots \\ \delta_{d,1} \end{pmatrix} = (u - Y_0)\begin{pmatrix} \mu_1 \\ \vdots \\ \mu_d \end{pmatrix} + \Delta_{1,2}\mathcal{G}^{-1}\begin{pmatrix} X_0^{r_1}\mu_1^p \\ \vdots \\ X_0^{r_d}\mu_d^p \end{pmatrix}$$

où $\mu_j \in \mathcal{O}_{1,\pi}^{\mathrm{cris}}(\mathfrak{A}_\infty)$. Un calcul donne:

$$\begin{pmatrix} \phi_1(u^{r_1}\delta_{1,1}) \\ \vdots \\ \phi_1(u^{r_d}\delta_{d,1}) \end{pmatrix} - \mathcal{G}\begin{pmatrix} \delta_{1,1} \\ \vdots \\ \delta_{d,1} \end{pmatrix} = Y_0^{p-1}(u-Y_0)\begin{pmatrix} X_0^{r_1}\mu_1^p \\ \vdots \\ X_0^{r_d}\mu_d^p \end{pmatrix} - (u-Y_0)\mathcal{G}\begin{pmatrix} \mu_1 \\ \vdots \\ \mu_d \end{pmatrix}$$

et il suffit de résoudre, si on pose $c_{j,1} = (u - Y_0)d_j$, $d_j \in A_\infty$:

$$\begin{pmatrix} \mu_1 \\ \vdots \\ \mu_d \end{pmatrix} = \mathcal{G}^{-1} \begin{pmatrix} d_1 \\ \vdots \\ d_d \end{pmatrix} + Y_0^{p-1} \mathcal{G}^{-1} \begin{pmatrix} Y_0^{\mathrm{pr}_1} \mu_1^p \\ \vdots \\ Y_0^{\mathrm{pr}_d} \mu_d^p \end{pmatrix}.$$

En réinjectant cette équation dans les $\mu_j^p$ à droite, puis dans les $\mu_j^{p^2}$ qui apparaissent, etc... et puisque $Y_0^n = 0$ pour $n >> 0$, on voit qu'on construit ainsi de proche en proche une solution (qui est du coup unique).

*Résolution de* $(*_0)$. Les éléments de $S_1$ étant des éléments de $\mathcal{O}_{1,\pi}^{\mathrm{cris}}(\mathfrak{A}_\infty)$, on peut également écrire la matrice $\mathcal{G}$ sous la forme: $\mathcal{G} = \mathcal{G}_0 + \mathcal{G}_1 + \mathcal{G}_2$ où $\mathcal{G}_0 \in \mathrm{GL}_d(A_\infty)$, $\mathcal{G}_1 = (u - Y_0)\mathcal{H}$, $\mathcal{H}$ à coefficients dans $A_\infty$ et $\mathcal{G}_2$ à coefficients dans $\mathcal{J}_{1,\pi}^{cris,[2]}(\mathfrak{A}_\infty)$ (écrire $u = (u - Y_0) + Y_0$ dans $\mathcal{G}$). Un calcul simple montre qu'il existe des $\gamma_j$ dans $A_\infty$ tels que $\phi_1(u^{r_j} - Y_0^{r_j}) - \gamma_j(u - Y_0) \in \mathcal{J}_{1,\pi}^{cris,[2]}(\mathfrak{A}_\infty)$; on note $\Delta_{r_j,2}$ cette différence. On cherche les $\delta_{j,0}$ sous la forme:

$$\begin{pmatrix} \delta_{1,0} \\ \vdots \\ \delta_{d,0} \end{pmatrix} = \begin{pmatrix} Y_0^{e-r_1} \nu_1 \\ \vdots \\ Y_0^{e-r_d} \nu_d \end{pmatrix} + (u - Y_0) \begin{pmatrix} \mu_1 \\ \vdots \\ \mu_d \end{pmatrix}$$
$$+ \mathcal{G}^{-1} \left[ \gamma_p(Y_0^e) \begin{pmatrix} \nu_1^p \\ \vdots \\ \nu_d^p \end{pmatrix} + \begin{pmatrix} \Delta_{r_1,2} Y_0^{p(e-r_1)} \nu_1^p \\ \vdots \\ \Delta_{r_d,2} Y_0^{p(e-r_d)} \nu_d^p \end{pmatrix} + \Delta_{1,2} \begin{pmatrix} Y_0^{\mathrm{pr}_1} \mu_1^p \\ \vdots \\ Y_0^{\mathrm{pr}_d} \mu_d^p \end{pmatrix} \right]$$
$$- \mathcal{G}^{-1} \left[ (u - Y_0)^2 \mathcal{H} \begin{pmatrix} \mu_1 \\ \vdots \\ \mu_d \end{pmatrix} + \mathcal{G}_2 \begin{pmatrix} Y_0^{e-r_1} \nu_1 \\ \vdots \\ Y_0^{e-r_d} \nu_d \end{pmatrix} + (u - Y_0)\mathcal{G}_2 \begin{pmatrix} \mu_1 \\ \vdots \\ \mu_d \end{pmatrix} \right]$$

où $\nu_j, \mu_j \in \mathcal{O}_{1,\pi}^{\mathrm{cris}}(\mathfrak{A}_\infty)$. Posons $w = \phi_1(Y_0^e) - \gamma_p(Y_0^e) \in A_\infty^*$, un calcul donne:

$$\begin{pmatrix} \phi_1(u^{r_1} \delta_{1,0}) \\ \vdots \\ \phi_1(u^{r_d} \delta_{d,0}) \end{pmatrix} - \mathcal{G} \begin{pmatrix} \delta_{1,0} \\ \vdots \\ \delta_{d,0} \end{pmatrix}$$
$$= Y_0^{p-1}(u - Y_0) \begin{pmatrix} Y_0^{\mathrm{pr}_1} \mu_1^p \\ \vdots \\ Y_0^{\mathrm{pr}_d} \mu_d^p \end{pmatrix} - (u - Y_0)\mathcal{G}_0 \begin{pmatrix} \mu_1 \\ \vdots \\ \mu_d \end{pmatrix}$$
$$+ (u - Y_0) \begin{pmatrix} \gamma_1 Y_0^{p(e-r_1)} \nu_1^p \\ \vdots \\ \gamma_d Y_0^{p(e-r_d)} \nu_d^p \end{pmatrix} - (u - Y_0)\mathcal{H} \begin{pmatrix} Y_0^{e-r_1} \nu_1 \\ \vdots \\ Y_0^{e-r_d} \nu_d \end{pmatrix}$$

$$+ w \begin{pmatrix} \nu_1^p \\ \vdots \\ \nu_d^p \end{pmatrix} - \mathcal{G}_0 \begin{pmatrix} Y_0^{e-r_1} \nu_1 \\ \vdots \\ Y_0^{e-r_d} \nu_d \end{pmatrix}$$

et pour résoudre $(*_0)$, il suffit de résoudre successivement les deux systèmes:

$$(*_{0,1}) \qquad w \begin{pmatrix} \nu_1^p \\ \vdots \\ \nu_d^p \end{pmatrix} = \mathcal{G}_0 \begin{pmatrix} Y_0^{e-r_1} \nu_1 \\ \vdots \\ Y_0^{e-r_d} \nu_d \end{pmatrix} - \begin{pmatrix} c_{1,0} \\ \vdots \\ c_{d,0} \end{pmatrix},$$

$(*_{0,2})$
$$\begin{pmatrix} \mu_1 \\ \vdots \\ \mu_d \end{pmatrix} = \mathcal{G}_0^{-1} \left[ \begin{pmatrix} \gamma_1 Y_0^{p(e-r_1)} \nu_1^p \\ \vdots \\ \gamma_d Y_0^{p(e-r_d)} \nu_d^p \end{pmatrix} - \mathcal{H} \begin{pmatrix} Y_0^{e-r_1} \nu_1 \\ \vdots \\ Y_0^{e-r_d} \nu_d \end{pmatrix} \right] + Y_0^{p-1} \mathcal{G}_0^{-1} \begin{pmatrix} Y_0^{\mathrm{pr}_1} \mu_1^p \\ \vdots \\ Y_0^{\mathrm{pr}_d} \mu_d^p \end{pmatrix}.$$

Il est clair que $(*_{0,1})$ se résoud sur un recouvrement syntomique de $A_\infty$ (ou $\mathfrak{A}_\infty$) et $(*_{0,2})$ se résoud alors comme pour $(*_1)$, en remplaçant $\mathfrak{A}_\infty$ par son recouvrement, $\mathcal{G}$ par $\mathcal{G}_0$ et les $d_j$ par le nouveau terme constant entre crochets. $\square$

COROLLAIRE 4.1.4. *Pour tout objet $\mathcal{M}$ de $(\mathrm{Mod}/S_1)$,*

$$\mathcal{E}\mathrm{xt}^1_{\prime(\mathrm{Mod}/S)}(\mathcal{M}, \mathcal{O}^{\mathrm{cris}}_{\infty,\pi}) = 0.$$

4.2. *$p$-groupes finis et plats sur $\mathcal{O}_K$.* Pour $p \neq 2$, on construit une anti-équivalence de catégories entre la catégorie des $p$-groupes finis et plats sur $\mathcal{O}_K$ et la catégorie $(\mathrm{Mod}/S)$ puis entre la sous-catégorie pleine des $p$-groupes dont le noyau de la multiplication par $p^n$ est plat pour tout $n$ et la sous-catégorie pleine $(\mathrm{Mod}\,\mathrm{FI}/S)$ (2.1.1.4) puis entre la catégorie des groupes $p$-divisibles sur $\mathcal{O}_K$ et la catégorie des $S$-modules fortement divisibles.

4.2.1. *$p$-groupes généraux.* A tout objet $\mathcal{M}$ de $(\mathrm{Mod}/S)$, on associe un faisceau $\mathrm{Gr}(\mathcal{M})$ sur $\mathrm{Sp}f(\mathcal{O}_K)_{\mathrm{syn}}$ en posant

$$\mathrm{Gr}(\mathcal{M})(\mathfrak{X}) = \mathrm{Hom}_{\prime(\mathrm{Mod}/S)}(\mathcal{M}, \mathcal{O}^{\mathrm{cris}}_{\infty,\pi}(\mathfrak{X}))$$

et on obtient ainsi un foncteur contravariant:

$$\mathrm{Gr} : (\mathrm{Mod}/S) \to (\mathrm{Ab}/\mathcal{O}_K), \ \mathcal{M} \mapsto \mathrm{Gr}(\mathcal{M}).$$

PROPOSITION 4.2.1.1. *Pour tout $\mathcal{M}$ dans $(\mathrm{Mod}/S)$, $\mathrm{Gr}(\mathcal{M})$ est représentable par un schéma en groupes fini et plat sur $\mathcal{O}_K$, tué par les puissances de $p$ qui annulent $\mathcal{M}$, et le foncteur $\mathrm{Gr} : (\mathrm{Mod}/S) \to (p\text{-}\mathrm{Gr}/\mathcal{O}_K)$ préserve les suites exactes courtes.*

*Preuve.* Soit $\mathcal{M}$ un objet de $(\mathrm{Mod}/S)$ tel qu'il existe une suite exacte dans $'(\mathrm{Mod}/S)$: $0 \to \mathcal{M}'_1 \to \mathcal{M} \to \mathcal{M}''_1 \to 0$ où $\mathcal{M}'_1$ et $\mathcal{M}''_1$ sont dans $(\mathrm{Mod}/S_1)$. Par (4.1.4) et (4.1.1), la suite de faisceaux $0 \to \mathrm{Gr}(\mathcal{M}''_1) \to \mathrm{Gr}(\mathcal{M}) \to \mathrm{Gr}(\mathcal{M}''_1) \to 0$ est exacte dans $(\mathrm{Ab}/\mathcal{O}_K)$ et $\mathcal{E}\mathrm{xt}^1_{'(\mathrm{Mod}/S)}(\mathcal{M}, \mathcal{O}^{\mathrm{cris}}_{\infty,\pi}) = 0$. Comme tout objet de $(\mathrm{Mod}/S)$ s'obtient par extensions successives à partir d'objets de $(\mathrm{Mod}/S_1)$, on en déduit $\mathcal{E}\mathrm{xt}^1_{'(\mathrm{Mod}/S)}(\mathcal{M}, \mathcal{O}^{\mathrm{cris}}_{\infty,\pi}) = 0$ pour tout $\mathcal{M}$ dans $(\mathrm{Mod}/S)$. Le foncteur $\mathrm{Gr} : (\mathrm{Mod}/S) \to (\mathrm{Ab}/\mathcal{O}_K)$ préserve donc les suites exactes courtes. La représentabilité résulte alors de (2.2.4). □

A tout groupe $G$ de $(p\text{-}\mathrm{Gr}/\mathcal{O}_K)$, on associe un objet $\mathrm{Mod}(G)$ de $'(\mathrm{Mod}/S)$ défini par:

- $\mathrm{Mod}(G) = \mathrm{Hom}_{(\mathrm{Ab}/\mathcal{O}_K)}(G, \mathcal{O}^{\mathrm{cris}}_{\infty,\pi})$,
- $\mathrm{Fil}^1\mathrm{Mod}(G) = \mathrm{Hom}_{(\mathrm{Ab}/\mathcal{O}_K)}(G, \mathcal{J}^{\mathrm{cris}}_{\infty,\pi})$,
- $\phi_1 : \mathrm{Fil}^1\mathrm{Mod}(G) \to \mathrm{Mod}(G)$ est induit par $\phi_1 : \mathcal{J}^{\mathrm{cris}}_{\infty,\pi} \to \mathcal{O}^{\mathrm{cris}}_{\infty,\pi}$.

Si $H$ est un groupe $p$-divisible sur $\mathcal{O}_K$, on rappelle que
$$H(n) = \mathrm{Ker}(p^n : H \to H).$$

PROPOSITION 4.2.1.2. *Soit $H$ un groupe $p$-divisible sur $\mathcal{O}_K$, alors $\mathrm{Mod}(H(n))$ est un $S_n$-module libre de type fini et on a des suites exactes dans $'(\mathrm{Mod}/S)$ pour tout $i \in \{1, \ldots, n-1\}$:*
$$0 \to \mathrm{Mod}(H(n-i)) \to \mathrm{Mod}(H(n)) \to \mathrm{Mod}(H(i)) \to 0.$$

*Preuve.* Par (3.2.6) (et la suite exacte $0 \to \mathcal{O}^{\mathrm{cris}}_{n,\pi} \to \mathcal{O}^{\mathrm{cris}}_{\infty,\pi} \xrightarrow{p^n} \mathcal{O}^{\mathrm{cris}}_{\infty,\pi} \to 0$), la suite est exacte sur les $S$-modules sous-jacents, ce qui entraîne facilement $\mathrm{Mod}(H(n-i)) \xrightarrow{\sim} p^i\mathrm{Mod}(H(n))$ (isomorphisme de $S$-modules). On en déduit alors la liberté par un dévissage à partir du cas tué par $p$ en utilisant que $\mathrm{Mod}(H(1))$ est libre sur $S_1$ (3.2.5) ainsi que (4.2.2.1) (on peut aussi utiliser directement [BBM, 3.3.10] et procéder comme en (3.2.5)). Reste à voir la surjectivité des flèches $\mathrm{Fil}^1\mathrm{Mod}(H(n)) \to \mathrm{Fil}^1\mathrm{Mod}(H(i))$. Par récurrence, il suffit de traiter le cas $i = n-1$ et on raisonne comme dans la preuve de (3.2.11) avec le diagramme commutatif:

$$\begin{array}{ccccccccc}
& & 0 & & 0 & & 0 & & \\
& & \downarrow & & \downarrow & & \downarrow & & \\
0 & \to & \mathcal{J}^{\mathrm{cris}}_{1,\pi} & \to & \mathcal{O}^{\mathrm{cris}}_{1,\pi} & \to & \mathcal{O}_1 & \to & 0 \\
& & \downarrow & & \downarrow & & \downarrow & & \\
0 & \to & \mathcal{J}^{\mathrm{cris}}_{n,\pi} & \to & \mathcal{O}^{\mathrm{cris}}_{n,\pi} & \to & \mathcal{O}_n & \to & 0 \\
& & \downarrow & & \downarrow & & \downarrow & & \\
0 & \to & \mathcal{J}^{\mathrm{cris}}_{n-1,\pi} & \to & \mathcal{O}^{\mathrm{cris}}_{n-1,\pi} & \to & \mathcal{O}_{n-1} & \to & 0 \\
& & \downarrow & & \downarrow & & \downarrow & & \\
& & 0 & & 0 & & 0 & &
\end{array}$$

et avec $H(n)$ au lieu de $H(2)$. □

COROLLAIRE 4.2.1.3. *Soit $G$ un schéma en groupes de $(p\text{-Gr}/\mathcal{O}_K)$, alors $\phi_1(\text{Fil}^1\text{Mod}(G))$ engendre $\text{Mod}(G)$ sur $S$.*

*Preuve.* Le résultat est vrai pour $G = H(n)$ par (4.2.1.2), (2.1.1.1) et (3.2.14). On le déduit pour $G$ quelconque comme en (3.2.15) en utilisant (3.2.6). □

Nous aurons besoin du lemme:

LEMME 4.2.1.4. *Soit $\mathcal{M}$ un objet de $(\text{Mod}/S_1)$ et soit $\text{Fil}^{1'}\mathcal{M} \subset \text{Fil}^1\mathcal{M}$ un sous $S_1$-module contenant $\text{Fil}^1 S_1 \mathcal{M}$ et tel que $\phi_1(\text{Fil}^{1'}\mathcal{M})$ engendre encore $\mathcal{M}$ sur $S_1$, alors $\text{Fil}^{1'}\mathcal{M} = \text{Fil}^1\mathcal{M}$.*

*Preuve.* Par (2.1.2.3), on a un isomorphisme de $k[u]/u^e$-modules:

$$\text{Fil}^{1'}\mathcal{M}/(u^e\text{Fil}^{1'}\mathcal{M} + \text{Fil}^p S_1\mathcal{M}) \xrightarrow{\sim} \text{Fil}^1\mathcal{M}/(u^e\text{Fil}^1\mathcal{M} + \text{Fil}^p S_1\mathcal{M})$$

qui entraîne facilement la surjectivité de $\text{Fil}^{1'}\mathcal{M} \to \text{Fil}^1\mathcal{M}$. □

PROPOSITION 4.2.1.5. *Le foncteur $\text{Mod} : (p\text{-Gr}/\mathcal{O}_K) \to {}'(\text{Mod}/S)$ préserve les suites exactes courtes, et tombe donc dans $(\text{Mod}/S)$.*

*Preuve.* Le problème est l'exactitude sur les $\text{Fil}^1$, celle sur les modules découlant de (3.2.6). J'ignore si l'on peut faire plus simple que la preuve qui suit.

*Premier cas.* On part d'une suite exacte: $0 \to G'_1 \to G \to G''' \to 0$ où $G'_1$ est tué par $p$, on a un diagramme commutatif de suites exactes de $S$-modules:

$$\begin{array}{ccccccccc}
0 & \to & \text{Fil}^1\text{Mod}(G''') & \to & \text{Fil}^1\text{Mod}(G) & \to & \text{Fil}^1\text{Mod}(G'_1) & & \\
 & & \downarrow & & \downarrow & & \downarrow & & \\
0 & \to & \text{Mod}(G''') & \to & \text{Mod}(G) & \to & \text{Mod}(G'_1) & \to & 0
\end{array}$$

et on pose $\text{Fil}^{1'}\text{Mod}(G'_1) = \text{Fil}^1\text{Mod}(G)/\text{Fil}^1\text{Mod}(G''') \hookrightarrow \text{Fil}^1\text{Mod}(G'_1)$. De plus, $\text{Fil}^1 S_1 \text{Mod}(G'_1) \subset \text{Fil}^{1'}\text{Mod}(G'_1)$ et $\phi_1(\text{Fil}^{1'}\text{Mod}(G'_1))$ engendre $\text{Mod}(G'_1)$ sur $S_1$ puisque c'est vrai pour $\text{Fil}^1\text{Mod}(G)$ (4.2.1.3). Par (4.2.1.4), on a donc $\text{Fil}^{1'}\text{Mod}(G'_1) = \text{Fil}^1\text{Mod}(G'_1)$ d'où l'exactitude dans ce cas.

*Deuxième cas.* On part d'une suite exacte $0 \to G' \to G \to G''_1 \to 0$ où $G''_1$ est tué par $p$. Soit $n$ tel que $p^n G = 0$, en procédant comme en (3.2.15), mais de façon duale (pour la dualité de Cartier), on trouve un groupe $p$-divisible $H$ sur $\mathcal{O}_K$ et un épimorphisme $H(n) \to G \to 0$. Notons $G^{(3)}$ le noyau de cet épimorphisme et $G^{(4)}$ le noyau de l'épimorphisme $H(n) \to G''_1$. On a un diagramme commutatif de suites exactes dans $(p\text{-Gr}/\mathcal{O}_K)$:

$$\begin{array}{ccccccccc}
& & 0 & & 0 & & & & \\
& & \uparrow & & \uparrow & & & & \\
0 & \to & G' & \to & G & \to & G''_1 & \to & 0 \\
& & \uparrow & & \uparrow & & \| & & \\
0 & \to & G^{(4)} & \to & H(n) & \to & G''_1 & \to & 0 \\
& & \uparrow & & \uparrow & & & & \\
& & G^{(3)} & = & G^{(3)} & & & & \\
& & \uparrow & & \uparrow & & & & \\
& & 0 & & 0 & & & &
\end{array}$$

qui donne un diagramme commutatif de suites exactes de $S$-modules:

$$\begin{array}{ccccccc}
0 & & 0 & & & & \\
\downarrow & & \downarrow & & & & \\
\text{Fil}^1\text{Mod}(G') & \leftarrow & \text{Fil}^1\text{Mod}(G) & \leftarrow & \text{Fil}^1\text{Mod}(G''_1) & \leftarrow & 0 \\
\downarrow & & \downarrow & & \| & & \\
\text{Fil}^1\text{Mod}(G^{(4)}) & \leftarrow & \text{Fil}^1\text{Mod}(H(n)) & \leftarrow & \text{Fil}^1\text{Mod}(G''_1) & \leftarrow & 0 \\
\downarrow & & \downarrow & & & & \\
\text{Fil}^1\text{Mod}(G^{(3)}) & = & \text{Fil}^1\text{Mod}(G^{(3)}), & & & &
\end{array}$$

et si l'on a la surjectivité de $\text{Fil}^1\text{Mod}(H(n)) \to \text{Fil}^1\text{Mod}(G^{(4)})$, une chasse au diagramme facile donne celle de $\text{Fil}^1\text{Mod}(G) \to \text{Fil}^1\text{Mod}(G')$. Comme $pG''_1 = 0$, le morphisme $H(n) \to G''_1$ se factorise par un épimorphisme $H(1) \to G''_1$ dont on note $G^{(5)}$ le noyau, on a un diagramme commutatif de suites exactes dans $(p\text{-Gr}/\mathcal{O}_K)$:

$$\begin{array}{ccccccccc}
& & 0 & & 0 & & & & \\
& & \uparrow & & \uparrow & & & & \\
0 & \to & G^{(5)} & \to & H(1) & \to & G''_1 & \to & 0 \\
& & \uparrow & & \uparrow & & \| & & \\
0 & \to & G^{(4)} & \to & H(n) & \to & G''_1 & \to & 0 \\
& & \uparrow & & \uparrow & & & & \\
& & H(n-1) & = & H(n-1) & & & & \\
& & \uparrow & & \uparrow & & & & \\
& & 0 & & 0 & & & &
\end{array}$$

qui donne un diagramme commutatif de suites exactes de $S$-modules:

$$\begin{array}{ccccccccc}
& & 0 & & 0 & & & & \\
& & \downarrow & & \downarrow & & & & \\
0 & \leftarrow & \mathrm{Fil}^1\mathrm{Mod}(G^{(5)}) & \leftarrow & \mathrm{Fil}^1\mathrm{Mod}(H(1)) & \leftarrow & \mathrm{Fil}^1\mathrm{Mod}(G''_1) & \leftarrow & 0 \\
& & \downarrow & & \downarrow & & \| & & \\
& & \mathrm{Fil}^1\mathrm{Mod}(G^{(4)}) & \leftarrow & \mathrm{Fil}^1\mathrm{Mod}(H(n)) & \leftarrow & \mathrm{Fil}^1\mathrm{Mod}(G''_1) & \leftarrow & 0 \\
& & \downarrow & & \downarrow & & & & \\
& & \mathrm{Fil}^1\mathrm{Mod}(H(n-1)) & = & \mathrm{Fil}^1\mathrm{Mod}(H(n-1)) & & & & \\
& & & & \downarrow & & & & \\
& & & & 0, & & & &
\end{array}$$

où les surjectivités proviennent du premier cas et de (4.2.1.2). Une chasse au diagramme donne alors la surjectivité de $\mathrm{Fil}^1\mathrm{Mod}(H(n)) \to \mathrm{Fil}^1\mathrm{Mod}(G^{(4)})$.

*Cas général.* On part d'une suite exacte $0 \to G' \to G \to G'' \to 0$ quelconque dans $(p\text{-}\mathrm{Gr}/\mathcal{O}_K)$ et on note $G''(1)$ l'adhérence schématique (au sens [Ra, 2.1]) du noyau de la multiplication par $p$ sur $G''$ et $G^{(2)} = G''/G''(1)$. Soit $G^{(3)}$ le produit fibré $G \times_{G''} G''(1)$, comme $G^{(3)}$ s'insère dans une suite exacte $0 \to G' \to G^{(3)} \to G''(1) \to 0$, il est représentable dans $(p\text{-}\mathrm{Gr}/\mathcal{O}_K)$ et on a un diagramme commutatif de suites exactes dans $(p\text{-}\mathrm{Gr}/\mathcal{O}_K)$:

$$\begin{array}{ccccccccc}
& & 0 & & 0 & & & & \\
& & \uparrow & & \uparrow & & & & \\
0 & \to & G''(1) & \to & G'' & \to & G^{(2)} & \to & 0 \\
& & \uparrow & & \uparrow & & \| & & \\
0 & \to & G^{(3)} & \to & G & \to & G^{(2)} & \to & 0 \\
& & \uparrow & & \uparrow & & & & \\
& & G' & = & G' & & & & \\
& & \uparrow & & \uparrow & & & & \\
& & 0 & & 0 & & & &
\end{array}$$

qui donne un diagramme commutatif de suites exactes de $S$-modules:

$$\begin{array}{ccccccccc}
& & 0 & & 0 & & & & \\
& & \downarrow & & \downarrow & & & & \\
& & \mathrm{Fil}^1\mathrm{Mod}(G''(1)) & \leftarrow & \mathrm{Fil}^1\mathrm{Mod}(G'') & \leftarrow & \mathrm{Fil}^1\mathrm{Mod}(G^{(2)}) & \leftarrow & 0 \\
& & \downarrow & & \downarrow & & \| & & \\
& & \mathrm{Fil}^1\mathrm{Mod}(G^{(3)}) & \leftarrow & \mathrm{Fil}^1\mathrm{Mod}(G) & \leftarrow & \mathrm{Fil}^1\mathrm{Mod}(G^{(2)}) & \leftarrow & 0 \\
& & \downarrow & & \downarrow & & & & \\
& & \mathrm{Fil}^1\mathrm{Mod}(G') & = & \mathrm{Fil}^1\mathrm{Mod}(G') & & & & \\
& & \downarrow & & & & & & \\
& & 0, & & & & & &
\end{array}$$

où la surjectivité provient du deuxième cas. Il suffit donc de montrer la surjectivité de $\mathrm{Fil}^1\mathrm{Mod}(G) \to \mathrm{Fil}^1\mathrm{Mod}(G^{(3)})$ et on est ramené à la situation initiale avec $(G', G'')$ remplacés par $(G^{(3)}, G^{(2)})$. On note alors $G^{(2)}(1)$ l'adhérence schématique du noyau de la multiplication par $p$ sur $G^{(2)}$, $G^{(4)} = G^{(2)}/G^{(2)}(1)$ et $G^{(5)} = G \times_{G^{(2)}} G^{(2)}(1)$, etc.: au bout d'un nombre fini d'itérations, on est ramené à une situation type deuxième cas, d'où le résultat. □

THÉORÈME 4.2.1.6. *Supposons* $p \neq 2$, *le foncteur* Mod *établit une anti-équivalence de catégories entre la catégorie* ($p$-Gr/$\mathcal{O}_K$) *et la catégorie* (Mod/$S$). *Le foncteur* Gr *en est un quasi-inverse. De plus, cette* (*anti-*)*équivalence préserve les suites exactes courtes des deux catégories.*

*Preuve.* Par ce qui précède, il reste à montrer que les morphismes canoniques $\mathcal{M} \to \text{Mod}(\text{Gr}(\mathcal{M}))$ et $G \to \text{Gr}(\text{Mod}(G))$ sont des isomorphismes pour tout $\mathcal{M}$ dans (Mod/$S$) et tout $G$ dans ($p$-Gr/$\mathcal{O}_K$). Par (4.2.1.1) et (4.2.1.5), si $0 \to \mathcal{M}' \to \mathcal{M} \to \mathcal{M}'' \to 0$ est une suite exacte dans (Mod/$S$) (resp. $0 \to G' \to G \to G'' \to 0$ une suite exacte dans ($p$-Gr/$\mathcal{O}_K$)), on a un diagramme commutatif de suites exactes:

$$\begin{array}{ccccccccc}
0 & \to & \mathcal{M}' & \to & \mathcal{M} & \to & \mathcal{M}'' & \to & 0 \\
& & \downarrow & & \downarrow & & \downarrow & & \\
0 & \to & \text{Mod}(\text{Gr}(\mathcal{M}')) & \to & \text{Mod}(\text{Gr}(\mathcal{M})) & \to & \text{Mod}(Gr\mathcal{M}'')) & \to & 0
\end{array}$$

$$\left(\text{resp.} \begin{array}{ccccccccc}
0 & \to & G' & \to & G & \to & G'' & \to & 0 \\
& & \downarrow & & \downarrow & & \downarrow & & \\
0 & \to & \text{Gr}(\text{Mod}(G')) & \to & \text{Gr}(\text{Mod}(G)) & \to & \text{Gr}(\text{Mod}(G'')) & \to & 0
\end{array}\right)$$

et par dévissage, on se ramène à montrer ces isomorphismes pour un objet tué par $p$, ce qui est fait dans (3.3). □

4.2.2. *$p$-groupes dont le noyau de $p^n$ est plat pour tout $n$ et groupes $p$-divisibles.* On commence par un lemme d'algèbre linéaire:

LEMME 4.2.2.1. *Soit $\mathcal{M}$ un $S$-module tel qu'on ait des isomorphismes de $S$-modules: $\mathcal{M}/p\mathcal{M} \simeq (S/pS)^d$ et $p\mathcal{M} \simeq \oplus_{i=1}^{d'} S/p^{n_i}S$ $(d, d', n_i \in \mathbf{N}^*)$. Alors on a aussi $\mathcal{M} \simeq \oplus_{j \in J} S/p^{m_j}S$.*

*Preuve.* C'est la même qu'en [Br3, 2.3.1.1] en remplaçant le $S$ de [Br3] par le $S$ "ramifié" du présent article, ce qui ne change pas les arguments. □

LEMME 4.2.2.2. *Soit $G$ un objet de ($p$-Gr/$\mathcal{O}_K$) tel que*

$$G(1) = \text{Ker}(p : G \to G)$$

*est plat et soit $\mathcal{M} = \text{Mod}(G)$, alors $p\text{Fil}^1\mathcal{M} = \text{Fil}^1\mathcal{M} \cap p\mathcal{M}$, $\text{Mod}(G/G(1)) \simeq (p\mathcal{M}, p\text{Fil}^1\mathcal{M}, \phi_1)$ et $\text{Mod}(G(1)) \simeq (\mathcal{M}/p\mathcal{M}, \text{Fil}^1\mathcal{M}/p\text{Fil}^1\mathcal{M}, \phi_1)$.*

*Preuve.* On a une suite exacte dans ($p$-Gr/$\mathcal{O}_K$): $0 \to G(1) \to G \xrightarrow{p} G$ d'où, par (4.2.1.5), une suite exacte dans (Mod/$S$): $\text{Mod}(G) \xrightarrow{p} \text{Mod}(G) \to \text{Mod}(G(1)) \to 0$ qui entraîne $\text{Fil}^1\mathcal{M}/p\text{Fil}^1\mathcal{M} \hookrightarrow \mathcal{M}/p\mathcal{M}$. Le reste est clair. □

COROLLAIRE 4.2.2.3. *Soit $G$ un objet de $(p\text{-}\mathrm{Gr}/\mathcal{O}_K)$ dont le noyau de la multiplication par $p^n$ est plat pour tout $n$, alors $\mathrm{Mod}(G)$ est un objet de $(\mathrm{Mod}\,\mathrm{FI}/S)$.*

*Preuve.* On fait une récurrence sur la (plus petite) puissance $n$ de $p$ qui annule $G$: c'est évidemment vrai pour $n = 1$ par (3.3.7). Par (4.2.2.2) et l'hypothèse de récurrence, $p\mathrm{Mod}(G) \simeq \mathrm{Mod}(G/G(1))$ est dans $(\mathrm{Mod}\,\mathrm{FI}/S)$ et, par le cas $n = 1$, $\mathrm{Mod}(G)/p\mathrm{Mod}(G) \simeq \mathrm{Mod}(G(1))$ est aussi dans $(\mathrm{Mod}\,\mathrm{FI}/S)$, d'où le résultat par (4.2.2.1). □

On rappelle que si $\mathcal{M}$ est un objet de $(\mathrm{Mod}\,\mathrm{FI}/S)$, $(p^r\mathcal{M}, p^r\mathrm{Fil}^1\mathcal{M}, \phi_1)$ et $(\mathcal{M}/p^r\mathcal{M}, \mathrm{Fil}^1\mathcal{M}/p^r\mathrm{Fil}^1\mathcal{M}, \phi_1)$ sont encore des objets de $(\mathrm{Mod}\,\mathrm{FI}/S)$ pour tout entier $r \geq 1$ (cf. 2.1.1.3).

LEMME 4.2.2.4. *Soit $\mathcal{M}$ un objet de $(\mathrm{Mod}\,\mathrm{FI}/S)$ et $G = \mathrm{Gr}(\mathcal{M})$, alors, pour tout entier $r \geq 1$, $G(r) = \mathrm{Ker}(p^r : G \to G)$ est un objet de $(p\text{-}\mathrm{Gr}/\mathcal{O}_K)$ et $\mathrm{Mod}(G(r)) \simeq \mathcal{M}/p^r\mathcal{M}$ dans $(\mathrm{Mod}\,\mathrm{FI}/S)$.*

*Preuve.* La suite $0 \to p^r\mathcal{M} \to \mathcal{M} \to \mathcal{M}/p^r\mathcal{M} \to 0$ est exacte dans $(\mathrm{Mod}/S)$, d'où une suite exacte dans $(p\text{-}\mathrm{Gr}/\mathcal{O}_K)$ par (4.2.1.1):

$$0 \to \mathrm{Gr}(\mathcal{M}/p^r\mathcal{M}) \to \mathrm{Gr}(\mathcal{M}) \to \mathrm{Gr}(p^r\mathcal{M}) \to 0.$$

Soit $\mathcal{M}^{(p^r)}$ le noyau de la multiplication par $p^r$ sur $\mathcal{M}$ muni des $\mathrm{Fil}^1$ et $\phi_1$ induits, la même démonstration qu'en [Br3, 2.3.1.2] à partir de (2.1.1.3) (voir 2.1.1.4) montre que $\mathcal{M}^{(p^r)}$ est encore un objet de $(\mathrm{Mod}\,\mathrm{FI}/S)$, on a donc une suite exacte dans $(\mathrm{Mod}\,\mathrm{FI}/S)$: $0 \to \mathcal{M}^{(p^r)} \to \mathcal{M} \to p^r\mathcal{M} \to 0$ qui induit une injection par (4.2.1.1) $\mathrm{Gr}(p^r\mathcal{M}) \hookrightarrow \mathrm{Gr}(\mathcal{M})$ d'où une suite exacte dans $(p\text{-}\mathrm{Gr}/\mathcal{O}_K)$: $0 \to \mathrm{Gr}(\mathcal{M}/p^r\mathcal{M}) \to \mathrm{Gr}(\mathcal{M}) \xrightarrow{p^r} \mathrm{Gr}(\mathcal{M})$ qui donne le résultat. □

On a donc démontré:

THÉORÈME 4.2.2.5. *Supposons $p \neq 2$, l'anti-équivalence de catégories $(p\text{-}\mathrm{Gr}/\mathcal{O}_K) \xrightarrow{\sim} (\mathrm{Mod}/S)$ (4.2.1.6) induit une anti-équivalence de catégories entre la sous-catégorie pleine de $(p\text{-}\mathrm{Gr}/\mathcal{O}_K)$ formée des schémas en groupes dont le noyau de la multiplication par $p^n$ est plat pour tout $n$ et la sous-catégorie pleine $(\mathrm{Mod}\,\mathrm{FI}/S)$ de $(\mathrm{Mod}/S)$ formée des modules "à Facteurs Invariants".*

*Remarque* 4.2.2.6. Si $e \leq p-2$, tout $p$-groupe fini et plat satisfait la condition ci-dessus de platitude des noyaux ([Ra, 3.3.6]) et la catégorie $(\mathrm{Mod}\,\mathrm{FI}/S)$ se retrouve dans ce cas équivalente à la catégorie $(\mathrm{Mod}/S)$, et même abélienne (voir aussi Remarque 2.1.1.6).

Rappelons ([BBM, 3.3.8]) qu'un groupe de Barsotti-Tate tronqué d'échelon $n \geq 2$ est un objet de $(p\text{-}\mathrm{Gr}/\mathcal{O}_K)$ annulé par $p^n$ tel que pour tout $r \in \{0, \ldots, n\}$, $\mathrm{Ker}(p^{n-r} : G \to G) = \mathrm{Im}(p^r : G \to G)$ (par exemple, si $A$ est un schéma abélien sur $\mathcal{O}_K$, $A[p^n] = \text{Ker}(p^n : A \to A)$ est un groupe de Barsotti-Tate tronqué d'échelon $n$ pour tout $n \geq 2$). Cela entraîne que le noyau de la multiplication par $p^r$ sur $G$ est plat pour tout $r$ et que la bigèbre de $G$ est un $\mathcal{O}_K$-module libre de rang $p^{nd}$ pour un entier $d$ appelé hauteur de $G$. De ce qui précède et de (3.1.1) et (3.2.5), on déduit par un dévissage facile:

COROLLAIRE 4.2.2.7. *Supposons $p \neq 2$, l'anti-équivalence de catégories $(p\text{-Gr}/\mathcal{O}_K) \xrightarrow{\sim} (\text{Mod}/S)$ induit une anti-équivalence de catégories entre la sous-catégorie pleine de $(p\text{-Gr}/\mathcal{O}_K)$ formée des groupes de Barsotti-Tate tronqués d'échelon $n \geq 2$ et de hauteur $d$ et la sous-catégorie pleine de $(\text{Mod}/S)$ formée des $S_n$-modules libres de rang $d$.*

Soit $H = \cup_{n \in \mathbf{N}^*} H(n)$ un groupe $p$-divisible sur $\text{Spec}(\mathcal{O}_K)$ et $\mathcal{M}_n = \text{Mod}(H(n))$. On pose $\mathcal{M} = \varprojlim_n \mathcal{M}_n$, muni de $\text{Fil}^1 \mathcal{M} = \varprojlim_n \text{Fil}^1 \mathcal{M}_n$ et de $\phi_1 : \text{Fil}^1 \mathcal{M} \to \mathcal{M}$ induit par $\phi_1 : \text{Fil}^1 \mathcal{M}_n \to \mathcal{M}_n$.

LEMME 4.2.2.8. *$\mathcal{M}$ est un module fortement divisible* (2.1.1.8).

*Preuve.* De (4.2.2.7), on déduit que $\mathcal{M}_n$ est un $S_n$-module libre de type fini et $\mathcal{M}_{n-1} \simeq \mathcal{M}_n/p^{n-1}\mathcal{M}_n$, d'où le fait que $\mathcal{M}$ est un $S$-module libre de type fini. Si $x \in \mathcal{M}$ est tel que $px \in \text{Fil}^1 \mathcal{M}$, on a pour tout $n \geq 2$, $p\overline{x} \in p\text{Fil}^1 \mathcal{M}_n$ puisque $p\text{Fil}^1 \mathcal{M}_n = p\mathcal{M}_n \cap \text{Fil}^1 \mathcal{M}_n$ (4.2.2.2), d'où $\overline{x} \in \text{Fil}^1 \mathcal{M}_{n-1}$ pour tout $n \in \mathbf{N}^*$; i.e. $x \in \text{Fil}^1 \mathcal{M}$. On en déduit un isomorphisme $\text{Fil}^1 \mathcal{M}/p\text{Fil}^1 \mathcal{M} \xrightarrow{\sim} \text{Fil}^1 \mathcal{M}_1$. Comme $\phi_1(\text{Fil}^1 \mathcal{M}_1)$ engendre $\mathcal{M}_1$ sur $S_1$ par (3.3.7), un dévissage facile utilisant la complétude de $\mathcal{M}$ pour la topologie $p$-adique montre que $\phi_1(\text{Fil}^1 \mathcal{M})$ engendre $\mathcal{M}$ sur $S$. □

THÉORÈME 4.2.2.9. *Supposons $p \neq 2$, le foncteur:*

$$H \mapsto \left( \varprojlim_n \text{Mod}(H(n)), \varprojlim_n \text{Fil}^1 \text{Mod}(H(n)), \varprojlim_n \phi_1 \right)$$

*établit une anti-équivalence de catégories entre la catégorie des groupes $p$-divisibles de hauteur $d$ sur $\text{Spec}(\mathcal{O}_K)$ et la catégorie des modules fortement divisibles de rang $d$. Le foncteur:*

$$\mathcal{M} \mapsto \cup_{n \in \mathbf{N}^*} \text{Gr}(\mathcal{M}/p^n \mathcal{M})$$

*en est un quasi-inverse.*

*Remarque* 4.2.2.10. Lorsque $e = 1$, il y deux manières (équivalentes) de retrouver la classification originale de Fontaine-Laffaille. La première consiste à utiliser [Br3, 2.3.1.2] qui montre que $(\text{Mod}/S) = (\text{Mod FI}/S)$, puis [Br4, 4.4.1] qui montre que le foncteur naturel de la catégorie de Fontaine-Laffaille $\underline{MF}_{tor}^{f,1}$ dans la catégorie $(\text{Mod}/S)$ est une équivalence (l'opérateur $N$ ne jouant aucun rôle dans la preuve). La deuxième consiste à reprendre toutes les démonstrations de cet article en remplaçant la base $E_n = \text{Spec}(W_n\langle u\rangle)$ par la base $\text{Spec}(W_n)$ (ce qui les rend plus simples !). On voit que les deux manières sont équivalentes en passant d'une base à l'autre par la flèche évidente de $W_n$-algèbres: $W_n \to W_n\langle u\rangle$ ce qui revient, du côté modules, à tensoriser par $W_n\langle u\rangle$ et prendre les structures "produit tensoriel".

*Remarque* 4.2.2.11. On peut bien sûr choisir d'autres relèvements du Frobenius sur $S$ que $\phi(u) = u^p$ (qui a été choisi pour simplifier les calculs). Attention cependant qu'un relèvement quelconque ne marche pas toujours (du moins si l'on veut garder des objets de $'(\text{Mod}/S)$). Par exemple, si $e = 1$, $\pi = p$ et $\phi(u) = u^p + p$, on a $\phi_1(u - p) = \gamma_p(u) \notin S^*$ et ce relèvement ne convient pas. Par contre, tout relèvement de la forme $\phi(u) = u^p(1 + ps(u))$ où $s(u) \in S$ marche.

*Remarque* 4.2.2.12. Les amateurs de log-structures pourront s'essayer à rajouter un opérateur de monodromie $N$ sur les objets de $(\text{Mod}/S)$, dans le style de [Br3], et voir si l'on peut étendre les équivalences ci-dessus en remplaçant les catégories de schémas en groupes par des catégories convenables de "log-schémas en groupes" au sens de Kato ([Ka]).

*Remarque* 4.2.2.13. On peut également se demander si ces équivalences ne se "faisceautisent" pas sur des bases plus générales que $\mathcal{O}_K$, par exemple des bases lisses sur $\mathcal{O}_K$. Le cas $\mathcal{O}_K = W$ a été fait avec cette généralité par Faltings ([Fa1, 7.1]).

## 5. Groupes $p$-divisibles et modules faiblement admissibles

On montre que pour $k \subseteq \overline{\mathbf{F}}_p$ et $p \neq 2$ toute représentation cristalline de $\text{Gal}(\overline{K}/K)$ à poids de Hodge-Tate entre 0 et 1 provient d'un groupe $p$-divisible sur $\mathcal{O}_K$ et que certains $(\phi, N)$-modules filtrés faiblement admissibles de Fontaine sont admissibles.

5.1. *Quelques rappels.* On renvoie à [Fo5] pour tout ce qui concerne les représentations semi-stables et cristallines et les $(\phi, N)$-modules filtrés. Rappelons simplement qu'un $(\phi, N)$-module filtré $D$ est dit faiblement admissible si, pour tout sous-module $D'$ de $D$ stable par $\phi$ et $N$ et muni de la filtration induite, le polygone de Hodge de $D'$ est en dessous du polygone de Newton et si pour $D$ ces deux polygones ont de plus mêmes extrémités. Le $(\phi, N)$-module filtré est dit admissible s'il "provient" d'une représentation semi-stable. Colmez et Fontaine viennent de montrer que ces deux notions sont équivalentes ([CF]). Nous redémontrons ci-dessous (différemment) une partie de ce résultat.

On reprend l'anneau $S$ de (2.1.1) en enrichissant ses structures: on le munit de l'unique $W$-dérivation $N$ continue pour la topologie $p$-adique telle que $N(u^i) = -iu^i$ et, pour $r$ entier $\geq 1$, on définit $\text{Fil}^r S$ comme la complétion $p$-adique de l'idéal engendré par les $\gamma_i(E(u)) = \frac{E(u)^i}{i!}$ pour $i \geq r$. On vérifie que $N\phi = p\phi N$, $N(\text{Fil}^r S) \subset \text{Fil}^{r-1} S$ et, pour $0 \leq r \leq p-1$, $\phi(\text{Fil}^r S) \subset p^r S$. On pose $S_{K_0} = S \otimes_W K_0$ et on étend par $K_0$-linéarité (ou semi-linéarité) toutes ces structures à $S_{K_0}$. On définit aussi $\Sigma = W[[u, \frac{u^{ep}}{p}]] \subset S$: $\Sigma$ est stable par $\phi$ et $N$ et on le munit de la filtration induite $\text{Fil}^r \Sigma = \Sigma \cap \text{Fil}^r S$. On pose $\Sigma_{K_0} = K_0 \otimes_W \Sigma$ muni des $\phi$ et $N$ déduits par extension des scalaires et $\text{Fil}^r \Sigma_{K_0} = K_0 \otimes_W \text{Fil}^r \Sigma$.

Soit $D$ un $(\phi, N)$-module filtré tel que $\text{Fil}^0 D_K = D_K$ et $\text{Fil}^{r+1} D_K = 0$ avec $r < p-1$ où $D_K = K \otimes_{K_0} D$. On lui associe un $S_{K_0}$-module filtré $\mathcal{D}$ avec Frobenius $\phi$ et monodromie $N$ de la façon suivante: $\mathcal{D} = S_{K_0} \otimes_{K_0} D$, $\phi = \phi \otimes \phi$, $N = N \otimes \text{Id} + \text{Id} \otimes N$, $\text{Fil}^0 \mathcal{D} = \mathcal{D}$ et, si $i \in \mathbf{N}$, $\text{Fil}^{i+1} \mathcal{D} = \{x \in \mathcal{D} \mid N(x) \in \text{Fil}^i \mathcal{D} \text{ et } f_\pi(x) \in \text{Fil}^{i+1} D_K\}$ où $f_\pi : \mathcal{D} \to D_K$, $s(u) \otimes x \mapsto s(\pi)x$ si $s(u) \in S_{K_0}$ et $x \in D$. On pose également $\mathcal{D}_\Sigma = \Sigma_{K_0} \otimes_{K_0} D \subset S_{K_0} \otimes_{K_0} D = \mathcal{D}$ muni des structures induites par $\mathcal{D}$ ($\phi$, $N$ et filtration).

*Définition* 5.1.1. *Un pseudo-module fortement divisible de* $\mathcal{D}$ (resp. de $\mathcal{D}_\Sigma$) est un sous-$S$-module $\mathcal{M}$ de $\mathcal{D}$ (resp. de $\mathcal{D}_\Sigma$) stable par $\phi$ et $N$ et vérifiant les trois propriétés:

(1) $\mathcal{M}$ est de type fini sur $S$ (resp. sur $\Sigma$),

(2) $K_0 \otimes_W \mathcal{M} \stackrel{\sim}{\to} \mathcal{D}$ (resp. $\mathcal{D}_\Sigma$),

(3) $\phi(\mathcal{M} \cap \text{Fil}^r \mathcal{D}) \subset p^r \mathcal{M}$ et $\frac{\phi}{p^r}(\mathcal{M} \cap \text{Fil}^r \mathcal{D})$ engendre $\mathcal{M}$ sur $S$ (resp. sur $\Sigma$).

Nous étendons maintenant la définition 2.1.1.8 des modules fortement divisibles:

*Définition* 5.1.2. *Un module fortement divisible de* $\mathcal{D}$ (resp. de $\mathcal{D}_\Sigma$) est un pseudo-module fortement divisible de $\mathcal{D}$ (resp. de $\mathcal{D}_\Sigma$) qui est de plus libre en tant que $S$ (resp. $\Sigma$)-module.

C'est bien une généralisation des modules considérés en (2.1.1.8) car:

PROPOSITION 5.1.3. (1) *Tout module fortement divisible* $\mathcal{M}$ *au sens de* (2.1.1.8) *peut être muni d'un unique endomorphisme* $W$-*linéaire* $N$ *tel que* $N(sx) = N(s)x + sN(x)$ ($s \in S$, $x \in \mathcal{M}$), $N\phi_1 = \phi N = \frac{1}{c}\phi_1 \circ (E(u)N)$ *et* $N(\mathcal{M}) \subset \sum_{i \geq 1} \frac{u^i}{q(i)!} \mathcal{M}$ (*cf.* 2.1.1 *pour* $q(i)!$).

(2) *Pour $\mathcal{M}$ comme en (1) muni du $N$ canonique, il existe un unique $\phi$-module filtré $D$ tel que le $S_{K_0}$-module $\mathcal{D}$ associé à $D$ par la recette précédente s'identifie, avec toutes ses structures, au module $\mathcal{M} \otimes_W K_0$.*

*Preuve.* (2) est démontré dans [Br7, 6], on montre donc uniquement (1). Prouvons l'existence de $N$. Soit $(e_1, \ldots, e_d)$ une base de $\mathcal{M}$ sur $S$ telle que $\mathrm{Fil}^1 \mathcal{M} = (\oplus_{i=1}^{d_1} \mathrm{Fil}^1 Se_i) \oplus (\oplus_{i=d_1+1}^{d} Se_i)$ (cf. 2.1.1.9). Notons $x_i = E(u)e_i$ si $0 \leq i \leq d_1$ et $x_i = e_i$ si $d_1 + 1 \leq i \leq d$. Puisque $\mathcal{M}$ est fortement divisible, $(\phi_1(x_i))_{1 \leq i \leq d}$ est encore une base de $\mathcal{M}$ sur $S$. Soit $I$ le complété $p$-adique de l'idéal de $S$ engendré par les $\frac{u^i}{q(i)!}$ pour $i \geq 1$, i.e. le noyau de la surjection $S \to W$ qui envoie $u$ et ses puissances divisées sur $0$. Définissons $N_0 : \mathcal{M} \to \mathcal{M}$ comme l'unique dérivation telle que $N_0(\phi_1(x_i)) = 0$ pour tout $i$ (i.e., $N_0(\sum s_i \phi_1(x_i)) = \sum N(s_i) \phi_1(x_i)$). On a clairement $N_0(\mathcal{M}) \subset I\mathcal{M}$. Définissons par récurrence $N_n : \mathcal{M} \to \mathcal{M}$ pour $n \geq 1$ comme l'unique dérivation telle que $N_n(\phi_1(x_i)) = \phi(N_{n-1}(x_i))$ pour tout $i$. On vérifie (par récurrence) que $(N_n - N_{n-1})(\mathcal{M}) \subset \phi^n(I)\mathcal{M}$ où $\phi^n(I)$ est le complété $p$-adique de l'idéal de $S$ engendré par les $\frac{u^{ip^n}}{q(i)!}$ pour $i \geq 1$. Comme $\phi^n(I) \subset p^{i(n)} S$ avec $i(n) \to +\infty$ quand $n \to +\infty$, $N_n$ converge pour la topologie $p$-adique vers une dérivation $N : \mathcal{M} \to \mathcal{M}$ telle que $N(\mathcal{M}) \subset I\mathcal{M}$ et $N\phi_1(x_i) = \phi(N(x_i))$ pour tout $i$, i.e. $N\phi_1 = \phi N$. L'unicité d'une telle dérivation est facile et est laissée au lecteur (s'il y en a deux, considérer leur différence). □

Ainsi, les modules de (2.1.1.8) correspondent au cas particulier $r = 1$ et $N = 0$ sur $D$ dans la définition (5.1.2). Si on pense à $\mathcal{M}$ comme provenant d'un groupe $p$-divisible sur $\mathcal{O}_K$, $D$ n'est autre que le module de Dieudonné de ce groupe $p$-divisible. Dans [Br5], il est démontré le résultat suivant:

THÉORÈME 5.1.4 ([Br5, 2.3.2.5, 3.2.1.5, 3.2.4.8]). *Soient $D$, $\mathcal{D}$ et $\mathcal{D}_\Sigma$ comme au début de cette section. Si $D$ est faiblement admissible, alors $\mathcal{D}$ (resp. $\mathcal{D}_\Sigma$) contient un pseudo-module fortement divisible. Si $\mathcal{D}$ (resp. $\mathcal{D}_\Sigma$) contient un module fortement divisible, alors $D$ est admissible.*

Je conjecture que $\mathcal{D}$ et $\mathcal{D}_\Sigma$ contiennent toujours un module fortement divisible dès que $D$ est faiblement admissible et $\mathrm{Fil}^{p-1} D_K = 0$ (et $k$ seulement supposé parfait). C'est un théorème si $\mathrm{Fil}^{r+1} D_K = 0$ où $er < p - 1$ ([Br5, 1.3]) et aussi si $\mathrm{Fil}^2 D_K = 0$ et $k \subseteq \overline{\mathbf{F}}_p$ comme conséquence de (5.3.2) et (4.2.2.9).

5.2. *Modules faiblement admissibles et modules admissibles.* Soit $D$ un $(\phi, N)$-module filtré faiblement admissible tel que $\mathrm{Fil}^0 D_K = D_K$ et $\mathrm{Fil}^{r+1} D_K = 0$ avec $r < p - 1$. Soient $\mathcal{D}$ (resp. $\mathcal{D}_\Sigma$) le $S_{K_0}$-module (resp. $\Sigma_{K_0}$-module) associé à $D$, $M$ un réseau de $D$ stable par $N$ et $\mathcal{M}_0 = \Sigma \otimes_W M$. Définissons par récurrence pour $i$ entier positif ou nul $\mathcal{M}_{i+1}$ comme le sous-$\Sigma$-module de $\mathcal{D}_\Sigma$ engendré par $\frac{\phi}{p^r}(\mathcal{M}_i \cap \operatorname{Fil}^r \mathcal{D}_\Sigma)$. Pour tout $i$, $\mathcal{M}_{i+1}$ est stable par $N$ car, si $x \in \mathcal{M}_i \cap \operatorname{Fil}^r \mathcal{D}_\Sigma$, on a $E(u)N(x) \in \mathcal{M}_i \cap \operatorname{Fil}^r \mathcal{D}_\Sigma$ et:

$$N\frac{\phi}{p^r}(x) = \frac{p}{\phi(E(u))}\frac{\phi}{p^r}(E(u)N(x)) \in \mathcal{M}_{i+1} \qquad (\frac{p}{\phi(E(u))} = \frac{1}{c} \in \Sigma^*).$$

Par ([Br5, 3.2.3.2]), $\mathcal{M}_{i+1}$ est libre sur $\Sigma$ et $K_0 \otimes_W \mathcal{M}_{i+1} \overset{\sim}{\to} \mathcal{D}_\Sigma$.

LEMME 5.2.1. *Supposons $k \subseteq \overline{\mathbf{F}}_p$. Avec les notations précédentes, l'ensemble $\{S \otimes_\Sigma \mathcal{M}_i, i \in \mathbf{N}\}$ est fini.*

*Preuve.* Soit $\mathcal{M}$ un pseudo-module fortement divisible dans $\mathcal{D}_\Sigma$ (cf. 5.1.4), quitte à faire une homothétie sur $\mathcal{M}$, on peut supposer $p^{n_0}\mathcal{M} \subset \mathcal{M}_0 \subset \mathcal{M}$ pour un $n_0 >> 0$. Comme $\mathcal{M}$ est engendré par $\frac{\phi}{p^r}(\mathcal{M} \cap \operatorname{Fil}^r \mathcal{D}_\Sigma)$, on a encore pour tout $i$: $p^{n_0}\mathcal{M} \subset \mathcal{M}_i \subset \mathcal{M}$. Soit $\mathcal{N}$ le sous-$S$-module de $\mathcal{D}$ engendré par $\mathcal{M}$, il suffit de montrer que l'ensemble $\{\overline{\mathcal{M}}_i, i \in \mathbf{N}\}$ est fini, où $\overline{\mathcal{M}}_i$ désigne l'image de $\mathcal{M}_i$ dans $\mathcal{N}/p^{n_0}\mathcal{N}$. Remarquons que si l'on munit $\mathcal{N}/p^{n_0}\mathcal{N}$ des $\operatorname{Fil}^r$ et $\phi_r$ quotients, on a $\overline{\mathcal{M}}_{i+1}$ engendré par $\phi_r(\overline{\mathcal{M}}_i \cap \operatorname{Fil}^r(\mathcal{N}/p^{n_0}\mathcal{N}))$. Si $k$ est fini, le résultat est clair puisque les $\overline{\mathcal{M}}_i$ sont tous des sous-modules de l'image de $\mathcal{M} \otimes_\Sigma \Sigma/p^{n_0}\Sigma$ dans $\mathcal{N}/p^{n_0}\mathcal{N}$, image qui n'a qu'un nombre fini d'éléments. Si $k \subset \overline{\mathbf{F}}_p$ n'est pas fini, on peut toujours se ramener au cas d'un corps fini en remarquant que dans $\mathcal{N}/p^{n_0}\mathcal{N}$, $\operatorname{Fil}^r$, $\phi_r$ et $\overline{\mathcal{M}}_0$ se définissent à partir d'un nombre fini d'éléments de $\mathcal{N}/p^{n_0}\mathcal{N}$, ne faisant donc intervenir qu'un nombre fini d'éléments de $\overline{\mathbf{F}}_p$, c'est-à-dire une extension finie (suffisamment grande) de $\mathbf{F}_p$. $\square$

*Exemple* 5.2.2. Supposons $p \neq 2$ et soient $K_0$ une extension finie non ramifiée de $\mathbf{Q}_p$ et $K = K_0[\pi]$ où $\pi^{p+2} = p$. Soit $D = K_0 e_1 \oplus K_0 e_2$ avec $N(e_1) = N(e_2) = 0$, $\phi(e_1) = e_2$, $\phi(e_2) = pe_1$, $\operatorname{Fil}^0 D_K = D_K$, $\operatorname{Fil}^1 D_K = K(\pi e_1 + e_2)$ et $\operatorname{Fil}^2 D_K = 0$: on vérifie que $D$ est bien faiblement admissible et que $\operatorname{Fil}^1 \mathcal{D}_\Sigma = \Sigma_{K_0}(ue_1 + e_2) + (u^{p+2} - p)\mathcal{D}_\Sigma$. On choisit $M = We_1 \oplus We_2$, donc $S \otimes_\Sigma \mathcal{M}_0 = Se_1 \oplus Se_2$. L'algorithme ci-dessus donne:

$$S \otimes_\Sigma \mathcal{M}_1 = S(e_1 + \frac{u^p}{p}e_2) \oplus Se_2,$$

$$S \otimes_\Sigma \mathcal{M}_2 = S(e_1 + \frac{u^{2p}}{pw_2}e_2) \oplus Se_2 \text{ où } w_2 = -1 + u^p + \frac{u^{p(p+2)}}{p} \in S^*,$$

$$S \otimes_\Sigma \mathcal{M}_3 = S(e_1 + \frac{u^p}{pw_3}e_2) \oplus Se_2 \text{ où } w_3 = 1 - \frac{u^{p(p-1)}}{u^{p^2}-1} + \frac{u^{p(2p+1)}}{p(u^{p^2}-1)} \in S^*,$$

$$S \otimes_\Sigma \mathcal{M}_4 = S(e_1 + \frac{u^{2p}}{pw_2}e_2) \oplus Se_2 = S \otimes_\Sigma \mathcal{M}_2.$$

Un calcul montre que $S \otimes_\Sigma \mathcal{M}_2 + S \otimes_\Sigma \mathcal{M}_3 = Se_1 + Se_2 + S\frac{u^p}{p}e_2 \subset \mathcal{D}$ est un pseudo-module fortement divisible. Il n'est clairement pas libre.

*Remarque* 5.2.3. Le lemme (5.2.1) est trivialement faux (en général) si on enlève l'hypothèse $k \subseteq \overline{\mathbf{F}}_p$.

THÉORÈME 5.2.4. *Supposons $k \subseteq \overline{\mathbf{F}}_p$ et soit $D$ un $(\phi, N)$-module filtré faiblement admissible tel que* $\mathrm{Fil}^0(D \otimes_{K_0} K) = D \otimes_{K_0} K$ *et* $\mathrm{Fil}^{p-1}(D \otimes_{K_0} K) = 0$, *alors $D$ est admissible.*

*Preuve.* Soit $\mathcal{N}_i = S \otimes_\Sigma \mathcal{M}_i$. Comme $S = \Sigma + \mathrm{Fil}^r S$ ([Br5, 3.2.1.1]) et $\phi(\mathrm{Fil}^r S) \subset p^r S$, on voit que $\mathcal{N}_{i+1}$ est encore le sous-$S$-module de $\mathcal{D}$ engendré par $\frac{\phi}{p^r}(\mathcal{N}_i \cap \mathrm{Fil}^r \mathcal{D})$. Par (5.2.1), l'algorithme $\mathcal{N}_i \to \mathcal{N}_{i+1}$ se résume donc à $\mathcal{N}_1 \to \mathcal{N}_2 \to \cdots \to \mathcal{N}_C \to \mathcal{N}_{i_0}$ pour un $C >> 0$ et un $i_0 \in \{1, \ldots, C\}$. Quitte à partir de $\mathcal{N}_{i_0}$ au lieu de $\mathcal{N}_1$ et à changer la numérotation, on a donc $\mathcal{N}_1 \to \mathcal{N}_2 \to \cdots \to \mathcal{N}_C \to \mathcal{N}_1$. Soit $D^{(C)} = D \oplus D \oplus \cdots \oplus D$ ($C$ fois) muni de $\mathrm{Fil}^i D_K^{(C)} = \mathrm{Fil}^i D_K \oplus \mathrm{Fil}^i D_K \oplus \cdots \oplus \mathrm{Fil}^i D_K$, de $N$ défini par

$$N(d_1 \oplus \cdots \oplus d_C) = N(d_1) \oplus \cdots \oplus N(d_C)$$

et de $\phi$ défini par

$$\phi(d_1 \oplus d_2 \oplus \cdots \oplus d_C) = \phi(d_C) \oplus \phi(d_1) \oplus \cdots \oplus \phi(d_{C-1}).$$

On vérifie que $N\phi = p\phi N$ de sorte que $D^{(C)}$ est encore un $(\phi, N)$-module filtré avec $\mathrm{Fil}^0 D_K^{(C)} = D_K^{(C)}$ et $\mathrm{Fil}^{r+1} D_K^{(C)} = 0$. Soit $\mathcal{D}^{(C)}$ le $S_{K_0}$-module filtré associé à $D^{(C)}$. Comme la filtration sur $\mathcal{D}$ ne dépend que de $N$ et de la filtration sur $D$, on voit que si on oublie le Frobenius: $\mathcal{D}^{(C)} \simeq \mathcal{D} \oplus \cdots \oplus \mathcal{D}$ en tant que $(N, S_{K_0})$-module filtré. Il est alors clair que $\mathcal{N}_1 \oplus \cdots \oplus \mathcal{N}_C$ est un $S$-module fortement divisible de $\mathcal{D}^{(C)}$ et donc que $D^{(C)}$ est admissible par (5.1.4). Par ailleurs, on a une surjection $s : D^{(C)} \to D$ de $(\phi, N)$-modules filtrés définie par $s(d_1 \oplus \cdots \oplus d_C) = d_1 + \cdots + d_C$. Comme $D$ est faiblement admissible, le noyau de cette surjection est faiblement admissible ([Fo5, 4.4.4.iii]), donc admissible ([Fo5, 5.6.7.vii]), d'où on déduit que $D$ est admissible. □

5.3. *Groupes p-divisibles et modules faiblement admissibles.* Rappelons que $A_\mathrm{cris} = \varprojlim W_n(\mathcal{O}_{\overline{K}}/p\mathcal{O}_{\overline{K}})^{\mathrm{DP}}$ (cf. [FM, I.1.3]), $B_\mathrm{cris}^+ = A_\mathrm{cris} \otimes_W K_0$ et qu'on a une injection canonique $K \otimes_{K_0} B_\mathrm{cris}^+ \hookrightarrow B_{dR}^+$ ce qui permet de munir $K \otimes_{K_0} B_\mathrm{cris}^+$ de la filtration induite par $B_{dR}^+$ (voir [Fo6, 4.1]).

LEMME 5.3.1. *Soit $D$ un $\phi$-module filtré faiblement admissible tel que $\mathrm{Fil}^0 D_K = D_K$ et $\mathrm{Fil}^2 D_K = 0$. Soit $\mathcal{D}$ son $S_{K_0}$-module associé et supposons que $\mathcal{D}$ contienne un module fortement divisible $\mathcal{M}$. Soit $H$ le groupe p-divisible associé à $\mathcal{M}$ par (4.2.2.9) (en oubliant l'opérateur $N$), alors $T_p H \otimes_{\mathbf{Z}_p} \mathbf{Q}_p \simeq \mathrm{Hom}_{(\phi, \mathrm{Fil}^1)}(D, B_\mathrm{cris}^+)$ où l'indice signifie qu'on prend les applications* $K_0$-*linéaires qui commutent au Frobenius et respectent le* $\mathrm{Fil}^1$ *après extension des scalaires à* $K$.

*Preuve.* Soit $\widehat{A_{\mathrm{cris}}} = \varprojlim \mathcal{O}_{n,\pi}^{\mathrm{cris}}(\mathcal{O}_{\overline{K}})$. En utilisant la formule précédente pour $A_{\mathrm{cris}}$ et un argument similaire à la preuve de (2.3.2), on voit que $\widehat{A_{\mathrm{cris}}}$ s'identifie au complété pour la topologie $p$-adique de $A_{\mathrm{cris}}[\frac{(u-\underline{\pi})^i}{i!}]_{i\in\mathbf{N}}$ où $\underline{\pi}$ est un élément de $A_{\mathrm{cris}}$ fabriqué à partir d'un système compatible de racines $p^{n^{\mathrm{ièmes}}}$ de $\pi$ (cf. [Br5, 2.2.2]). Par ailleurs, on a un plongement naturel $A_{\mathrm{cris}}$-linéaire de $\widehat{A_{\mathrm{cris}}}$ dans l'anneau $\widehat{A_{st}}$ de [Br5, 2.2.2] en envoyant $\frac{(u-\underline{\pi})^i}{i!}$ sur $\frac{(u-\underline{\pi})^i}{i!}$ (!). On s'en sert pour munir $\widehat{A_{\mathrm{cris}}}$ de la filtration et du Frobenius induits. En particulier, $\frac{(u-\underline{\pi})^i}{i!} \in \mathrm{Fil}^i \widehat{A_{\mathrm{cris}}}$. Si $H = (H(n))_{n\in\mathbf{N}}$, on a par construction du foncteur en (4.2.2.9) et (4.2.1):

$$\begin{aligned} T_p H = \varprojlim H(n)(\mathcal{O}_{\overline{K}}) &= \varprojlim \mathrm{Hom}_{(\phi_1, \mathrm{Fil}^1)}(\mathcal{M}/p^n\mathcal{M}, \mathcal{O}_{n,\pi}^{\mathrm{cris}}(\mathcal{O}_{\overline{K}})) \\ &\simeq \mathrm{Hom}_{(\phi_1, \mathrm{Fil}^1)}(\mathcal{M}, \widehat{A_{\mathrm{cris}}}) \end{aligned}$$

où l'indice signifie qu'on prend les applications $S$-linéaires qui commutent à $\phi_1$ et respectent le $\mathrm{Fil}^1$. Donc $T_p H \otimes_{\mathbf{Z}_p} \mathbf{Q}_p \simeq \mathrm{Hom}_{(\phi, \mathrm{Fil}^1)}(\mathcal{D}, \widehat{A_{\mathrm{cris}}}[\frac{1}{p}])$. En utilisant la bijectivité du Frobenius sur $\mathcal{D}$ et le fait que $\phi(u - \underline{\pi}) = u^p - \underline{\pi}^p$, on vérifie que se donner un élément de $\mathrm{Hom}_{(\phi, \mathrm{Fil}^1)}(\mathcal{D}, \widehat{A_{\mathrm{cris}}}[\frac{1}{p}])$, c'est se donner une application $f$ de $D$ dans $B_{\mathrm{cris}}^+$ $K_0$-linéaire et compatible au Frobenius telle que l'application $S_{K_0}$-linéaire qu'on en déduit de $\mathcal{D} = S_{K_0} \otimes_{K_0} D$ dans $\widehat{A_{\mathrm{cris}}}[\frac{1}{p}]$ est compatible au $\mathrm{Fil}^1$. Mais on a un diagramme commutatif:

$$\begin{array}{ccccc} S_{K_0} \otimes_{K_0} D & \xrightarrow{\mathrm{Id}\otimes f} & S_{K_0} \otimes_{K_0} A_{\mathrm{cris}}[\frac{1}{p}] & \hookrightarrow & \widehat{A_{\mathrm{cris}}}[\frac{1}{p}] \\ f_\pi \downarrow & & \downarrow f_\pi & & \downarrow f_\pi \\ K \otimes_{K_0} D & \xrightarrow{\mathrm{Id}\otimes f} & K \otimes_{K_0} B_{\mathrm{cris}}^+ & \hookrightarrow & B_{dR}^+ \end{array}$$

où les flèches verticales envoient $u$ sur $\pi$ et on vérifie que $x \in \mathrm{Fil}^1(\widehat{A_{\mathrm{cris}}}[\frac{1}{p}])$ (resp. $x \in \mathrm{Fil}^1(S_{K_0} \otimes_{K_0} D)$) si et seulement si $f_\pi(x) \in \mathrm{Fil}^1 B_{dR}^+$ (resp. $f_\pi(x) \in \mathrm{Fil}^1(K \otimes_{K_0} D)$). Il est donc équivalent d'être compatible au $\mathrm{Fil}^1$ en haut ou en bas, d'où le résultat. □

THÉORÈME 5.3.2. *Supposons* $p \neq 2$, $k \subseteq \overline{\mathbf{F}}_p$ *et soit* $V$ *une représentation* $p$-*adique cristalline de* $\mathrm{Gal}(\overline{K}/K)$ *à poids de Hodge-Tate entre* $0$ *et* $1$. *Alors il existe un groupe* $p$-*divisible* $H$ *sur* $\mathcal{O}_K$ *tel que* $V \simeq T_p H \otimes_{\mathbf{Z}_p} \mathbf{Q}_p$ *où* $T_p H$ *est le module de Tate de* $H$.

*Preuve.* Soit $V$ une représentation cristalline comme dans l'énoncé et $D = \mathrm{Hom}_{Gal}(V, B_{\mathrm{cris}}^+)$ son module filtré admissible associé. On reprend la preuve du théorème (5.2.4) avec ce $D$: d'après (5.1.4) le module $D^{(C)}$ est admissible, d'après (4.2.2.9) le module fortement divisible $\mathcal{N}_1 \oplus \cdots \oplus \mathcal{N}_C$ correspond à un groupe $p$-divisible sur $\mathcal{O}_K$, et d'après (5.3.1) la représentation cristalline $V^{(C)} = \text{Hom}_{(\phi,\text{Fil}^1)}(D^{(C)}, B_{\text{cris}}^+)$ que la théorie de Fontaine associe au module admissible $D^{(C)}$ est le module de Tate tensorisé par $\mathbf{Q}_p$ de ce groupe $p$-divisible. Mais $V$ est naturellement une sous-représentation de $V^{(C)}$, donc $V$ est aussi le module de Tate tensorisé par $\mathbf{Q}_p$ d'un groupe $p$-divisible sur $\mathcal{O}_K$ par [Ra, 2.1 et 2.3.1] ou [Ta, Prop. 12]. □

On a donc finalement, puisque le module de Dieudonné d'un groupe $p$-divisible est toujours faiblement admissible:

COROLLAIRE 5.3.3. *Supposons $p \neq 2$ et $k \subseteq \overline{\mathbf{F}}_p$. Le foncteur "module de Dieudonné" induit une (anti-)équivalence de catégories entre les groupes $p$-divisibles sur $\mathcal{O}_K$ à isogénie près et les $\phi$-modules filtrés faiblement admissibles $D$ tels que $\text{Fil}^0(D \otimes_{K_0} K) = D \otimes_{K_0} K$ et $\text{Fil}^2(D \otimes_{K_0} K) = 0$.*

Signalons enfin le corollaire:

COROLLAIRE 5.3.4. *Supposons $p \neq 2$ et $k \subseteq \overline{\mathbf{F}}_p$. Soient $A$ une variété abélienne définie sur $K$ et $V_p(A)$ son module de Tate tensorisé par $\mathbf{Q}_p$. Alors $A$ a réduction semi-stable (resp. bonne réduction) sur $\mathcal{O}_K$ si et seulement si $V_p(A)$ est une représentation semi-stable (resp. cristalline) de $\text{Gal}(\overline{K}/K)$.*

*Preuve.* Les implications "$A$ a réduction semi-stable (resp. bonne réduction) entraîne $V_p(A)$ est une représentation semi-stable (resp. cristalline)" sont classiques et bien connues (et dues à Fontaine). Soit $T_p(A)$ le module de Tate de $A$ et supposons $V_p(A)$ cristalline, alors par (5.3.2) et [Ra, 2.1, 2.3.1], $T_p(A)$ provient d'un groupe $p$-divisible sur $\mathcal{O}_K$ et par (SGA$_7$ exposé IX cor.5.10) $A$ a bonne réduction sur $\mathcal{O}_K$. Supposons $V_p(A)$ semi-stable non cristalline et soit $D$ le module admissible à pentes positives associé à $V_p(A)$. Comme $V_p(A)$ n'est pas cristalline, l'opérateur $N$ sur $D$ est par définition non nul et comme les pentes du Frobenius sur $D$ sont entre 0 et 1, on vérifie facilement que les pentes sur $N(D)$ sont toutes nulles et que $N^2(D) = 0$. Par la condition de faible admissibilité, $\text{Fil}^1(N(D) \otimes_{K_0} K) = 0$ et $N(D)$ est admissible, donc aussi $D/N(D)$. Soit $V' \subset V_p(A)$ la sous-représentation cristalline correspondant à $D/N(D)$, alors $V_p(A)/V'$ correspond à $N(D)$ et est aussi cristalline. Soit $T' = V' \cap T_p(A)$, alors par (5.3.2) et [Ra, 2.1, 2.3.1], $T'$ et $T_p(A)/T'$ sont les modules de Tate de groupes $p$-divisibles sur $\mathcal{O}_K$ et par (SGA$_7$ exposé IX prop.5.13.c) $A$ a réduction semi-stable sur $\mathcal{O}_K$. □

*Remarque* 5.3.5. Le corollaire ci-dessus était déjà connu (pour $k$ parfait et $p$ quelconque) dans le cas $e \leq p-1$ comme conséquence de [Laf1] et dans le cas de bonne réduction ([Mo] si $A$ est potentiellement un produit de jacobiennes, [CI] pour le cas général de bonne réduction).

C.N.R.S., Université Paris-Sud, Orsay, France
*Adresse e-mail*: breuil@math.u-psud.fr